\documentclass[onefignum,onetabnum]{siamart251216}
\usepackage{ifpdf}
\usepackage{booktabs}

\usepackage{amsmath,amssymb,amsfonts,latexsym,stmaryrd}
\usepackage{graphicx,float,epsfig} 
\usepackage{mathrsfs} 
\usepackage{color}
\usepackage{setspace} 
\usepackage{lipsum}
\usepackage{placeins} 
\usepackage{graphicx,float}
\usepackage{color}
\usepackage{epstopdf}
\usepackage{subcaption}
\usepackage{longtable}
\usepackage{multirow}
\usepackage{svg}
\usepackage{url}
\newcommand{\bd}[1]{\boldsymbol{#1}}
\usepackage{diagbox}
\ifpdf
\DeclareGraphicsExtensions{.eps,.pdf,.png,.jpg}
\else
\DeclareGraphicsExtensions{.eps}
\fi
\newcommand{\Om}{\O}

\def\ds{\displaystyle}

\def\O{\Omega}

\renewcommand\sp{\mathop{\mathrm{Sp}}\nolimits}

\newcommand{\jump}[1]{\llbracket #1 \rrbracket}
\newcommand{\mean}[1]{\{#1\}}

\usepackage{booktabs}
\usepackage{float}

\newcommand{\N}{\mathbb N}

\newcommand\bu{\boldsymbol{u}}
\newcommand\bv{\boldsymbol{v}}
\newcommand\bw{\boldsymbol{w}}

\newcommand{\bH}{\bd{\mathcal{H}}}

\newcommand\bF{\boldsymbol{f}}

\newcommand\bS{\boldsymbol{S}}
\newcommand\bT{\boldsymbol{T}}

\newcommand{\cF}{\mathcal{F}}

\def\CT{{\mathcal T}}

\newcommand\bsig{\boldsymbol{\sigma}}
\newcommand\btau{\boldsymbol{\tau}}

\newcommand\R{\mathbb{R}}

\renewcommand\H{\mathrm{H}}

\renewcommand\L{\mathrm{L}}

\renewcommand\O{\Omega}

\newcommand{\vertiii}[1]{{\left\vert\kern-0.25ex\left\vert\kern-0.25ex\left\vert #1 
		\right\vert\kern-0.25ex\right\vert\kern-0.25ex\right\vert}}
\newcommand\bdiv{\mathop{\mathbf{div}}\nolimits}

\newcommand\tr{\mathop{\mathrm{tr}}\nolimits}

\renewcommand\sp{\mathop{\mathrm{sp}}\nolimits}

\newcommand{\jumpt}[1]{\boldsymbol{[} #1 \boldsymbol{]}}

\newsiamremark{remark}{Remark}
\newsiamremark{hypothesis}{Hypothesis}
\crefname{hypothesis}{Hypothesis}{Hypotheses}
\newsiamthm{problem}{Problem}
\newsiamthm{claim}{Claim}

\headers{ DG method the Stokes-Brinkman equations }{Isabella Estrada and Felipe Lepe}

\title{A DG method for the  pseudostress-velocity formulation  of the  Stokes-Brinkman  eigenvalue problem\thanks{Submitted to the editors DATE.
		\funding{The first and second authors were  supported by Universidad del B\'io- B\'io through Proyecto Regular RE2514703. 
}}}

\author{Isabella Estrada\thanks{GIMNAP-Departamento de Matem\'atica, Universidad del B\'io - B\'io, Casilla 5-C, Concepci\'on, Chile. \email{iestrada@ubiobio.cl
}.}
\and Felipe Lepe\thanks{GIMNAP-Departamento de Matem\'atica, Universidad del B\'io - B\'io, Casilla 5-C, Concepci\'on, Chile. \email{flepe@ubiobio.cl}.}}
\usepackage{amsopn}

\ifpdf
\hypersetup{
  pdftitle={The Stokes-Brinkman eigenvalue problem},
  pdfauthor={Isabella Estrada, Felipe Lepe}
}
\fi

\def\CT{{\mathcal T}}

\begin{document}

\maketitle

\begin{abstract}
We introduce a discontinuous Galerkin method for the approximation of the eigenvalues and eigenfunctions of the Stokes-Brinkman eigenvalue problem. The method is based on a mixed formulation in which the pseudostress tensor and the fluid velocity are the main unknowns, while the pressure can be recovered through a post-processing procedure applied to the pseudostress tensor. Using discontinuous polynomial spaces of degree $k\geq 1$, we propose a consistent symmetric interior penalty discretization for the approximation of the spectrum. Relying on the spectral theory for non-compact operators, we establish convergence of the discrete spectral approximation and derive error estimates for the method. A series of numerical experiments is presented to assess its computational performance. On the one hand, we investigate the convergence rates for different polynomial degrees, computational domains, and values of the permeability. Furthermore, we study the influence of the stabilization parameter on the discrete spectrum by comparing different choices of this parameter and examining their impact on the occurrence of spurious eigenvalues.
\end{abstract}

\begin{keywords}
  porous media, fluid equations,  mixed formulations, Stokes-Brinkman equations, discontinuous Galerkin methods, a priori error estimates
\end{keywords}

\begin{AMS}
  35P15, 35Q35,  65N15, 65N25, 65N30, 65N50, 76M10, 76S99
\end{AMS}

\section{Introduction}
\label{sec:intro}
Let $\Omega\subset \R^d$, with $d\in\{2,3\}$, be an open and bounded domain with Lipschitz boundary $\partial \Omega$.  We are interested in the following eigenvalue problem (see \cite{MR4937821}): 
\begin{equation}
\label{eq:brinkman_stokes_PDE}
    -\nu\Delta \bd{u}+ \mathbf{K}^{-1}\bd{u}+\nabla p = \lambda \bu \,\, \text{in } \Omega,\quad
    \mathrm{div}\,\bd{u} = 0 \,\,  \text{in } \Omega, \quad
    \displaystyle\int_\Omega p = 0 \,\, \text{in } \Omega,\quad
    \bd{u} =0 \,\,\text{on } \partial\Omega,
\end{equation} 
where $\lambda\in \R$ is the eigenvalue, $\bd{u}$ represents the fluid velocity, $p$ is the pressure, $\nu>0$ represents the kinematic viscosity, and $\mathbf{K}^{-1}$ is a tensor associated with the permeability of the medium. This problem simultaneously addresses two different phenomena: on the one hand, we recover the classical Stokes eigenvalue problem when the permeability is sufficiently small, whereas, on the other hand, when the kinematic viscosity is small compared with the permeability, Darcy's law becomes the predominant phenomenon and describes the motion of a fluid through a porous medium. This is one of the most relevant features represented by system \eqref{eq:brinkman_stokes_PDE}. From the source-problem perspective, this problem has been addressed, for instance, in \cite{hannukainen2011computations,juntunen2010analysis, konno2011h}, where special care must be taken when the physical parameters are chosen. Indeed, when the viscosity tends to zero, the Darcy problem must be equipped with a boundary condition associated with the normal component of the velocity, whereas, in the Stokes case, the Dirichlet boundary condition is the appropriate one. Hence, different functional spaces are required according to the nature of the boundary conditions. In our case, for \eqref{eq:brinkman_stokes_PDE}, the case $\nu=0$ is not allowed, since, to the best of our knowledge, a proper physical interpretation of a pure Darcy eigenvalue problem is not evident.

System \eqref{eq:brinkman_stokes_PDE} has been recently studied by means of a finite element method for the velocity-pressure formulation in \cite{MR4937821}, where it was proved both theoretically and numerically that the method is capable of approximating the spectrum while avoiding spurious eigenvalues. This reference also illustrates the physical and geometrical implications of \eqref{eq:brinkman_stokes_PDE}. In particular, when the permeability is sufficiently large, the porous regions of the computational domain begin to behave as obstacles and, consequently, the regularity of the eigenfunctions is affected by the singularities induced by the highly permeable regions. This phenomenon, which does not arise in the classical Stokes eigenvalue problem, is of particular interest and must be accounted for in applications involving the stability of fluid flow in porous media. This is especially relevant  when considering more complex systems than \eqref{eq:brinkman_stokes_PDE}, such as coupled fluid–structure systems.

Regarding the literature, mixed formulations for flow eigenvalue problems have been extensively studied. In particular, for \eqref{eq:brinkman_stokes_PDE}, mixed finite element methods have been introduced in \cite{MR5023159,MR5105993}, where formulations involving either the pseudostress tensor or the vorticity have been analyzed. These articles require a specific analytical approach, since the permeability term in \eqref{eq:brinkman_stokes_PDE} leads to a perturbed saddle-point problem that must be treated accordingly. In this setting, the corresponding source problem plays a fundamental role in establishing the well-posedness of the associated solution operators. Our aim is to further develop our research program on the eigenvalue problem \eqref{eq:brinkman_stokes_PDE} by investigating an alternative numerical approach based on the discontinuous Galerkin (DG) method.

The DG method is a well-established numerical approach based on piecewise polynomial approximations, which has been successfully applied to the numerical solution of a wide range of problems; see, for instance, \cite{MR2882148}. In particular, its application to eigenvalue problems is well documented in, for example, \cite{MR2220929,MR2324460,MR2594995,MR3962898,MR4077220,MR4407474,MR4564097,MR4090799}, where the method has shown good accuracy in the approximation of the spectrum of linear selfadjoint and non-selfadjoint operators under different formulations. For the problem \eqref{eq:brinkman_stokes_PDE}, a particular class of DG methods, known as interior penalty discontinuous Galerkin (IPDG) methods, has been introduced for the velocity-pressure formulation in two and three dimensions. The authors demonstrate that IPDG methods provide a suitable numerical approach for this formulation of \eqref{eq:brinkman_stokes_PDE}. However, their application to mixed formulations of the problem, such as the one proposed in \cite{MR5023159}, remains open.

Our work aims to apply the IPDG method to the mixed formulation introduced in \cite{MR5023159}. The source version of this problem has been studied in \cite{MR4108273} for the Brinkman source problem. That paper considers a symmetric IPDG method, since this version is precisely the consistent one for the problem. Moreover, the pseudostress employed there is symmetric, since it is defined in terms of the symmetric part of the velocity gradient. Our proposal is instead to use the non-symmetric pseudostress tensor and to analyze the IPDG method with this unknown. An immediate difference appears in the polynomial space used to approximate the pseudostress. In \cite{MR4108273}, the authors seek the tensor in a piecewise polynomial space in which symmetry is strongly imposed, whereas in our case no such restriction is required, since we use tensor-valued polynomials of degree $k$ without additional symmetry constraints. Another key aspect that changes is the vectorial jump, since in our case this definition does not need to be symmetrized.

A further relevant aspect is the symmetry of the method. IPDG methods can also be formulated in non-symmetric versions, such as the non-symmetric interior penalty (NIP) and the incomplete interior penalty (IIP) method. These methods have been studied for eigenvalue problems, for instance, in \cite{MR2220929,MR2324460,MR4077220}, where significant differences among them are observed, particularly with respect to the behavior of spurious eigenvalues. Since consistency of the mixed IPDG formulation considered in this work is achieved only by the SIP method, we restrict our attention to this variant and do not consider the NIP and IIP methods. An additional advantage of the symmetric formulation is its ability to achive optimal orders of convergence for the eigenvalues, which NIP and IIP methods are generally unable to achieve; see, for instance, the numerical experiments reported in \cite{MR2324460,MR4077220}. Nevertheless, a detailed analysis of the influence of the penalization parameter in the IPDG method required, following the approach of \cite{MR4177014}, since this parameter directly affects the matrices of the generalized eigenvalue problem. Consequently, the convergence rates may also be affected and, more importantly, spurious eigenvalues may emerge if the stabilization parameter is not chosen appropriately.
\subsection{Outline}
The organization of the paper is as follows: In Section \ref{sec:variational}, we introduce the mixed formulation of \eqref{eq:brinkman_stokes_PDE} and summarize several results of the continuous formulation, such as regularity, well-definition of the solution operators, and spectral characterization. The core of the paper begins in Section \ref{sec:dg}, where we introduce the DG method. Here, we define mesh elements, norms and seminorms related to the method, discrete bilinear forms and, hence, the DG discretization of the continuous problem. Section \ref{sec:convegence} is devoted to the convergence analysis of the proposed discretization. In particular, we prove the well-posedness of the discrete solution operators, convergence, and a priori error estimates. Finally, in Section \ref{sec:numerics} we report a series of numerical experiments in two and three dimensions where we study the influence of the stabilization parameter on the computation of the spectrum and orders of convergence.

\section{Variational formulation}
\label{sec:variational}
From now and on, we will assume that $\textbf{K}^{-1}$ is a symmetric, bounded, and  positive definite matrix,  i.e., there exists two constants  $\kappa_*,\kappa^*>0$ such that  
\begin{equation*}
	\label{eq:kappa-bounds}
	0<\kappa_*\Vert\bv\Vert_{0,\Omega}^2 \leq (\textbf{K}^{-1}\bv,\bv)_{0,\Omega}\leq \kappa^*\Vert{\bv}\Vert_{0,\Omega}^2.
\end{equation*}
%
%
Following the analysis of  \cite{MR4937821},  let us introduce the pseudostress tensor $\bd{\sigma}$, defined by   $\bd{\sigma} := \nu\nabla \bd{u}-p\mathbf{I}$ where  $\mathbf{I}\in\mathbb{R}^{d\times d}$ represents the identity   tensor. Then, applying $\bdiv$ on   the momentum equation in  \eqref{eq:brinkman_stokes_PDE}, we obtain  $\bdiv\bsig=\mathbf{K}^{-1}\bu-\lambda\bu$ in $\Omega$. Moreover, we define the deviator of   $\bsig$ by 
$
\bsig^{\texttt{d}}:=\bsig-\frac{1}{d}\tr(\bsig)\mathbf{I}.
$
 
 Using the incompressibility condition of \eqref{eq:brinkman_stokes_PDE}, we have  $\bd{\sigma}^\texttt{d}=\nu\nabla\bd{u}$ in $\Omega$. Hence, problem \eqref{eq:brinkman_stokes_PDE}  is rewritten as follows: find a non-vanishing pair $(\lambda,(\bd{\sigma},\bd{u}))$, and $\lambda\in\mathbb{R}^+$ such that
\begin{equation}
\label{eq:eigenproblem_bs}
\left\{
\begin{aligned}
\bdiv\bsig-\mathbf{K}^{-1}\bu
&= -\lambda\bu
&& \text{in } \Omega,\\
\bsig^{\texttt{d}}-\nu\nabla\bu
&= \boldsymbol{0}
&& \text{in } \Omega,\\
\bu
&= \boldsymbol{0}
&& \text{on } \partial\Omega,\\
\int_{\Om}\tr(\bsig)
&=0.
\end{aligned}
\right.
\end{equation}
We remark that the pressure can be recovered via a post-process of the pseudostress tensor using the relation $p:=d^{-1}\tr(\bsig)$. 

The variational formulation of \eqref{eq:eigenproblem_bs} is the following:  find  $\lambda\in\R^+$ and   $(\mathbf{0},\mathbf{0})\neq (\bd{\sigma},\bd{u})\in \bd{\mathcal{H}}\times \mathbf{Q}$ such that 
\begin{align*}
\frac{1}{\nu}\int_{\Omega}\bsig^{\texttt{d}}:\btau^{\texttt{d}}+\int_{\Omega}\bu\cdot\bdiv\btau&=0\quad\quad\quad\quad\,\,\,\forall\btau\in\bd{\mathcal{H}},\\
\int_{\Omega}\bdiv\bsig\cdot\bv-\int_{\Omega}\mathbf{K}^{-1}\bu\cdot\bv&=-\lambda\int_{\Omega}\bu\cdot\bv\quad\forall\bv\in\mathbf{Q},
\end{align*}
where  $\bd{\mathcal{H}}:=\bd{\mathcal{H}}(\bdiv,\Omega)$ and $\mathbf{Q}:=\L^2(\Omega)^d$.  With these spaces at hand, we define the bounded bilinear forms $a:\bd{\mathcal{H}}\times \bd{\mathcal{H}}\to\R$, $b:\bd{\mathcal{H}}\times \mathbf{Q}\to\R$ and  $c:\mathbf{Q}\times\mathbf{Q}\rightarrow\mathbb{R}$ as follows
\begin{align*}
a(\boldsymbol{\rho},\btau)&:=\dfrac{1}{\nu}\int_\Omega\boldsymbol{\rho}^{\texttt{d}}:\btau^{\texttt{d}}\qquad &\boldsymbol{\rho},\btau\in\bd{\mathcal{H}},\\
b(\btau,\bv)&:=\int_\Omega\bdiv\btau\cdot\bv\qquad &(\btau,\bv)\in \bd{\mathcal{H}}\times\mathbf{Q},\\
c(\bv,\bw)&:=\int_\Omega\mathbf{K}^{-1}\bv\cdot\bw\qquad &\bv,\bw\in \mathbf{Q}.
\end{align*}
Hence, the variational  problem can be written as follows:  find $\lambda\in\R^+$ and  $(\boldsymbol{0},\boldsymbol{0})\neq (\bd\sigma,\bd{u})\in \bd{\mathcal{H}}\times \mathbf{Q}$ such that
\begin{equation}
\label{eq:brinkman_fv_1_not}
	\left\{
	\begin{array}{rcll}
a(\bsig,\btau)+ b(\btau,\bu)&=&0&\forall \btau\in \bd{\mathcal{H}},\\

 b(\bsig, \bv)-c(\bv,\bu)&=&-\lambda(\bu,\bv)_{0,\Omega}&\forall \bv\in \mathbf{Q}.
\end{array}
	\right.
	\end{equation}
	
	Let us consider the following  decomposition  $\bd{\mathcal{H}}=\bd{\mathcal{H}}_0\oplus \R\mathbf{I}$, where the space $\bd{\mathcal{H}}_0$ is defined by
$$\bd{\mathcal{H}}_0:=\left\{\btau\in \bd{\mathcal{H}}: \displaystyle\int_\Omega \tr(\bd\tau)=0 \quad\mathrm{in}\quad\Omega\right\}.$$ 
We observe that this space accounts for the zero-mean constraint on the pressure in \eqref{eq:brinkman_stokes_PDE}, which is expressed in terms of the pseudostress by the third equation of \eqref{eq:eigenproblem_bs}.

With this decomposition at hand, for any  $\bd{\tau}\in \bd{\mathcal{H}}$, there exist $\btau_0\in\bd{\mathcal{H}}_0$ and a constant $\overline{c}\in\mathbb{R}$ such that
\begin{equation*}
\label{eq:def_m}
\btau=\btau_0+\overline{c}\mathbf{I},\quad\text{with}\quad \btau_0\in \bd{\mathcal{H}}_0\quad \text{and}\quad \overline{c}:=\dfrac{1}{d|\Omega|}\int_{\Omega}\tr(\btau)\in\mathbb{R},
\end{equation*}
and also, problem \eqref{eq:brinkman_fv_1_not} is equivalently written as follows: find  $\lambda\in\R^+$,  $(\boldsymbol{0},\boldsymbol{0})\neq(\bsig_0,\bd{u})\in \bd{\mathcal{H}}_0\times \mathbf{Q}$, such that
\begin{equation}
\label{eq:brinkman_fv_1_reduced}
	\left\{
	\begin{array}{rcll}
a(\bsig_0,\btau)+ b(\btau,\bu)&=&0&\forall \btau\in \bd{\mathcal{H}}_0,\\

 b(\bsig_0, \bv)-c(\bv,\bu)&=&-\lambda(\bu,\bv)_{0,\Omega}&\forall \bv\in \mathbf{Q}.
\end{array}
	\right.
	\end{equation}

For the analysis of the eigenvalue problem, we introduce the so-called solution operators, which we define by
\begin{equation*}
\bT:\mathbf{Q}\rightarrow\mathbf{Q},\,\,\,\,\,\bF\mapsto\bT\bF:=\widehat{\bu}\,\,\,\text{and}\,\,\,\,\,\bS:\mathbf{Q}\rightarrow\bd{\mathcal{H}}_0\,,\,\,\bF\mapsto \bS\bF:=\widehat{\bsig}_0,
\end{equation*}
where the pair  $(\widehat{\bsig}_0,\widehat{\bu})\in\bd{\mathcal{H}}_0\times\mathbf{Q}$ is the solution of the following source problem
\begin{equation}
\label{eq:source_continuous_red}
	\left\{
	\begin{array}{rcll}
a(\widehat{\bsig}_0,\btau)+ b(\btau,\widehat{\bu})&=&0&\forall \btau\in \bd{\mathcal{H}}_0,\\

 b(\widehat{\bsig}_0, \bv)-c(\bv,\widehat{\bu})&=&-(\bF,\bv)_{0,\Omega}&\forall \bv\in \mathbf{Q}.
\end{array}
	\right.
\end{equation}
We now state the following technical result, which will be used to establish the coercivity of the bilinear form $a(\cdot,\cdot)$; see, for instance, \cite[Chapter 9, Proposition 9.1.1]{MR3097958}.

\begin{lemma}
\label{lmm:deviator}
There exists $C_1>0$, depending on $\Omega$, such that  
\begin{equation*}
	\label{eq:des_deviator}
	C_1\|\btau\|_{0,\Omega}^2\leq \|\btau^{\texttt{d}}\|_{0,\Omega}^2+\|\bdiv\btau\|_{0,\Omega}^2\quad\forall\btau\in\bd{\mathcal{H}}_0.
\end{equation*}
\end{lemma}

Let us define the kernel of $b(\cdot,\cdot)$ as follows
\begin{equation*}
	\boldsymbol{\mathcal{V}}:=\{\btau\in \bd{\mathcal{H}}_0\,:\, b(\btau,\bv)=0\,\,\,\forall \bv\in \mathbf{Q}\}=\{\btau\in \bd{\mathcal{H}}_0\,:\,\, \bdiv\btau=\boldsymbol{0}\}.
\end{equation*}

Now, if  $\bd\tau\in\boldsymbol{\mathcal{V}}$, then for $a(\cdot,\cdot)$ we have
\begin{equation}
\label{eq:a_coercive}
a(\btau,\btau)=\frac{1}{\nu}\|\btau^{\texttt{d}}\|_{0,\Omega}^2\geq \frac{C_1}{\nu}\|\btau\|_{0,\Omega}^2=\frac{C_1}{\nu}\|\btau\|_{\bdiv,\Omega}^2,
\end{equation}
where $C_1>0$ is the constant involved in Lemma  \ref{lmm:deviator}, proving that  $a(\cdot,\cdot)$ is  $\boldsymbol{\mathcal{V}}$-coercive  with constant  $\alpha=C_1/\nu$.  

Additionally, by applying \cite[Theorem 2.1]{MR2594823}, we obtain the existence of a constant $\beta>0$ such that the following inf-sup condition holds
\begin{equation}
\label{eq:inf_sup_b}
	\displaystyle\sup_{\boldsymbol{0}\neq\btau\in \bd{\mathcal{H}}_{0}}\frac{b(\btau,\bv)}{\|\btau\|_{\bdiv,\Omega}}\geq\beta\|\bv\|_{0,\Omega}\quad\forall\bv\in \mathbf{Q}.
\end{equation}

On the other hand, for  $\bv\in\mathbf{Q}$, we have
\begin{equation}
\label{eq:c_coercive}
c(\bv,\bv)=\int_{\Omega}\mathbf{K}^{-1}\bv\cdot\bv\geq \kappa_*\|\bv\|_{0,\Omega}^2\quad\forall\bv\in\mathbf{Q}.
\end{equation}

Then, from   \eqref{eq:a_coercive}, \eqref{eq:inf_sup_b} and \eqref{eq:c_coercive}, due to the Babu\v ska-Brezzi theory   (see for instance  \cite[Proposition 4.3.1]{MR3097958}), we deduce that problem \eqref{eq:source_continuous_red}  has a unique solution  $(\widehat{\bd\sigma}_0,\widehat{\bd{u}})\in \bd{\mathcal{H}}_0\times \mathbf{Q}$, implying that  $\bT$ and $\bS$ are well defined. 	Moreover, there exists $C>0$ such that the following estimate holds
\begin{equation*}
\|\widehat{\bsig}_0\|_{\bdiv,\Omega}+\|\widehat{\bu}\|_{0,\Omega}\leq C\|\bF\|_{0,\Omega}.
\end{equation*}

Since problem \eqref{eq:brinkman_fv_1_reduced} is symmetric, an immediate consequence is that  $\bT$ is selfadjoint with respect to the $\L^2(\O)^d$ inner product.
Moreover,  $(\lambda, (\bsig_0,\bu))\in\mathbb{R}^+\times (\bd{\mathcal{H}}_0\times \mathbf{Q})$  is solution of problem \eqref{eq:brinkman_fv_1_reduced} if and only if  $(1/\lambda,\bu)$ is an eigenpair of  $\bT$, on the sense that $\bT\bu=\ds\frac{1}{\lambda}\bu$.

Let us focus now on the regularity properties of the source solution and the eigenfunctions of the eigenvalue problem. Since the term  $(\mathbf{K}^{-1}\bu,\bv)_{0,\Omega}$ is well defined due to the hypothesis on the permeability tensor $\mathbf{K}^{-1}$, we take advantage of the well known results for the regularity of the Stokes problem provided by, for instance, in \cite{MR975121,MR1600081} to establish the following result.
\begin{theorem}
\label{th:reg_velocity}
If $(\lambda, (\widehat{\bu},\widehat{p}))\in \mathbb{R}^+\times\H_0^1(\Omega)^d\times \L_0^2(\Omega)$  solves \eqref{eq:source_continuous_red}, then there exists  $s\in(0,1]$ such that  $\widehat{\bu}\in \H^{1+s}(\Omega)^d$ and $\widehat{p}\in \H^s(\Omega)$.
\end{theorem}

It is worth to noticing that if $\Omega$ is convex, then $s=1$. Consequently, the velocity belongs to $\H^2(\Omega)^d$, while the pseudostress and its divergence belong to $\H^1(\O)^{d\times d}$ and $\H^2(\Omega)^d$, respectively. Moreover, the following estimate holds 
\begin{equation}
	\label{eq:reg_u_if}
\|\widehat{\bsig}_0\|_{s,\Omega}+\|\bdiv\widehat{\bsig}_0\|_{1+s,\Omega}+\|\widehat{\bu}\|_{1+s,\Omega}\leq C\|\boldsymbol{f}\|_{0,\Omega},
\end{equation}
where $C>0$ depends on $\O$ and the physical parameters of the problem.

Now, for the eigenfunctions we have the following regularity. If the pair $(\bsig_0,\bu)\in\bd{\mathcal{H}}_0\times\mathbf{Q}$ is  solution of   \eqref{eq:brinkman_fv_1_reduced}, there exists a constant $C_{\lambda}>0$,  depending on $\lambda$ and a  parameter  $r>0$ such that 
\begin{equation}
	\label{eq:reg_u_if2}
	\|\bsig_0\|_{r,\Omega}+\|\bdiv\bsig_0\|_{1+r,\Omega}+\|\bu\|_{1+r,\Omega}\leq C_{\lambda}\|\bu\|_{0,\Omega}.
\end{equation}

Invoking Theorem \ref{th:reg_velocity} and  \eqref{eq:reg_u_if} we conclude that $\bT$ is compact. Moreover, since $\bT$ is selfadjoint, the spectrum of $\bT$ which we denote by $\sp(\bT)$, is given by  $\text{sp}(\bT):=\{\mu_k\}_{k\in\N}\cup\{0\}$, where $\{\mu_k\}_{k\in\N}$ is a sequence of eigenvalues of  $\bT$ that converge to zero as  $k\to\infty$.
\subsection{The discontinuous Galerkin method}
\label{sec:dg}
In the following section, we introduce the discontinuous Galerkin method that we will use to approximate the spectrum and eigenfunctions of the model problem. We begin with some preliminary definitions. We consider a regular triangulation $\{\mathcal{T}_h\}_h$ of $\overline{\Omega}$, which may include hanging nodes. Let $h$ be the mesh parameter defined by $h=\max_{T\in\mathcal{T}_h}h_T$, where $h_T$ denotes the diameter of each element $T\in\mathcal{T}_h$. Let $T^+,T^-\in\mathcal{T}_h$ be two adjacent elements (not necessarily matching). We define an interior face of $\mathcal{T}_h$ as the nonempty interior of $\partial T^+\cap \partial T^-$. On the other hand, a boundary face is defined as the nonempty interior of $\partial T\cap \partial\Omega$, where $T$ is a boundary element of $\mathcal{T}_h$. We denote by $\mathcal{F}_h^i$ the set of all interior faces of $\mathcal{T}_h$, by $\mathcal{F}_h^\partial$ the set of all boundary faces, and define $\mathcal{F}_h=\mathcal{F}_h^i\cup\mathcal{F}_h^\partial$.

In what follows, we denote by  $\bdiv_h$  the broken divergence operator that acts on the rows  of tensorial fields. When we restrict this operator to a particular element $T\in\CT_h$ we only write $\bdiv$.

Let  $\bv$ and  $\btau$ sufficiently smooth vectorial and tensorial functions, respectively. For each interior edge/face $F\in\mathcal{F}_h^i$ which is shared by two elements  $T^+$ and $T^-$, we denote by  $\mathbf{n}^+$ and $\mathbf{n}^-$  the normal unit vectors of $T^+$ and $T^-$, respectively.  Also we define the vectorial and tensorial traces by   $\bv^\pm=\bv|_{\partial T^\pm}$ and  $\btau^\pm=\btau|_{\partial T^\pm}$, respectively.

On the other hand,  we define the averages with respect to a edge/face $F$ for vectorial and tensorial fields are defined, respectively, as follows
$$\mean{\bv}:=\frac{1}{2}(\bv^++\bv^-),\quad \mean{\btau}=\frac12(\btau^++\btau^-),$$ 
whereas the interelement jumps for vectorial fields are defined by    $$\jump{\bv}:=\bv^+ \otimes \mathbf{n}^++\bv^-\otimes \mathbf{n}^- ,$$
where $\bv\otimes\bu$ is a tensor whose $(i,j)-$th entry is given by  $v_i u_j$ and  for  tensorial fields, the jump  is defined by 
$\jumpt{\btau}=\btau^+\mathbf{n}^++\btau^-\mathbf{n}^-$.

Now, if $F\in\mathcal{F}_h^\partial$ is a boundary edge/face, the previous definitions are given by  $\mean{\bv}:=\bv$ and $\mean{\btau}=\btau$ for the averages of vector and tensor, respectively,  whereas  the jumps  for vectors and tensors are given respectively  by  $\jump{\bv}:=\bv\otimes\mathbf{n}$ and $\jumpt{\btau}=\btau\mathbf{n}$.

For any vector $\bv$ and tensor  $\btau$, the following identity holds  (see   \cite{MR4108273})
\begin{equation*}\label{id_sum}
    \sum_{T\in\mathcal{T}_h}\int_{\partial T} (\btau\mathbf{n}_T)\cdot\bv=\int_{\mathcal{F}_h^i}\jumpt{\btau}\cdot\mean{\bv}+\int_{\mathcal{F}_h}\mean{\btau}:\jump{\bv}.
\end{equation*}

For a fixed integer  $k\geq 1$,we introduce the following discrete spaces  $\bd{\mathcal{H}}_h$  and  $\mathbf{Q}_h$  defined by 
$$\bd{\mathcal{H}}_h:=\left\{\btau_h\in \L^2(\Om)^{d\times d}:\btau_h|_T\in {\mathcal{P}}_{k}(T)^{d\times d} \quad\forall T\in\mathcal{T}_h\right\},$$
and
 $$\mathbf{Q}_h:=\{\bv_h\in\mathbf{Q}:\bv_h\in \mathcal{P}_{k-1}(T)^d\quad \forall T\in\mathcal{T}_h\},$$
where the pseudostress tensor is approximated with the elements of $\bd{\mathcal{H}}_h$ and the velocity by $\mathbf{Q}_h$. Moreover, we endow the product space  $\bH_h\times \mathbf{Q}_h$ with the following norm
$$\vertiii{(\btau_h,\bv_h)}^2:=\|\btau_h\|_{\bH_h}^2+\|\bv_h\|_{0,\Omega}^2 \quad\forall \btau_h\in\bH_h,\forall\bv_h\in\mathbf{Q}_h,$$
where  $\|\cdot\|_{\bd{\mathcal{H}}_h}$ is defined by  $$\|\btau_h\|^2_{\bd{\mathcal{H}}_h}=\| \btau_h \|_{0,\Om}^2+\|\bdiv_h\btau_h\|^2_{0,\Om}+\|h_{\mathcal{F}}^{-1/2}\jumpt{\btau_h}\|^2_{0,\mathcal{F}_h^i}.$$

Finally, the discrete version of $\bd{\mathcal{H}}_{0}$ is defined by 
$$\bd{\mathcal{H}}_{0,h}:=\left\{\btau_h\in\bd{\mathcal{H}}_h\,:\, \int_\Om\tr{(\btau_h)}=0\right\}.$$

For  our analysis,   we introduce the following trace inequality  (see \cite[Section 3]{MR4673997} and  \cite[Proposition 4.1]{MR3376135})
\begin{equation}\label{eq:dg_trace}
\|h_{\mathcal{F}}^{1/2}\mean{\xi(x)\bv_h}\|_{0,\cF_h}\leq C_{\texttt{tr}}\|\bv_h\|_{0,\O}\qquad\forall \bv_h\in\bd{\mathcal{H}}_h, 
\end{equation}
where $\xi(x)\in \L^{\infty}(\O)$ and  $C_{\texttt{tr}}>0$ is a constant depending on $\|\xi(x)\|_{0,\infty}$, the  polynomial degree, the  spatial dimension, and the shape regularity of the mesh. In particular for our work, we need to consider $\xi(x)\equiv 1$ in \eqref{eq:dg_trace}.

Let us introduce the IPDG discretization of problem \eqref{eq:brinkman_fv_1_not}: find  $\lambda_h\in\R^+$, and   $(\boldsymbol{0},\boldsymbol{0})\neq(\bsig_h,\bu_h)\in \bd{\mathcal{H}}_{0,h}\times \mathbf{Q}_h$, such that
\begin{equation}
\label{eq:brinkman_dg}
	\left\{
	\begin{array}{rcll}
a_h(\bsig_h,\btau_h)+ b_h(\btau_h,\bu_h)&=&0&\forall \btau_h\in \bd{\mathcal{H}}_{0,h},\\

 b_h(\bsig_h, \bv_h)-c(\bv_h,\bu_h)&=&-\lambda_h(\bu_h,\bv_h)_{0,\Om}&\forall \bv_h\in \mathbf{Q}_h,
\end{array}
	\right.
\end{equation}
where the discrete bilinear forms  $a_h:\bd{\mathcal{H}}_{0,h}\times \bd{\mathcal{H}}_{0,h}\to\R$ y $b_h:\bd{\mathcal{H}}_{0,h}\times \mathbf{Q}_h\to\R$ are defined by 
\begin{equation*}
    a_h(\bsig_h,\btau_h)= \frac{1}{\nu}\int_\Om \bsig_h^\texttt{d}:\btau_h^\texttt{d} + \int_{\mathcal{F}_h^i} \frac{\texttt{a}_S}{h_{\mathcal{F}}} \jumpt{\bsig_h}\cdot\jumpt{\btau_h} ,
    \end{equation*}
    and
    \begin{equation*}
    b_h(\btau_h,\bv_h) = \int_\Om \bdiv_h\btau_h\cdot\bv_h- \int_{\mathcal{F}_h^i}\jumpt{\btau_h}\cdot \mean{\bv_h}.
\end{equation*}

On the other hand $c:\mathbf{Q}_h\times\mathbf{Q}_h\rightarrow\mathbb{R}$ is defined by $c(\bv,\bw):=\displaystyle\int_{\O}\mathbf{K}^{-1}\bv\cdot\bw$. Let us remark that in $a_{h}(\cdot,\cdot)$ the parameter $\texttt{a}_S>0$ is the so-called \emph{penalization parameter}, which is essential on the IDPG discretization. We refer to \cite{MR2882148} for a further discussion on the construction of IPDG methods and the need of this parameter.

It is easy to check that forms $a_h(\cdot,\cdot)$, $b_h(\cdot,\cdot)$ and $c(\cdot,\cdot)$ are bounded.  Indeed, let  $\bsig_h,\btau_h \in \bd{\mathcal{H}}_{0,h}$ and  $\bv_h\in\mathbf{Q}_h$. For $a_h(\cdot,\cdot)$ we have
\begin{multline*}
    |a_h(\bsig_h,\btau_h)| = \frac{1}{\nu}\|\bsig_h^\texttt{d}\|_{0,\Om}\|\btau_h^\texttt{d}\|_{0,\Om}+\texttt{a}_S\|h_{\cF}^{-1/2}[\bsig_h]\|_{0,\mathcal{F}_h^i}\|h_F^{-1/2}[\btau_h]\|_{0,\mathcal{F}_h^i}\\
    \leq \frac{1}{\nu}\left(1+\frac{\sqrt{d}}{d}\right)^2\|\bsig_h\|_{0,\Om}\|\btau_h\|_{0,\Om}+\texttt{a}_S\|h_{\cF}^{-1/2}[\bsig_h]\|_{0,\mathcal{F}_h^i}\|h_{\cF}^{-1/2}[\btau_h]\|_{0,\mathcal{F}_h^i}\\
    + \|\bdiv_h\bsig_h\|_{0,\Om}\|\bdiv_h\btau_h\|_{0,\Om}\\
    \leq \max\left\{\frac{1}{\nu}\left(1+\frac{\sqrt{d}}{d}\right)^2,\texttt{a}_S,1\right\}\left(\|\bsig_h\|_{0,\Om}^2+\|h_{\cF}^{-1/2}[\bsig_h]\|_{0,\mathcal{F}_h^i}+\|\bdiv_h\bsig_h\|^2_{0,\Om}\right)^{1/2}\\
    \times\left(\|\btau_h\|_{0,\Om}^2+\|h_{\cF}^{-1/2}[\btau_h]\|_{0,\mathcal{F}_h^i}+\|\bdiv_h\btau_h\|^2_{0,\Om}\right)^{1/2}
    = C_{\nu,\texttt{a}_S,d} \|\bsig_h\|_{\bd{\mathcal{H}}_h}\|\btau_h\|_{\bd{\mathcal{H}}_h}.
\end{multline*}

For  $b_h(\cdot,\cdot)$, we obtain
\begin{multline*}
    |b_h(\btau_h,\bv_h)| = \left|\int_\Om \bdiv_h\btau_h\cdot\bv_h- \int_{\mathcal{F}_h^i}\jumpt{\btau_h}\cdot \mean{\bv_h}\right|\\
    \leq \|\bdiv_h\btau_h\|_{0,\Om}\|\bv_h\|_{0,\Om} + \|h_\mathcal{F}^{-1/2}\jumpt{\btau_h}\|_{0,\mathcal{F}_h^i} \|h_\mathcal{F}^{1/2}\mean{\bv_h}\|_{0,\mathcal{F}_h^i} \\
    \leq  \left(\|\bdiv_h\btau_h\|_{0,\Om}^2+\|h_\mathcal{F}^{-1/2}\jumpt{\btau_h}\|_{0,\mathcal{F}_h^i}^2\right)^{1/2}\left(\|\bv_h\|_{0,\Om}^2+\|h_\mathcal{F}^{1/2}\mean{\bv_h}\|_{0,\mathcal{F}_h^i}^2\right)^{1/2}\\
    \leq  \left(\|\bdiv_h\btau_h\|_{0,\Om}^2+\|h_\mathcal{F}^{-1/2}\jumpt{\btau_h}\|_{0,\mathcal{F}_h^i}^2\right)^{1/2}\left(\|\bv_h\|_{0,\Om}^2+C_{\texttt{tr}}^2\|\bv_h\|_{0,\Om}^2\right)^{1/2}\\
    \leq  \left(\|\btau_h\|^2_{0,\Om}+ \|\bdiv_h\btau_h\|_{0,\Om}^2+\|h_\mathcal{F}^{-1/2}\jumpt{\btau_h}\|_{0,\mathcal{F}_h^i}^2\right)^{1/2}\left((1+C_{\texttt{tr}}^2)\|\bv_h\|_{0,\Om}^2\right)^{1/2}\\
    = (1+C_{\texttt{tr}}^2)^{1/2}\|\btau_h\|_{\bd{\mathcal{H}}_h}\|\bv_h\|_{0,\Om},
\end{multline*}
where $C_{\texttt{tr}}>0$ is the constant given by \eqref{eq:dg_trace}. For $c(\cdot,\cdot)$, the boundedness is direct because the treatment for this term  is analogous to the conforming finite element approach, which immediately implies that $|c(\bv,\bw)|\leq \kappa^*\|\bv\|_{0,\Omega}\|\bw\|_{0,\Omega}$.

Additionally, for $\bsig_h,\btau_h\in\bd{\mathcal{H}}_{0,h}$ and  $\bu_h,\bv_h\in\mathbf{Q}_h$, we define the symmetric and bounded bilinear form $A_h(\cdot,\cdot)$ by 
$$A_h((\bsig_h,\bu_h),(\btau_h,\bv_h))=a_h(\bsig_h,\btau_h)+b_h(\btau_h,\bu_h)+b_h(\bsig_h,\bv_h)-c(\bu_h,\bv_h).$$ 
Hence, problem \eqref{eq:brinkman_dg} is written as follows: find  $\lambda_h\in\mathbb{R}^+$  and  $(\boldsymbol{0},\boldsymbol{0})\neq(\bsig_h,\bu_h)\in\bd{\mathcal{H}}_{0,h}\times\mathbf{Q}_h$ such that
\begin{equation} \label{eq:formA_M_star}
    A_h((\bsig_h,\bu_h),(\btau_h,\bv_h))=-\lambda_h(\bu_h,\bv_h)_{0,\Omega}\qquad\forall(\btau_h,\bv_h)\in\bd{\mathcal{H}}_{0,h}\times\mathbf{Q}_h.
\end{equation}

\begin{remark}
Let us define the space $\boldsymbol{\mathcal{H}}(h):=\boldsymbol{\mathcal{H}}+\boldsymbol{\mathcal{H}}_{0,h}$. If we suppose that the pairs $(\bsig,\bu),(\btau,\bv)\in\boldsymbol{\mathcal{H}}(h)\times\mathbf{Q}$ satisfies   $\bdiv\btau,\bdiv\bsig\in\H^t(\O)^d$  with $t>1/2$, then the bilinear form $A_h(\cdot,\cdot)$ is still bounded in the sense that there exists a 
constant $M^*>0$, independent of $h$, such that
$$A_h((\bsig,\bu),(\btau,\bv))\leq M^*\vertiii{(\bsig,\bu)}_*\vertiii{(\btau,\bv)}_*,$$
where $\vertiii{(\btau,\bv)}_*^2:=\vertiii{(\btau,\bv)}^2+\|h_{\cF}^{1/2}\mean{\bdiv\btau}\|_{0,\cF_h^i}^2.$
\end{remark}

Let us introduce the   discrete solution operators, which are defined by 
$$\bT_h: \mathbf{Q} \to \mathbf{Q}_h,\quad \bF \mapsto \bT_h\bF:=\widehat{\bu}_h\quad\mathrm{and}\quad \bS_h:\mathbf{Q}\to \bd{\mathcal{H}}_{0,h},\quad \bF\mapsto \bS_h\bF:=\widehat{\bsig}_h,$$ where  $(\widehat{\bsig}_h,\widehat{\bu}_h)\in \bd{\mathcal{H}}_{0,h}\times \mathbf{Q}_h$ is the solution of the following discrete source problem
\begin{equation}\label{eq:brinkman_dgf}
    A_h((\widehat\bsig_h,\widehat\bu_h),(\btau_h,\bv_h))=-(\bF,\bv_h)_{0,\Om} \quad \forall (\btau_h,\bv_h)\in\bd{\mathcal{H}}_{0,h}\times \mathbf{Q}_h.
\end{equation}

Now the task is to prove that $\bT_h$ and  $\bS_h$ are well defined. Hence, we need to conclude that the discrete source problem is well-posed. To do this, we need the following technical result (see \cite[Lemma 3.3]{MR3417280}).
\begin{lemma}\label{lema:normL2}
    There exists a constant $C>0$, independent of  $h$, such that 
    \begin{equation*}
        \|\btau_h\|_{0,\Om}^2\leq C_1(\|\btau_h^{\texttt{d}}\|_{0,\Om}+\|\bdiv_h\btau_h\|_{0,\Om}^2+\|h_\mathcal{F}^{-1/2}\jumpt{\btau_h}\|^2_{0,\mathcal{F}_h^i})\quad \forall\btau_h\in \bd{\mathcal{H}}_{0,h}.
    \end{equation*}
\end{lemma}

Let us introduce the following auxiliary space  $\boldsymbol{\mathcal{H}}_h^c$ (see \cite{MR3723734} for instance)
\begin{equation*}
\begin{aligned}
\boldsymbol{\mathcal{H}}_h^c
&:=
\Bigl\{
\btau_h\in \L^2(\Omega)^{d\times d}:\;
\btau_h|_T\in\mathcal{P}_{k+1}(T)^{d\times d}
\quad \forall T\in\mathcal{T}_h,\\
&
\text{the moments of $\btau_h\mathbf{n}$ up to degree $k$ are continuous} \text{ across all interior faces}
\Bigr\}.
\end{aligned}
\end{equation*}

In addition, we define the following space
\begin{equation*}
\boldsymbol{\mathcal{H}}_{0,h}^c:=\Bigl\{ \btau_h\in \boldsymbol{\mathcal{H}}_h^c: \int_\Om \tr(\btau_h)=0\Bigr\}.
\end{equation*}
The following results holds for elements of $\boldsymbol{\mathcal{H}}_{0,h}^c$ (see \cite[Lemma 4.1]{MR3723734}).
\begin{lemma}\label{lema4_1}
   For every  $\btau\in \H^1(\Om)^{d\times d}$, there exists  $\btau_h\in \bH_h^c$ such that  $$\int_F (\btau-\btau_h)\mathbf{n}_F\cdot q=0 \quad \forall q\in \mathcal{P}_k(F)^d \quad\text{and}\quad \|\btau_h\|_{\bH_h}\leq C\|\btau\|_{1,\Om},$$ where  $C>0$ is independent of  $h$.
\end{lemma}


The following result establish an inf-sup  condition for $A_h(\cdot,\cdot)$, which will be useful to prove existence and uniqueness of solution for  \eqref{eq:brinkman_dgf} and consequently, the well definition of operators $\bT_h$ and $\bS_h$.

  \begin{lemma}
  \label{lema:infsupdisc}
 For each  $(\widehat{\bsig}_h,\widehat{\bu}_h)\in \bH_{0,h}\times \mathbf{Q}_h$, there exists $(\btau_h,\bv_h)\in \bH_{0,h}\times \mathbf{Q}_h$ satisfying    $\vertiii{(\btau_h,\bv_h)}\leq C_*\vertiii{(\widehat{\bsig}_h,\widehat{\bu}_h)}$, and a constant $C_A>0$,  depending on  $\texttt{a}_S$, $\nu$, $\mathbf{K}$,  $d$, and  the domain  $\Omega$, such that 
   \begin{equation*}
        \vertiii{(\widehat{\bsig}_h,\widehat{\bu}_h)}^2\leq C_A A_h((\widehat\bsig_h,\widehat\bu_h),(\btau_h,\bv_h)).
    \end{equation*}
\end{lemma}
\begin{proof}
    Let $\delta>0$ and define  $(\btau_h,\bv_h):=(\widehat{\bsig}_h,-\widehat{\bu}_h+\delta\bdiv\widehat{\bsig}_h)$. Replacing this choice in  $A_h(\cdot,\cdot)$, we obtain
    \begin{multline*}
        A_h((\widehat\bsig_h,\widehat\bu_h),(\widehat{\bsig}_h,-\widehat{\bu}_h+\delta\bdiv\widehat{\bsig}_h)) \\
        = a_h(\widehat{\bsig}_h,\widehat{\bsig}_h)+\delta b_h(\widehat{\bsig}_h,\bdiv_h\widehat{\bsig}_h)+c(\widehat{\bu}_h,\widehat{\bu}_h)-\delta c(\widehat{\bu}_h,\bdiv_h\widehat{\bsig}_h)\\
        \geq \frac{1}{\nu}\|\widehat{\bsig}_h^{\texttt{d}}\|^2_{0,\Omega}+\texttt{a}_S\|h_{\mathcal{F}}^{-1/2}\jumpt{\widehat{\bsig}_h}\|^2_{0,\mathcal{F}_h^i}+\delta\|\bdiv_h\widehat{\bsig}_h\|^2_{0,\Om}       \\
         -\delta\|h_{\mathcal{F}}^{-1/2}\jumpt{\widehat{\bsig}_h}\|_{0,\mathcal{F}_h^i}\|h_{\mathcal{F}}^{1/2}\mean{\bdiv_h\widehat{\bsig}_h}\|_{0,\mathcal{F}_h^i}
        +\kappa_*\|\widehat{\bu}_h\|_{0,\Om}-\delta\kappa^*\|\widehat{\bu}_h\|_{0,\Om}\|\bdiv_h\widehat{\bsig}_h\|_{0,\Om}.
    \end{multline*}
    Now, recalling the following well known inequality of real numbers  $ab\leq (2\epsilon_1)^{-1}a^2+2^{-1}\epsilon_1b^2$, with  $a,b\in\R$ and  $\epsilon_1>0$ arbitrary,   and invoking \eqref{eq:dg_trace}, we have
    \begin{multline*}
         -\delta\|h_{\mathcal{F}}^{-1/2}\jumpt{\widehat{\bsig}_h}\|_{0,\mathcal{F}_h^i}\|h_{\mathcal{F}}^{1/2}\mean{\bdiv_h\widehat{\bsig}_h}\|_{0,\mathcal{F}_h^i}\\
          \geq -\frac{\delta \|h_{\mathcal{F}}^{-1/2}\jumpt{\widehat{\bsig}_h}\|_{0,\mathcal{F}_h^i}^2}{2\epsilon_1}-\frac{\delta C_{\texttt{tr}}^2\epsilon_1\|\bdiv_h\widehat{\bsig}_h\|_{0,\Omega}^2}{2},\\
        -\delta\kappa^*\|\widehat{\bu}_h\|_{0,\Om}\|\bdiv_h\widehat{\bsig}_h\|_{0,\Om}
         \geq -\frac{\delta\kappa^*\|\widehat{\bu}_h\|_{0,\Om}^2}{2\epsilon_1}-\frac{\delta\kappa^*\epsilon_1\|\bdiv_h\widehat{\bsig}_h\|^2_{0,\Om}}{2}.
    \end{multline*}
   Gathering all the previous estimates we obtain
    \begin{multline*}
        A_h((\widehat\bsig_h,\widehat\bu_h),(\widehat{\bsig}_h,-\widehat{\bu}_h+\delta\bdiv\widehat{\bsig}_h))\\
         \geq \frac{1}{\nu}\|\widehat{\bsig}_h^{\texttt{d}}\|^2_{0,\Omega}+\delta\left(1-\frac{C^2_{\texttt{tr}}\epsilon_1}{2}-\frac{\kappa^*\epsilon_1}{2}\right)\|\bdiv_h\widehat{\bsig}_h\|^2_{0,\Om}\\
        +\left(\texttt{a}_S-\frac{\delta}{2\epsilon_1}\right)\|h_{\mathcal{F}}^{-1/2}\jumpt{\widehat{\bsig}_h}\|^2_{0,\mathcal{F}_h^i} +\left(\kappa_*-\frac{\delta\kappa^*}{2\epsilon_1}\right)\|\widehat{\bu}_h\|_{0,\Om}^2 \\
         \geq \min\left\{ \frac{1}{\nu },\delta-\frac{\delta C_{\texttt{tr}}^2\epsilon_1}{2}-\frac{\delta\kappa^*\epsilon_1}{2},\texttt{a}_S-\frac{\delta}{2\epsilon_1}\right\}\\
         \times\left(\|\widehat{\bsig}_h^\texttt{d}\|_{0,\Om}^2+\|\bdiv_h\widehat{\bsig}_h\|^2_{0,\Om}+\|h_{\mathcal{F}}^{-1/2}\jumpt{\widehat{\bsig}_h}\|^2_{0,\mathcal{F}_h^i}\right)^{1/2}
         +\left(\kappa_*-\frac{\delta\kappa^*}{2\epsilon_1}\right)\|\widehat{\bu}_h\|_{0,\Om}^2\\
         \geq \frac{1}{2C_1}\min\left\{ \frac{1}{\nu },\delta-\frac{\delta C_{\texttt{tr}}^2\epsilon_1}{2}-\frac{\delta\kappa^*\epsilon_1}{2},\texttt{a}_S-\frac{\delta}{2\epsilon_1}\right\}\|\widehat{\bsig}_h\|_{\bH_h}^2+\left(\kappa_*-\frac{\delta\kappa^*}{2\epsilon_1}\right)\|\widehat{\bu}_h\|_{0,\Om}^2.
    \end{multline*}
    
   Now, invoking Lemma  \ref{lema4_1}, we have that  for any $\widehat{\bu}_h\in \mathbf{Q}_h$, there exists a constant  $C_2>0$ independent of  $h$ and  $\bw_h\in \bH_{h}^c$ such that
    \begin{equation}\label{eq:sobre}
        \bdiv_h{\bw}_h=\widehat{\bu}_h \quad \text{and} \quad  \|\bw_h\|_{\bH_h} \leq  C_2 \|\widehat{\bu}_h\|_{0,\Om}.
    \end{equation}
   Let us define $\tilde{\bsig}_h=\bw_h-\dfrac{\gamma}{d}\mathbf{I}$, where $\gamma=\displaystyle\frac{1}{|\Om|}\int_\Omega \tr(\bw_h)$. The element $\tilde{\bsig}_h$, by construction, is such that  $\tilde{\bsig}_h\in\bH_{0,h}^c$  and since $\bdiv_h\mathbf{I}=0$, we have $\bdiv_h\tilde{\bsig}_h=\widehat{\mathbf{u}}_h$ and  $\|\tilde{\bsig}_h\|_{\bH_h}\leq C_2\|\bu_h\|_{0,\Omega}$.
   
   Next, we set  $(\btau_h,\bv_h)=(\tilde{\bsig}_h,\boldsymbol{0})$. Then, since  $a_h(\cdot,\cdot)$ is bounded,  we obtain
    \begin{multline*}
        A_h((\widehat{\bsig}_h,\widehat{\bu}_h),(\tilde{\bsig}_h,\boldsymbol{0})) = a_h(\widehat{\bsig}_h,\tilde{\bsig}_h)+b_h(\tilde{\bsig}_h,\widehat{\bu}_h)\\
        \geq -C_{\nu,\texttt{a}_S,d}\|\widehat{\bsig}_h\|_{\bH_h}\|\tilde{\bsig}_h\|_{\bH_h}+\int_{\Om}\bdiv_h\tilde{\bsig}_h\cdot\widehat{\bu}_h-\int_{\mathcal{F}_h^i}\jumpt{\tilde{\bsig}_h}\cdot\mean{\widehat{\bu}_h}.
    \end{multline*}
   Invoking  \eqref{eq:sobre}, using the fact that $\int_{\mathcal{F}_h^i}\jumpt{\tilde{\bsig}_h}\cdot\mean{\widehat{\bu}_h}=0$, together with the inequality  $ab\leq (2\epsilon_2)^{-1}a^2+2^{-1}\epsilon_2b^2$ with $\epsilon_2>0$, there holds
    \begin{align*}
        A_h((\widehat{\bsig}_h,\widehat{\bu}_h),(\tilde{\bsig}_h,\boldsymbol{0})) &\geq -C_{\nu,\texttt{a}_S,d}\|\widehat{\bsig}_h\|_{\bH_h}\|\tilde{\bsig}_h\|_{\bH_h}+\|\widehat{\bu}_h\|_{0,\Om}^2\\
        &\geq C_{\nu,\texttt{a}_S,d} \left(-\frac{\|\widehat{\bsig}_h\|^2_{\bH_h}}{2\epsilon_2}-\frac{\epsilon_2\|\tilde{\bsig}_h\|^2_{\bH_h}}{2}\right)+\|\widehat{\bu}_h\|_{0,\Om}^2\\
        &\geq C_{\nu,\texttt{a}_S,d} \left(-\frac{\|\widehat{\bsig}_h\|^2_{\bH_h}}{2\epsilon_2}-\frac{\epsilon_2C_2^2\|\widehat{\bu}_h\|_{0,\Om}^2}{2}\right)+\|\widehat{\bu}_h\|_{0,\Om}^2\\
        &= -\frac{C_{\nu,\texttt{a}_S,d}}{2\epsilon_2}\|\widehat{\bsig}_h\|^2_{\bH_h}+\left(1-\frac{\epsilon_2C_2^2C_{\nu,\texttt{a}_S,d}}{2}\right)\|\widehat{\bu}_h\|_{0,\Om}^2.
    \end{align*}
    
  Now, let  $\beta>0$  and consider $(\btau_h,\bv_h)=(\widehat{\bsig}_h+\beta\tilde{\bsig}_h,-\widehat{\bu}_h+\delta\bdiv_h\widehat{\bsig}_h)$. Then, the following estimate holds     
  \begin{multline*}
        A_h((\widehat{\bsig}_h,\widehat{\bu}_h),(\widehat{\bsig}_h+\beta\tilde{\bsig}_h,-\widehat{\bu}_h+\delta\bdiv_h\widehat{\bsig}_h))\\
         = A_h\bigl((\widehat{\bsig}_h,\widehat{\bu}_h),(\widehat{\bsig}_h,-\widehat{\bu}_h+\delta\bdiv_h\widehat{\bsig}_h)\bigr)+\beta A_h\bigl((\widehat{\bsig}_h,\widehat{\bu}_h),(\tilde{\bsig}_h,\boldsymbol{0})\bigr)\\
        \geq\left(\frac{1}{2C_1}\min\left\{ \frac{1}{\nu },\delta-\frac{\delta C_{\texttt{tr}}^2\epsilon_1}{2}-\frac{\delta\kappa^*\epsilon_1}{2},\texttt{a}_S-\frac{\delta}{2\epsilon_1}\right\}-\frac{\beta C_{\nu,\texttt{a}_S,d}}{2\epsilon_2}\right)\|\widehat{\bsig}_h\|_{\bH_h}^2\\
          +\left[\beta-\frac{\beta\epsilon_2C_2^2C_{\nu,\texttt{a}_S,d}}{2}+\kappa_*-\frac{\delta\kappa^*}{2\epsilon_1}\right]\|\widehat{\bu}_h\|_{0,\Om}^2.
    \end{multline*}
   Then, choosing the parameters $\epsilon_1,\epsilon_2,\beta$  and  $\delta$  as 
    $$\epsilon_1:=\frac{1}{C_{\texttt{tr}}^2+\kappa^*},
    \qquad\epsilon_2:=\frac{1}{C_2^2C_{\nu,\texttt{a}_S,d}},$$
    $$\beta:=\frac{\min\left\{\displaystyle\frac{1}{\nu},\frac{\texttt{a}_S}{2(C_{\texttt{tr}}^2+\kappa^*)},\frac{\texttt{a}_S}{2}\right\}}{2C_1C_2^2C_{\nu,\texttt{a}_S,d}^2},\quad\text{and}\quad
    \delta:=\frac{\texttt{a}_S}{C_{\texttt{tr}}^2+\kappa^*},$$
   there holds
    \begin{equation*}
        A_h((\widehat{\bsig}_h,\widehat{\bu}_h),(\widehat{\bsig}_h+\beta\tilde{\bsig}_h,-\widehat{\bu}_h+\delta\bdiv_h\widehat{\bsig}_h))\geq C_A^{-1}\vertiii{(\widehat{\bsig}_h,\widehat{\bu}_h)}^2,
    \end{equation*}
    where
    \begin{equation*}
\begin{aligned}
C_A^{-1}
:=
\min\Bigg\{&
\frac{
2-\min\left\{
\displaystyle\frac{1}{\nu},
\frac{\texttt{a}_S}{2(C_{\texttt{tr}}^2+\kappa^*)},
\frac{\texttt{a}_S}{2}
\right\}
}{4C_1},
\\[0.3cm]
&
\frac{2\kappa_*-\texttt{a}_S\kappa^*}{2}
+
\frac{
\min\left\{
\dfrac{1}{\nu},
\dfrac{\texttt{a}_S}{2(C_{\texttt{tr}}^2+\kappa^*)},
\dfrac{\texttt{a}_S}{2}
\right\}
}{
4C_1C_2^2C_{\nu,\texttt{a}_S,d}^2
}
\Bigg\}.
\end{aligned}
\end{equation*}

Moreover, it follows that
    \begin{multline*}
        \vertiii{(\widehat{\bsig}_h+\beta\tilde{\bsig}_h,-\widehat{\bu}_h+\delta\bdiv_h\widehat{\bsig}_h)}^2=\|\widehat{\bsig}_h+\beta\tilde{\bsig}_h\|_{\bH_h}^2+\|-\widehat{\bu}_h+\delta\bdiv_h\widehat{\bsig}_h\|_{0,\Om}^2\\
        \leq \max\left\{1+\dfrac{\texttt{a}_S^2}{(C_{\texttt{tr}}^2+\kappa^*)^2},1+\left(\dfrac{\min\left\{\dfrac{1}{\nu},\dfrac{\texttt{a}_S}{2(C_{\texttt{tr}}^2+\kappa^*)},\dfrac{\texttt{a}_S}{2}\right\}}{2C_1C_2C_{\nu,\texttt{a}_S,d}^2}\right)^2\right\}\vertiii{(\widehat{\bsig}_h,\widehat{\bu}_h)}^2.
    \end{multline*}
    Hence, we conclude the proof defining 
    $$C_*:=\sqrt{\max\left\{1+\dfrac{\texttt{a}_S^2}{(C_{\texttt{tr}}^2+\kappa^*)^2},1+\left(\dfrac{\min\left\{\displaystyle\frac{1}{\nu},\dfrac{\texttt{a}_S}{2(C_{\texttt{tr}}^2+\kappa^*)},\dfrac{\texttt{a}_S}{2}\right\}}{2C_1C_2C_{\nu,\texttt{a}_S,d}^2}\right)^2\right\}}.$$
\end{proof}

In order to derive error estimates for the unknowns of the problem, we need to introduce suitable approximations. Since the velocity is approximated with a piecewise polynomial space, as in the finite element method, we can consider the   $L^2$-orthogonal projection as an interpolant operator. On the other hand, for the pseudostress, we employ the techniques used in \cite[Section 4]{MR4108273} to obtain an interpolant that belongs to $\bH_{0,h}$. To this end, we begin considering the classic Scott-Zhang interpolant \cite{MR1011446} and proceed as follows.
Let $\widehat{\bsig} \in \H^s(\Omega)^{d\times d}$ with $\bdiv\widehat{\bsig} \in \H^{1+s}(\Omega)^d$. We denote by $\widehat{\bsig}_I^{SZ} \in \H^1(\Omega)^{d\times d}$ the standard Scott-Zhang interpolation operator, which satisfies the following approximation properties with positive constants $C_{sz1}, C_{sz2} > 0$, independents of $h$
\begin{equation}\label{est:SZ1}
	\|\widehat{\bsig} - \widehat{\bsig}_I^{SZ}\|_{0,\Omega} \le C_{sz1} h^s |\widehat{\bsig}|_{s,\Omega} \quad \text{and} \quad \|\mathbf{div}(\widehat{\bsig} - \widehat{\bsig}_I^{SZ})\|_{0,\Omega} \le C_{sz2} h^{1+s} |\mathbf{div}\,\widehat{\bsig}|_{1+s,\Omega}.
\end{equation}
We proceed by modifying the Scott-Zhang interpolator by
 $$\widehat{\bsig}_I = \widehat{\bsig}_I^{SZ} - \frac{1}{n|\Omega|} \int_{\Omega} \mathrm{tr}(\widehat{\bsig}_I^{SZ})\mathbf{I}  \in \H^1(\Omega)^{d\times d} \cap \bd{\mathcal{H}}_{0,h},$$
 which satisfies the following approximation property
\begin{equation}
	\|\widehat{\bsig} - \widehat{\bsig}_I\|_{\bH(h)} \le C_{app} h^{s} \left( |\widehat{\bsig}|_{s,\Omega} + |\mathbf{div}\,\widehat{\bsig}|_{1+s,\Omega} \right),\quad s>1/2.
\end{equation}


\section{Convergence and error estimates} \label{sec:convegence}
The main aim of this section is to prove convergence of the IPDG discretization. More precisely, we need to establish that the discrete solution operators $\bT_h$ and $\bS_h$ converge to the continuous versions $\bT$ and $\bS$, respectively. Since the numerical scheme is non-conforming,  we proceed as in \cite{MR2220929,MR3962898,MR4077220} where the theory of non-compact operators of  \cite{MR483400,MR483401}  is adapted to this case.

The analysis of convergence, according to \cite{MR483400}, is based in the validation of the following properties:
\begin{itemize}

\item P1.
$\|\boldsymbol{T}-\bT_h\|_{h}
:=
\displaystyle\sup_{\bd{0}\neq\bd{x}_h\in\mathbf{Q}_h}
\frac{
\|(\boldsymbol{T}-\bT_h)\bd{x}_h \|
}{
\|\bd{x}_h\|_{0,\Omega}
}
\to 0$
as  $h\to 0.$

\item P2.
$\forall \bv\in \mathbf{Q}_h$, there holds
\[
\inf_{\bv_h\in \mathbf{Q}_h}
\|\bv-\bv_h\|
\to 0
\quad\text{as } h\to 0.
\]

\end{itemize}
Property P2 follows immediately from the density of continuous piecewise polynomial functions of degree $k$ in $\mathbf{Q}_h$. Hence, the task is to prove P1. Let us begin with the following convergence result for sources in $\mathbf{Q}$.
\begin{theorem}
	\label{th:cea}
	Let  $\bF\in \mathbf{Q}$. If  $\widehat{\bu}:=\bT\bF$ and  $\widehat{\bsig}:=\bS\bF$ are the solutions of \eqref{eq:source_continuous_red} and  $\widehat{\bu}_h:=\bT_h\bF$ and  $\widehat{\bsig}_h:=\bS_h\bF$  are the solutions of  \eqref{eq:brinkman_dgf} then, there exists a constant  $C_\star>0$ independent of $h$ such that \begin{equation*}
		\vertiii{((\bS-\bS_h)\bF,(\bT-\bT_h)\bF)}\leq C_\star \inf_{(\btau_h,\bv_h)\in \bH_{0,h}\times \mathbf{Q}_h} \vertiii{(\widehat{\bsig}-\widehat{\bsig}_h,\widehat{\bu}-\widehat{\bu}_h)}_*.
	\end{equation*}
	Moreover, if  $s\in(0,1]$, there holds 
	\begin{equation*}
		\vertiii{((\bS-\bS_h)\bF,(\bT-\bT_h)\bF)}\leq C_\# h^{s}\left(|\widehat{\bsig}|_{s,\Omega}+|\widehat{\bu}|_{1+s,\Omega}+|\bdiv \widehat{\bsig}|_{{1+s,\Omega}}\right),
	\end{equation*} where  $C_\#>0$ is independent of  $h$.
\end{theorem}
\begin{proof}
	Let us denote by   $\widehat{\bsig}_I$ the modified version  of the Scott-Zhang operator, whereas  $\widehat{\bu}_I$  denotes the  $\L^2$ projection  of  $\widehat{\bu}$.
	
	Now, invoking  Lemma \ref{lema:infsupdisc}, given  $(\widehat{\bsig}_I-\widehat{\bsig}_h,\widehat{\bu}_I-\widehat{\bu}_h)\in\bH_{0,h}\times\mathbf{Q}_h$,  there exists a pair $(\btau_h,\bv_h)\in \bH_{0,h}\times \mathbf{Q}_h$ and a constant $C_*>0$ such that  
	\begin{equation}\label{eq:Lema11}
		\vertiii{(\btau_h,\bv_h)}\leq C_*\vertiii{(\widehat{\bsig}_I-\widehat{\bsig}_h,\widehat{\bu}_I-\widehat{\bu}_h)}.
	\end{equation}	Since the method is consistent, we have the following Galerkin orthogonality property 
	\begin{equation}\label{eq:orto_galerkin} 
	A_h((\widehat{\bsig}-\widehat{\bsig}_h,\widehat{\bu}-\widehat{\bu}_h),(\btau_h,\bv_h)) =0\quad \forall (\btau_h,\bv_h)\in \bH_{0,h}\times \mathbf{Q}_h.
	\end{equation}
	On the other hand, from triangle inequality we have 
	 \begin{multline}\label{eq:desigtrian}
			\vertiii{((\bS-\bS_h)\bF,(\bT-\bT_h)\bF)} = \vertiii{(\widehat{\bsig}-\widehat{\bsig}_h,\widehat{\bu}-\widehat{\bu}_h)}\\
			\leq \vertiii{(\widehat{\bsig}-\widehat{\bsig}_I,\widehat{\bu}-\widehat{\bu}_I)}+\vertiii{(\widehat{\bsig}_I-\widehat{\bsig}_h,\widehat{\bu}_I-\widehat{\bu}_h)}.
	\end{multline}
	
	Let us focus on the second term of the right-hand side of the inequality above.  From the linearity of  $A_h(\cdot,\cdot)$, Lemma \ref{lema:infsupdisc}, the Galerkin orthogonality given in \eqref{eq:orto_galerkin}, and  \eqref{eq:Lema11}, we obtain
	\begin{align*}
		C_A^{-1}\vertiii{(\widehat{\bsig}_I-\widehat{\bsig}_h,\widehat{\bu}_I-\widehat{\bu}_h)}^2&\leq A_h((\widehat{\bsig}_I-\widehat{\bsig}_h,\widehat{\bu}_I-\widehat{\bu}_h),(\btau_h,\bv_h))\\
		&= A_h((\widehat{\bsig}_I-\widehat{\bsig},\widehat{\bu}_I-\widehat{\bu}),(\btau_h,\bv_h))\\
		&\leq 2 C_{\max}\vertiii{(\widehat{\bsig}_I-\widehat{\bsig},\widehat{\bu}_I-\widehat{\bu})}\vertiii{(\btau_h,\bv_h)}\\
		&\leq 2C_*C_{\max}\vertiii{(\widehat{\bsig}_I-\widehat{\bsig},\widehat{\bu}_I-\widehat{\bu})}\vertiii{(\widehat{\bsig}_I-\widehat{\bsig}_h,\widehat{\bu}_I-\widehat{\bu}_h)},
	\end{align*}
	where $C_{\max}:=\max\left\{C_{\nu,\texttt{a}_S,d},(1+C^2_{\texttt{tr}})^{1/2},\kappa^*\right\}$. Therefore, 
	\begin{multline*}
		 \vertiii{(\widehat{\bsig}-\widehat{\bsig}_h,\widehat{\bu}-\widehat{\bu}_h)}\leq(1+ 2C_AC_*C_{\max})\vertiii{(\widehat{\bsig}-\widehat{\bsig}_I,\widehat{\bu}-\widehat{\bu}_I)}\\
		\leq(1+ 2C_AC_*C_{\max})\vertiii{(\widehat{\bsig}-\widehat{\bsig}_I,\widehat{\bu}-\widehat{\bu}_I)}_*.
	\end{multline*}
	Replacing this last estimate  in \eqref{eq:desigtrian}  and using the approximation properties of $\widehat{\bsig}_I$ and $\widehat{\bu}_I$ we obtain  
	\begin{multline*}
	\vertiii{((\bS-\bS_h)\bF,(\bT-\bT_h)\bF)}\\
	\leq(1+2C_AC_*C_{\max})\max\{C_{app},C\}h^{s}(|\widehat{\bd\sigma}|_{s,\Omega}+|\widehat{\bd{u}}|_{1+s,\Omega}+|\bdiv\widehat{\bd\sigma}|_{1+s,\Omega}).
	\end{multline*}
which in conjunction with Theorem \ref{th:reg_velocity} yields to
$$\vertiii{((\bS-\bS_h)\bF,(\bT-\bT_h)\bF)}\leq Ch^s\|\bF\|_{0,\Omega},$$
concluding the proof.
\end{proof}

It is possible to improve the order of convergence for the velocity error with the IPDG method in $\L^2(\O)^d$ norm. This is established in the following result.
\begin{theorem}
	There exists  $C>0$ independent of $h$, such that  
	$$\|\widehat{\bu}-\widehat{\bu}_h\|_{0,\Omega}\leq C^*h^{2s}\|\bF\|_{0,\Omega},$$ 
	where  $s\in (0,1]$ is as in Theorem \ref{th:reg_velocity}.
\end{theorem}
\begin{proof}
	The proof is based in a classic duality argument. First, we consider the following auxiliary source  problem: find $(\bd{\varphi},\bd{\psi})$ such that the following system holds
\begin{equation}
\label{eq:dual_problem}
\left\{
\begin{aligned}
\bdiv\bd{\varphi}-\mathbf{K}^{-1}\bd{\psi}
&= \widehat{\bu}-\widehat{\bu}_h
&& \text{in } \Omega,\\
\bd{\varphi}^{\texttt{d}}-\nu\nabla\bd{\psi}
&= \boldsymbol{0}
&& \text{in } \Omega,\\
\bd{\psi}
&= \bd{0}
&& \text{on } \partial\Omega.
\end{aligned}
\right.
\end{equation}
for which there exists   $s\in(0,1]$ such that the solution satisfies
	\begin{equation}
	\label{eq:dep_dual}
	 \|\bd{\varphi}\|_{s,\Omega}+\|\bdiv\bd{\varphi}\|_{1+s,\Omega}+\|\bd{\psi}\|_{1+s,\Omega}\leq C\|\widehat{\bu} -\widehat{\bu}_h\|_{0,\Omega},
	 \end{equation}
	and if $\O$ is convex, then $s=1$. From the first equation of \eqref{eq:dual_problem} we have 
	\begin{equation} \label{eq:sub}
		\|\widehat{\bu} -\widehat{\bu}_h\|_{0,\Omega}^2=\sum_{T \in \mathcal{T}_h}
		\int_T \bdiv\bd\varphi\cdot(\widehat{\bd{u}}-\widehat{\bd{u}}_h)
		-\sum_{T\in\mathcal{T}_h}
		\int_T \mathbf{K}^{-1}\bd{\psi}\cdot
		(\widehat{\bu} - \widehat{\bu}_h),
	\end{equation} 
	which, given the regularity assumed on $\bd\varphi$ that implies $\jumpt{\bd\varphi}=0$, \eqref{eq:sub}  can be rewritten as follows
	\begin{equation*}
		\|\widehat{\bu} -\widehat{\bu}_h\|_{0,\Omega}^2=\sum_{T \in \mathcal{T}_h}
		\int_T \bdiv\bd\varphi\cdot(\widehat{\bd{u}}-\widehat{\bd{u}}_h)
		-\int_{\mathcal{F}_h^i}\jumpt{\bd\varphi}\cdot\mean{\widehat{\bu} -\widehat{\bu}_h}-\sum_{T\in\mathcal{T}_h}
		\int_T \mathbf{K}^{-1}\bd{\psi}\cdot
		(\widehat{\bu} - \widehat{\bu}_h).
	\end{equation*}
	
	Now, testing $\btau_h=\bd\varphi_I$ and 
	$\bv_h=\bd\psi_I$ in \eqref{eq:orto_galerkin} we obtain
	\begin{equation}
	\label{eq:orto_dual}
	 A_h((\widehat{\bsig}-\widehat{\bsig}_h,\widehat{\bu}-\widehat{\bu}_h),(\bd\varphi_I,\bd\psi_I))=0,
	 \end{equation}
	and then, subtracting \eqref{eq:orto_dual} to  \eqref{eq:sub} yields to
	\begin{multline}
	\label{eq:duality_divergencia1}
		\|\widehat{\bu} - \widehat{\bu}_h\|_{0,\Omega}^2
		=
		\sum_{T \in \mathcal{T}_h}
		\int_T \bdiv(\bd\varphi-\bd\varphi_I)\cdot(\widehat{\bd{u}}-\widehat{\bd{u}}_h)
		 - \int_{\mathcal{F}_h^i}\jumpt{\bd\varphi-\bd\varphi_I}\cdot\mean{\widehat{\bu} -\widehat{\bu}_h}\\
		 -\sum_{T\in\mathcal{T}_h}
		\int_T \mathbf{K}^{-1}(\bd{\psi}-\bd{\psi}_I)\cdot
		(\widehat{\bu} -\widehat{\bu}_h)- \frac{1}{\nu}\int_\Omega (\widehat{\bsig}-\widehat{\bsig}_h)^\texttt{d}:\bd\varphi_I^\texttt{d}\\
		-\int_{\mathcal{F}_h^i}\frac{{\texttt{a}_S}}{h_\mathcal{F}}\jumpt{\widehat{\bsig}-\widehat{\bsig}_h}\cdot\jumpt{\bd\varphi_I}
		-\sum_{T \in \mathcal{T}_h}
		\int_T \bdiv(\widehat{\bsig}-\widehat{\bsig}_h)\cdot\bd\psi_I + \int_{\mathcal{F}_h^i} \jumpt{\widehat{\bsig}-\widehat{\bsig}_h} \cdot\mean{\bd\psi_I}.
	\end{multline}
	An application of  integration by parts shows the following identity
	\begin{multline}
	\label{eq:duality_divergencia2}
		-\sum_{T \in \mathcal{T}_h}
		\int_T \bdiv(\widehat{\bsig}-\widehat{\bsig}_h)\cdot\bd\psi_I
		= \sum_{T \in \mathcal{T}_h}
		\int_T (\widehat{\bsig}-\widehat{\bsig}_h) \cdot\nabla \bd\psi_I-\sum_{T \in \mathcal{T}_h}
		\int_{\partial T} ( (\widehat{\bsig}-\widehat{\bsig}_h)\mathbf{ n}_T)\cdot\bd\psi_I\\
		=\sum_{T \in \mathcal{T}_h}
		\int_T  (\widehat{\bsig}-\widehat{\bsig}_h) \cdot\nabla \bd\psi_I - \int_{\mathcal{F}_h^i}\jumpt{\widehat{\bsig}-\widehat{\bsig}_h}\cdot\mean{\bd\psi_I}-\int_{\mathcal{F}_h}\mean{\widehat{\bsig}-\widehat{\bsig}_h}:\jump{\bd\psi_I}.
	\end{multline}
	Hence, replacing \eqref{eq:duality_divergencia2} in  \eqref{eq:duality_divergencia1}  and using once again that  $\jumpt{\bd\varphi}=0$  we obtain
	\begin{multline*}
		\|\widehat{\bu} - \widehat{\bu}_h\|_{0,\Omega}^2
		=
		\sum_{T \in \mathcal{T}_h}
		\int_T \bdiv(\bd\varphi-\bd\varphi_I)\cdot(\widehat{\bd{u}}-\widehat{\bd{u}}_h)
		- \int_{\mathcal{F}_h^i}\jumpt{\bd\varphi-\bd\varphi_I}\cdot\mean{\widehat{\bu} -\widehat{\bu}_h}\\
		-\sum_{T\in\mathcal{T}_h}
		\int_T \mathbf{K}^{-1}(\bd{\psi}-\bd{\psi}_I)\cdot
		(\widehat{\bu} -\widehat{\bu}_h)- \frac{1}{\nu}\int_\Omega (\widehat{\bsig}-\widehat{\bsig}_h)^\texttt{d}:\bd\varphi_I^\texttt{d}
		\\
		+\int_{\mathcal{F}_h^i}\frac{{\texttt{a}_S}}{h_\mathcal{F}}\jumpt{\widehat{\bsig}-\widehat{\bsig}_h}\cdot\jumpt{\bd\varphi-\bd\varphi_I}
		+\sum_{T \in \mathcal{T}_h}
		\int_T  (\widehat{\bsig}-\widehat{\bsig}_h) \cdot\nabla \bd\psi_I -\int_{\mathcal{F}_h}\mean{\widehat{\bsig}-\widehat{\bsig}_h}:\jump{\bd\psi_I}.
	\end{multline*}
On the other hand, from the second equation of \eqref{eq:dual_problem} we have 
 $$\displaystyle\int_\Omega (\widehat{\bsig}-\widehat{\bsig}_h):\left(\frac{1}{\nu}\bd\varphi^\texttt{d}-\nabla\bd\psi\right)=0,$$ 
implying that
$$- \frac{1}{\nu}\int_\Omega (\widehat{\bsig}-\widehat{\bsig}_h)^\texttt{d}:\bd\varphi_I^\texttt{d}= \frac{1}{\nu}\int_\Omega (\widehat{\bsig}-\widehat{\bsig}_h)^\texttt{d}:(\bd\varphi-\bd\varphi_I)^\texttt{d}-\int_{\Omega}(\widehat{\bsig}-\widehat{\bsig}_h):\nabla \bd\psi.$$
Moreover, integration by parts, together with the fact that $\jump{\bd\psi}=0$, yields the following identity
	\begin{multline*}
		\sum_{T \in \mathcal{T}_h}
		\int_T  (\widehat{\bsig}-\widehat{\bsig}_h) \cdot\nabla (\bd\psi_I -\bd\psi)\\
		 = -\sum_{T \in \mathcal{T}_h}
		\int_T \bdiv  (\widehat{\bsig}-\widehat{\bsig}_h)\cdot (\bd\psi_I -\bd\psi) + \sum_{T \in \mathcal{T}_h}
		\int_{\partial T} ((\widehat{\bsig}-\widehat{\bsig}_h)\mathbf{ n}_T)\cdot (\widehat{\bsig}-\widehat{\bsig}_h)\\
		=\sum_{T \in \mathcal{T}_h}
		\int_T \bdiv  (\widehat{\bsig}-\widehat{\bsig}_h)\cdot (\bd\psi -\bd\psi_I) + \int_{\mathcal{F}_h^i}\jumpt{\widehat{\bsig}-\widehat{\bsig}_h}\cdot\mean{\bd\psi_I -\bd\psi}+\int_{\mathcal{F}_h}\mean{\widehat{\bsig}-\widehat{\bsig}_h}:\jump{\bd\psi_I -\bd\psi}\\
		=\sum_{T \in \mathcal{T}_h}
		\int_T \bdiv  (\widehat{\bsig}-\widehat{\bsig}_h)\cdot (\bd\psi -\bd\psi_I)  - \int_{\mathcal{F}_h^i}\jumpt{\widehat{\bsig}-\widehat{\bsig}_h}\cdot\mean{\bd\psi -\bd\psi_I}+\int_{\mathcal{F}_h}\mean{\widehat{\bsig}-\widehat{\bsig}_h}:\jump{\bd\psi_I}.
	\end{multline*}
Hence, \eqref{eq:sub} is stated now as follows
	\begin{multline*}
		\|\widehat{\bu} - \widehat{\bu}_h\|_{0,\Omega}^2
		=
		\sum_{T \in \mathcal{T}_h}
		\int_T \bdiv(\bd\varphi-\bd\varphi_I)\cdot(\widehat{\bd{u}}-\widehat{\bd{u}}_h)
		- \int_{\mathcal{F}_h^i}\jumpt{\bd\varphi-\bd\varphi_I}\cdot\mean{\widehat{\bu} -\widehat{\bu}_h}\\
		-\sum_{T\in\mathcal{T}_h}
		\int_T \mathbf{K}^{-1}(\bd{\psi}-\bd{\psi}_I)\cdot
		(\widehat{\bu} -\widehat{\bu}_h)
		+\frac{1}{\nu}\int_\Omega (\widehat{\bsig}-\widehat{\bsig}_h)^\texttt{d}:(\bd\varphi-\bd\varphi_I)^\texttt{d}	+\int_{\mathcal{F}_h^i}\frac{{\texttt{a}_S}}{h_\mathcal{F}}\jumpt{\widehat{\bsig}-\widehat{\bsig}_h}\cdot\jumpt{\bd\varphi-\bd\varphi_I}\\
		+\sum_{T \in \mathcal{T}_h}
		\int_T \bdiv  (\widehat{\bsig}-\widehat{\bsig}_h)\cdot (\bd\psi -\bd\psi_I)  - \int_{\mathcal{F}_h^i}\jumpt{\widehat{\bsig}-\widehat{\bsig}_h}\cdot\mean{\bd\psi -\bd\psi_I}\\
		=A_h((\bd\varphi-\bd\varphi_I,\bd\psi -\bd\psi_I),(\widehat{\bsig}-\widehat{\bsig}_h,\widehat{\bd{u}}-\widehat{\bd{u}}_h)).
	\end{multline*}
	Now, from the boundedness of  $A_h(\cdot,\cdot)$, the approximation properties of $\bd\varphi_I$ and $\bd\psi_I$, and \eqref{eq:dep_dual}, we obtain 
	\begin{equation*}
\begin{aligned}
\|\widehat{\bu}-&\widehat{\bu}_h\|_{0,\Omega}^2\\
&= C_{\nu,\texttt{a}_S,d}
\|\widehat{\bsig}-\widehat{\bsig}_h\|_{\bH_h}
\|\bd\varphi-\bd\varphi_I\|_{\bH_h}
+(1+C_{\texttt{tr}}^2)
\|\widehat{\bsig}-\widehat{\bsig}_h\|_{\bH_h}
\|\bd\psi-\bd\psi_I\|_{0,\Omega}
\\
&+(1+C_{\texttt{tr}}^2)
\|\bd\varphi-\bd\varphi_I\|_{\bH_h}
\|\widehat{\bu}-\widehat{\bu}_h\|_{0,\Omega}
+\kappa^*
\|\bd\psi-\bd\psi_I\|_{0,\Omega}
\|\widehat{\bu}-\widehat{\bu}_h\|_{0,\Omega}
\\
&\leq C_{\nu,\texttt{a}_S,d}Ch^s
\|\bd\varphi\|_{s,\Omega}
\|\widehat{\bsig}-\widehat{\bsig}_h\|_{\bH_h}
+(1+C_{\texttt{tr}}^2)Ch^s
\|\bd\psi\|_{1+s,\Omega}
\|\widehat{\bsig}-\widehat{\bsig}_h\|_{\bH_h}
\\
&+(1+C_{\texttt{tr}}^2)Ch^s
\|\bd\varphi\|_{s,\Omega}
\|\widehat{\bu}-\widehat{\bu}_h\|_{0,\Omega}
+\kappa^*Ch^s
\|\bd\psi\|_{1+s,\Omega}
\|\widehat{\bu}-\widehat{\bu}_h\|_{0,\Omega}
\\
&\leq C_{\nu,\texttt{a}_S,d}Ch^s
\|\widehat{\bu}-\widehat{\bu}_h\|_{0,\Omega}
\|\widehat{\bsig}-\widehat{\bsig}_h\|_{\bH_h}
+(1+C_{\texttt{tr}}^2)Ch^s
\|\widehat{\bu}-\widehat{\bu}_h\|_{0,\Omega}
\|\widehat{\bsig}-\widehat{\bsig}_h\|_{\bH_h}
\\
&+(1+C_{\texttt{tr}}^2)Ch^s
\|\widehat{\bu}-\widehat{\bu}_h\|_{0,\Omega}^2
+\kappa^*Ch^s
\|\widehat{\bu}-\widehat{\bu}_h\|_{0,\Omega}^2.
\end{aligned}
\end{equation*}

	Finally, invoking \eqref{eq:reg_u_if}, we have 
	\begin{equation*}
		\|\widehat{\bu} - \widehat{\bu}_h\|_{0,\Omega} 
		\leq
		C^* h^{2s}\left(|\widehat{\bsig}|_{s,\Omega}+|\widehat{\bu}|_{1+s,\Omega}+|\bdiv \widehat{\bsig}|_{{1+s,\Omega}}\right)
		\leq
		C^* h^{2s}\|\bF\|_{0,\Omega},
	\end{equation*}
	where the constant $C^*$ is independent of $h$ and is defined by  
	$$C^*:=\max\left\{C_{\nu,\texttt{a}_S,d}+1+C_{\texttt{tr}}^2,\kappa^*+1+C_{\texttt{tr}}^2\right\}C C_\#.$$ 
	This completes the proof.
\end{proof}

As a direct consequence of  Theorem \ref{th:cea}, we have the following result by assuming a discrete source. Since the proof is the same,  we do not incorporate further details.
\begin{corollary}\label{cor:norm}
	Given  $\bF_h\in \mathbf{Q}_h$, the following estimate holds
	$$\vertiii{((\bS-\bS_h)\bF_h,(\bT-\bT_h)\bF_h)}\leq C_{\#}h^s\|\bF_h\|_{0,\Omega},$$ where $s\in(0,1]$ as in Theorem \ref{th:reg_velocity}  and  $C_{\#}>0$ is independent of  $h$ as in Theorem \ref{th:cea}.
\end{corollary}
Now we are in position of establish property P1.
\begin{lemma}
	There exists  $C>0$, independent of $h$, such that  
	$$\|\bT-\bT_h\|_{h}\leq C_{\#} h^s,$$
	with $s\in(0,1]$ as in Theorem \ref{th:reg_velocity}.
\end{lemma}
\begin{proof}
	Given $\bF_h\in \mathbf{Q}_h$, we have 
	\begin{align*}
	\|\bT-\bT_h\|_{h}&:=\sup_{\bd{0}\neq {\bF_h\in\mathbf{Q}_h}}\frac{\|(\bT-\bT_h)\bF_h\|_{0,\Omega}}{{\|\bF_h\|_{0,\Omega}}}\\
	&\leq\sup_{\bd{0}\neq {\bF_h\in\mathbf{Q}_h}}\frac{\vertiii{((\bS-\bS_h)\bF_h,(\bT-\bT_h)\bF_h)}}{{\|\bF_h\|_{0,\Omega}}}\leq C_{\#} h^s,
	\end{align*}
	where we have used Corollary \ref{cor:norm}.  This concludes the proof.
\end{proof}

 Now that we have the convergence of  $\bT_h$  to  $\bT$ as $h\rightarrow 0$  according to property P1,  we   invoke  the results of  \cite[Chapter IV]{MR0203473},  \cite[Theorem 9.1]{MR2652780} and  \cite[Theorem 1]{MR483400}, to conclude that   the IPDG method do no introduce spurious eigenvalues.
\begin{theorem}\label{free}
Let $E\subset\mathbb{C}$ be a compact subset not intersecting $\sp(\bT)$.
Then, there exists $h_0>0$ such that, if $h\leq h_0$, then $E\cap\sp(\bT_h)=\emptyset.$ 
\end{theorem}

\subsection{Error estimates}
Now the aim is to derive a priori error estimates for the eigenvalues and eigenfunctions. We begin by introducing some definitions.   First, we define the  \emph{gap} $\hat{\delta}$ between two closed subspaces  $\bd{\mathcal{M}}$ and $\bd{\mathcal{N}}$ of  $\mathbf{Q}$ by  $$\hat{\delta}(\bd{\mathcal{M}},\bd{\mathcal{N}}):=\max\left\{\delta(\bd{\mathcal{M}},\bd{\mathcal{N}}),\delta(\bd{\mathcal{N}},\bd{\mathcal{M}})\right\},$$ where $$\delta(\bd{\mathcal{M}},\bd{\mathcal{N}}):=\sup_{\substack{\bd{x}\in\bd{\mathcal{M}}\\ \|\bd{x}\|_{0,\Omega}=1}}\left(\inf_{\bd{y}\in\bd{\mathcal{N}}}\|\bd{x}-\bd{y}\|_{0,\Omega}\right).$$

Let  $\mu$ be an isolated eigenvalue of  $\bT$ and let $D$ be an open disk of the complex plane with boundary  $\Gamma$ such that $\mu$ is the only eigenvalue of  $\bT$  contained by  $D$, satisfying  $\Gamma\cap\text{sp}(\bT)=\emptyset$. Let us introduce the spectral projection  of  $\bT$ and  $\bT_h$  which are defined, respectively, by
$$\bd{\mathcal{E}}:=\frac{1}{2\pi i}\int_\Gamma (z\mathbf{I}-\bT)^{-1}\,dz:\mathbf{Q}\to\mathbf{Q},$$
$$\bd{\mathcal{E}}_h:=\frac{1}{2\pi i}\int_\Gamma (z\mathbf{I}-\bT_h)^{-1}\,dz:\mathbf{Q}_h\to\mathbf{Q}_h.$$


Next, we need to consider that the invariant space  $\bd{\mathcal{E}}(\mathbf{Q})$ is sufficiently regular, in the sense that $\bu\in\bd{\mathcal{E}}(\mathbf{Q})\subset\mathbf{H}^{1+r}(\O)$ with  $r>0$ and also, $\bsig\in\boldsymbol{\mathcal{H}}^{r}(\O)$. Now, following \cite[Theorem 5.1]{MR2220929}, we have the following result.
\begin{lemma}
	There holds
	$$\lim_{h\to0}\|\bd{\mathcal{E}}-\bd{\mathcal{E}}_h\|_{\mathcal{L}(\mathbf{Q},\mathbf{Q}_h)}=0.$$
\end{lemma}

On the other hand, by adapting  \cite[Lemma 7]{MR4623018} to our symmetric IPDG method, we obtain the following result.
\begin{lemma}
\label{lmm:gap}
	There exists a constant $h_0$  such that for  $h<h_0$, there holds  $$\hat{\delta}(\bd{\mathcal{E}}(\mathbf{Q}),\bd{\mathcal{E}}_h(\mathbf{Q}_h))\leq Ch^{\min\{r,k\}},$$ where $r>0$ is as in \eqref{eq:reg_u_if2} and  $C>0$ is independent of  $h$.
\end{lemma}

An estimate for the error $\|\bsig-\bsig_h\|_{\boldsymbol{\mathcal{H}}_h}$ can be obtained in the same way as in Theorem \ref{th:cea}, but considering the regularity of $\bsig$ instead of $\widehat{\bsig}$. More precisely, if   $
(\bsig,\bu)\in \H^r(\Omega)^{d\times d}\times \H^{1+r}(\Omega)^d$ with $r>0$, together with the following estimate \begin{equation}\label{eq:cota_norm}
	\vertiii{(\bsig-\bsig_h,\bu-\bu_h)}\leq C h^{\min\{r,k\}} \left(\|\bsig\|_{r,\Omega} +\|\bu\|_{1+r,\Omega}\right),
\end{equation}
which shows that the gap in Lemma \ref{lmm:gap} is controlled by the velocity and pseudostress errors of the eigenvalue problem.

Additionally, the following estimate for the pressure can be established, as stated in the next corollary.
\begin{corollary}
	There exists a constant $h_0$ such that $h<h_0$ and the following estimate holds
	$$\|p - p_h\|_{0,\Omega}\leq C_{\texttt{pr}}h^{\min\{r,k\}},$$ where $r>0$ as in \eqref{eq:reg_u_if2} and $C_{\texttt{pr}}>0$ independent of $h$.
\end{corollary}
\begin{proof}
	Let us assume $\|\bu\|_{0,\Omega}=1$. Since $p=-\displaystyle\frac{1}{d}\tr(\bsig)$ and $p_h=-\displaystyle\frac{1}{d}\tr(\bsig_h)$, we have 
	\begin{multline*}
		\|p - p_h\|_{0,\Omega}\\
		 = \frac{1}{d}\|\tr(\bsig-\bsig_h)\|_{0,\Omega}
		\leq \frac{\sqrt{d}}{d}\|\bsig-\bsig_h\|_{\bH_h}
		\leq \frac{\sqrt{d}}{d}Ch^{\min\{r,k\}}(\|\bsig\|_{r,\Omega}+\|\bu\|_{1+r,\Omega})\\
		 \leq \frac{\sqrt{d}}{d}CC_{\lambda}h^{\min\{r,k\}}\|\bu\|_{0,\Omega}
		\leq \frac{\sqrt{d}}{d}CC_{\lambda}h^{\min\{r,k\}},
	\end{multline*}
	where we have used \eqref{eq:cota_norm}, \eqref{eq:reg_u_if2} and the fact that $\|\bu\|_{0,\O}=1$. Defining $C_{\texttt{pr}}:= \frac{\sqrt{d}}{d}CC_{\lambda}$, we conclude the proof.
\end{proof}

Now, it is  possible to improve the convergence order of the eigenvalues. This is established in the following result. 
\begin{theorem}
	There exists a constant $h_0$  such that  $h<h_0$ and the following estimate holds 
	\begin{equation}
		\label{orddob}
		|\lambda-\lambda_h|\leq \widehat{C}h^{2\min\{r,k\}},
	\end{equation} 
	where  $r>0$  and  $\widehat{C}>0$  depends on the physical parameters but not on  $h$.
\end{theorem}
\begin{proof}
Let $(\lambda,(\bsig,\bu))$  and $(\lambda_h,(\bsig_h,\bu_h))$  be the solutions of the respective continuous and discrete eigenvalue problems. Let us assume that $\|\bu_h\|_{0,\Omega}=\|\bu\|_{0,\Omega}=1$. Consider the following well known identity
	  \begin{equation*}
		A_h((\bsig-\bsig_h,\bu-\bu_h),(\bsig-\bsig_h,\bu-\bu_h))+\lambda (\bu-\bu_h,\bu-\bu_h)_{0,\Omega}=(\lambda-\lambda_h)(\bu_h,\bu_h)_{0,\Omega}.
	\end{equation*}
	Since the eigenfunctions are sufficiently regular, triangle inequality reveals immediately 
	\begin{align*}
		|\lambda-\lambda_h|&\leq |	A_h((\bsig-\bsig_h,\bu-\bu_h),(\bsig-\bsig_h,\bu-\bu_h))|+|\lambda|| (\bu-\bu_h,\bu-\bu_h)_{0,\Omega}|\\
		&\leq M^*\vertiii{(\bsig-\bsig_h,\bu-\bu_h)}_*^2+|\lambda|\|\bu-\bu_h\|_{0,\Omega}^2\\
		&\leq (M^*+|\lambda|)\vertiii{(\bsig-\bsig_h,\bu-\bu_h)}_*^2.
	\end{align*}
	where $M^*>0$ is the constant involved in \eqref{eq:formA_M_star} which depends on the dimension, the viscosity, the permeability, and the stabilization parameter $\texttt{a}_S$.
	From the definition of $\vertiii{\cdot}_*$  we have
	\begin{equation}
	\label{eq:padra0}
	\vertiii{(\bsig-\bsig_h,\bu-\bu_h)}_*^2:= \vertiii{(\bsig-\bsig_h,\bu-\bu_h)}^2+\|h_{\mathcal{F}}^{1/2}\mean{\bdiv_h(\bsig-\bsig_h)}\|_{0,\mathcal{F}_h^i}^2.
	\end{equation}
	
	Now the task is to estimate each contribution on the right hand side of \eqref{eq:padra0}. For the first, we invoke \eqref{eq:cota_norm} obtaining
	\begin{equation}
	\label{eq:padra1}
		\vertiii{(\bsig-\bsig_h,\bu-\bu_h)}\leq C h^{\min\{r,k\}} \left(\|\bsig\|_{r,\Omega} +\|\bu\|_{1+r,\Omega}\right) 
		\leq CC_{\lambda} h^{\min\{r,k\}}\|\bu\|_{0,\Omega}.
	\end{equation}
	
	For the second terms, if $\bsig_I$ is the best approximation of $\bsig$, we have
	\begin{equation*}
		\|h_{\mathcal{F}}^{1/2}\mean{\bdiv_h(\bsig-\bsig_h)}\|_{0,\mathcal{F}_h^i}^2 \leq \|h_{\mathcal{F}}^{1/2}\mean{\bdiv_h(\bsig-\bsig_I)}\|_{0,\mathcal{F}_h^i}^2+\|h_{\mathcal{F}}^{1/2}\mean{\bdiv_h(\bsig_I-\bsig_h)}\|_{0,\mathcal{F}_h^i}^2,
	\end{equation*}
	and then, using the approximation properties  $\bsig_I$ and the regularity of $\bsig$ given in   \eqref{eq:reg_u_if2} we obtain
	\begin{align*}
		\|h_{\mathcal{F}}^{1/2}\mean{\bdiv_h(\bsig-\bsig_I)}\|_{0,\mathcal{F}_h^i} &\leq \sum_{F\in\mathcal{F}_h^i} h_F\|\bdiv(\bsig-\bsig_I)\|_{0,F}\\
		&\leq \sum_{T \in \mathcal{T}_h} C h_T^{\min\{1+r,k\}}|\bdiv\bsig|_{1+r,T}\\
		&\leq C h^{2\min\{1+r,k\}}\|\bdiv\bsig\|_{1+r,\Omega}\\
		&\leq CC_{\lambda} h^{\min\{1+r,k\}}\|\bu\|_{0,\Omega}.
	\end{align*}
	For the remaining term, invoking trace inequality, the additional regularity and \eqref{eq:cota_norm} yields to 
	\begin{align*}
		\|h_{\mathcal{F}}^{1/2}\mean{\bdiv_h(\bsig_I-\bsig_h)}\|_{0,\mathcal{F}_h^i}&\leq C_{\texttt{tr}} \|\bdiv_h(\bsig_I-\bsig_h)\|_{0,\Omega}\\
		&\leq  C_{\texttt{tr}} \left(\|\bdiv_h(\bsig_I-\bsig)\|_{0,\Omega} + \|\bdiv_h(\bsig-\bsig_h)\|_{0,\Omega}\right)\\
		&\leq  C_{\texttt{tr}} \left(\|\bdiv_h(\bsig_I-\bsig)\|_{0,\Omega} + \vertiii{(\bsig-\bsig_h,\bu-\bu_h)}\right)\\
		&\leq  C_{\texttt{tr}} Ch^{\min\{r,k\}}\left(\|\bdiv\bsig\|_{1+r,\Omega} + \|\bsig\|_{r,\Omega}+\|\bu\|_{1+r,\Omega}\right)\\
		&\leq  C_{\texttt{tr}} C_{\lambda}Ch^{\min\{r,k\}}\|\bu\|_{0,\Omega}.
	\end{align*}
	Therefore
	\begin{equation}
	\label{eq:padra2}
		\|h_{\mathcal{F}}^{1/2}\mean{\bdiv_h(\bsig-\bsig_h)}\|_{0,\mathcal{F}_h^i}^2 \leq C_{\lambda}^2(C+C_{\texttt{tr}}^2)h^{2\min\{r,k\}}\|\bu\|_{0,\Omega}^2.
	\end{equation}
	Replacing \eqref{eq:padra1} and \eqref{eq:padra2} in   \eqref{eq:padra0} yields to
	\begin{align*}
		\vertiii{(\bsig-\bsig_h,\bu-\bu_h)}_*^2 &\leq CC_{\lambda}^2(2C+C_{\texttt{tr}}^2)h^{2\min\{r,k\}}\|\bu\|_{0,\Omega}^2.
	\end{align*}
	Finally, since  $\|\bu\|_{0,\Omega}=1$, we obtain
	\begin{align*}
		|\lambda-\lambda_h| &\leq (M^*+|\lambda|)CC_{\lambda}^2(2C+C_{\texttt{tr}}^2)h^{2\min\{r,k\}}\|\bu\|_{0,\Omega}^2\\
		& \leq (M^*+|\lambda|)CC_{\lambda}^2(2C+C_{\texttt{tr}}^2)h^{2\min\{r,k\}},
	\end{align*}
	and defining  $\widehat{C}:=(M^*+|\lambda|)CC_{\lambda}^2(2C+C_{\texttt{tr}}^2)$, we conclude the proof.
\end{proof}

\section{Numerical experiments}
\label{sec:numerics}
In this section we present a series of numerical experiments designed to assess the performance of the proposed finite element approximation. All computations were carried out within the \texttt{FEniCSx} framework~\cite{barrata2023dolfinx,scroggs2022basix}, employing the standard eigensolvers provided by \texttt{SLEPc}~\cite{hernandez2005slepc} to compute the discrete spectrum.


\subsection{Spurious analysis}
We begin our numerical study with the influence of the stabilization parameter on the emerging of pollution on the computed spectrum. To do this, the strategy is the following: we fix a mesh refinement level and a particular choice of $\textbf{K}^{-1}$ for which we compute some components of the spectrum. With this fixed mesh, we select different ways to define $\texttt{a}_S$ and compute the eigenvalues for different polynomial degrees. To determine whether a computed eigenvalue is a spurious eigenvalue or a correct approximation of an exact one, it is essential to examine the plot of the associated eigenfunction.

Let us define the inverse of the permeability tensor  for our tests
\begin{equation}
\label{eq:permea_tests}
\textbf{K}^{-1}:=
\begin{cases}
	\kappa\,\textbf{I}, & \text{in } \Omega_D,\\[4pt]
	\boldsymbol{0}, & \text{in } \Omega_S,
\end{cases}
\end{equation}
where $\kappa>0$ is a parameter that may change in such a way that $\textbf{K}^{-1}$ can produce an obstacle on the domain or not (where this last case is a Stokes type of eigenvalue problem).

\begin{figure}[!htbp]
	\centering
	\includegraphics[width=0.4\linewidth]{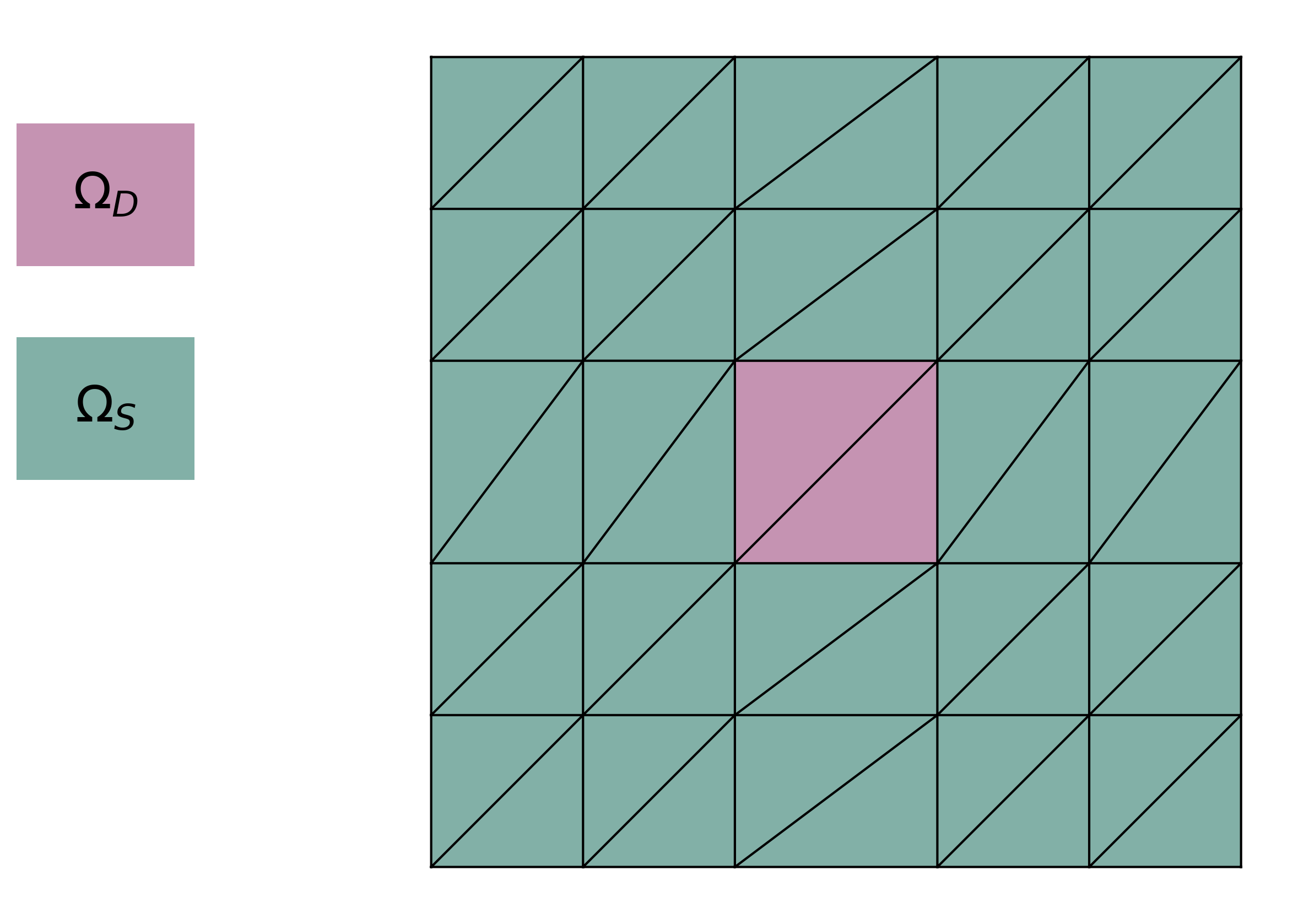}
	\includegraphics[width=0.4\linewidth]{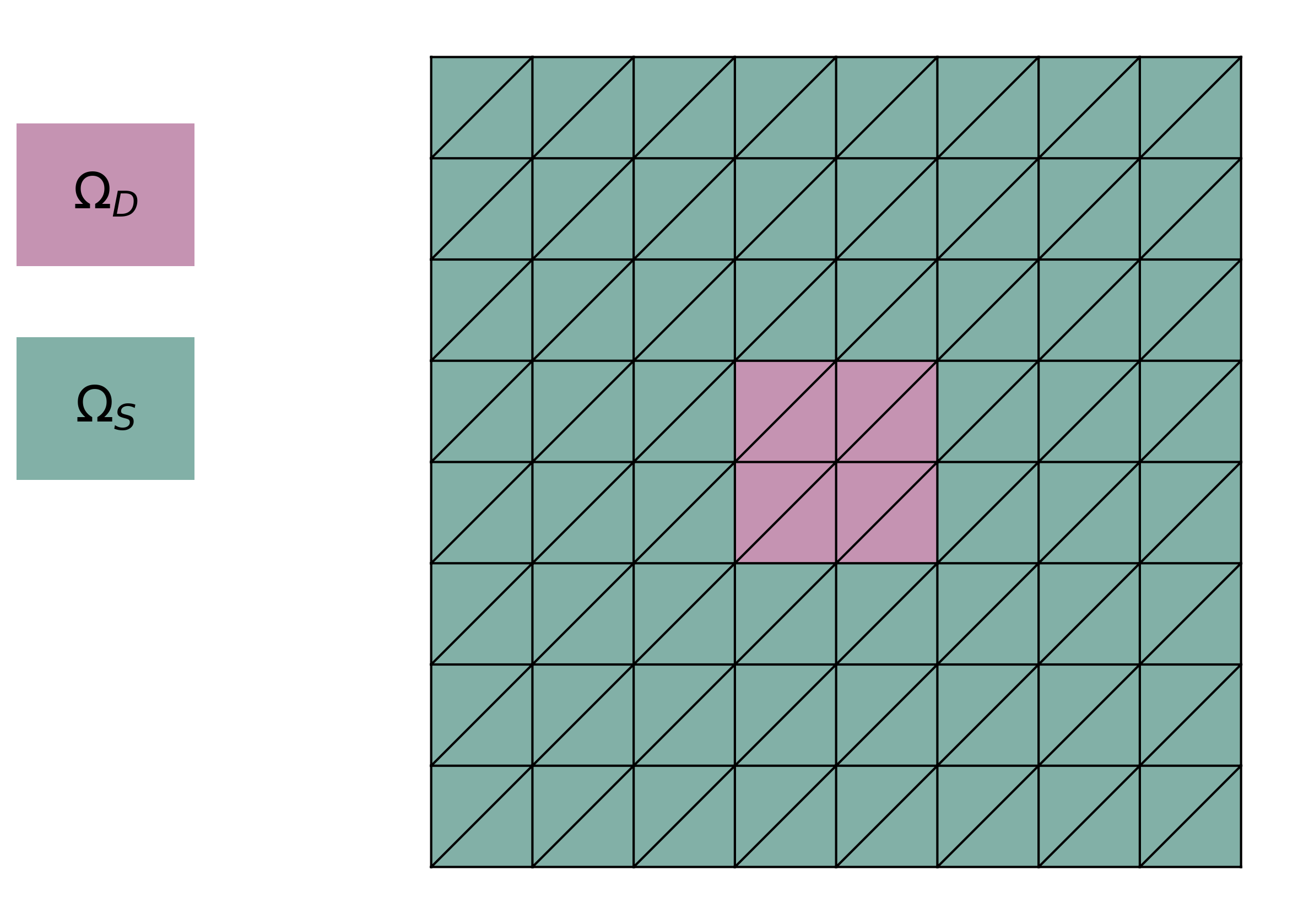}
	\caption{Sample of the computational domain for a mesh  with   $N=5$ (left) and $N=8$ (right).}
	\label{fig:malla}
\end{figure}

For the stabilization parameter, there is not a unique manner to define it. As discussed in \cite{MR2324460} and subsequently in \cite{MR3962898,MR4077220}, a common choice is $\texttt{a}_S=ak^2$, where $a>0$ and $k\geq 1$ is the polynomial degree of the approximation. For our tests, we aim to analyze other choices besides $\texttt{a}_S=ak^2$, in order to compare the behavior of the spectrum respect to the spurious values that may emerge. With this aim in mind, we introduce the following choices for the stabilization parameter in one hand $\texttt{a}_S=ak(k+1)$ and on the other $\texttt{a}_S= a(k+1)^2$.  For each of these choices of $\texttt{a}_S$, we need to define the values of $a>0$ for the experiments. On this case, we choose $a:=2^{i}/\kappa$, where $i\in\{-3,-2,-1,0,1,2,3\}$ and $\kappa>0$ is the parameter involved in \eqref{eq:permea_tests}. On the following results, we need to precise that the computational domain is the unit square $\Omega:=(0,1)^2$ with homogeneous Dirichlet boundary conditions, $\bu=\boldsymbol{0}$ on $\partial\Omega$. The computational domain is partitioned into a Stokes region $\Omega_S$ and a porous subregion $\Omega_D$, the latter defined as $\Omega_D:=(3/8,5/8)^2$, while $\Omega_S:=\Omega\setminus\Omega_D$.  A sketch of the domain and a sample of a mesh is presented in Figure \ref{fig:malla}. 

From now on, the parameter $N$ will represent the refinement level  it is related to the number of elements on each edge of the domain, as we observe from Figure \ref{fig:malla}. Let us remark that the relation between $N$ and the mesh-size is given by $N\sim h^{-1}$.

On the following tables, all the reported results were computed with $N=8$ and $\mathbf{K}^{-1}= 10^3\mathbf{I}$. 
Moreover,  the numbers enclosed in boxes represent spurious eigenvalues, identified by checking the plots of the associated functions.  We begin by reporting the results for $k=1$.

\begin{table}[!htbp]
\centering
\small
\setlength{\tabcolsep}{3pt}
\renewcommand{\arraystretch}{1.1}
\resizebox{\linewidth}{!}{%
\begin{tabular}{|c|c|c|c|c|c|c|}
\hline
$2^{-3}/\kappa$ & $2^{-2}/\kappa$ & $2^{-1}/\kappa$ & $1/\kappa$ & $2/\kappa$ & $4/\kappa$ & $8/\kappa$ \\
\hline
183.0044 & 169.5755 & 152.9317 & 135.0139 & 117.4854 & 101.6098 & 88.6569 \\
387.2616 & 372.8240 & 352.9674 & 328.5764 & 296.4376 & 260.4908 & 231.3170 \\
430.4875 & 407.4491 & 375.0036 & 336.1620 & 301.1152 & 272.9582 & 247.7343 \\
435.6670 & 418.8618 & 395.6796 & 367.5911 & 332.7328 & 295.4713 & 262.6055 \\
454.7370 & 433.5485 & 404.6547 & 370.1419 & 336.7113 & 305.5210 & 277.5423 \\
$\boxed{537.3598}$ & 506.7123 & 465.4440 & 419.2579 & 374.9141 & 336.1838 & 303.5393 \\
$\boxed{544.7538}$ & $\boxed{512.8245}$ & 471.8694 & 427.16517 & 382.6124 & 340.2777 & 304.8416 \\
$\boxed{549.9283}$ & $\boxed{527.8588}$ & $\boxed{495.1438}$ & $\boxed{452.9768}$ & 406.1780 & 361.0795 & 323.0318 \\
$\boxed{622.1696}$ & $\boxed{604.5873}$ & $\boxed{579.0504}$ & $\boxed{545.8858}$ & 489.0079 & $\boxed{437.0095}$ & 394.3032 \\
$\boxed{679.8253}$ & $\boxed{647.6487}$ & $\boxed{601.5022}$ & $\boxed{545.9849}$ & $\boxed{507.2284}$ & 455.2797 & 409.6292 \\
\hline
\end{tabular}%
}
\caption{Computed eigenvalues for  $\mathbf{K}^{-1}=10^{3}\mathbf{I}$, $k=1$, $N=8$, different values of $a$, and  $\texttt{a}_S=ak^2$.}
\label{tb:sp_ak2}
\end{table}

\begin{table}[!htbp]
\centering
\small
\setlength{\tabcolsep}{3pt}
\renewcommand{\arraystretch}{1.1}
\resizebox{\linewidth}{!}{%
\begin{tabular}{|c|c|c|c|c|c|c|}
\hline
$2^{-3}/\kappa$ & $2^{-2}/\kappa$ & $2^{-1}/\kappa$ & $1/\kappa$ & $2/\kappa$ & $4/\kappa$ & $8/\kappa$ \\
\hline
169.5755 & 152.9317 & 135.0139 & 117.4854 & 101.6098 & 88.6569 & 79.3912 \\
372.8240 & 352.9674 & 328.5764 & 296.4376 & 260.4908 & 231.3170 & 210.3561 \\
407.4491 & 375.0036 & 336.1620 & 301.1152 & 272.9582 & 247.7343 & 228.5324 \\
418.8618 & 395.6796 & 367.5911 & 332.7328 & 295.4713 & 262.6055 & 237.8065 \\
433.5485 & 404.6547 & 370.1419 & 336.7113 & 305.5210 & 277.5423 & 256.0172 \\
506.7123 & 465.4440 & 419.2579 & 374.9141 & 336.1838 & 303.5393 & 275.9383 \\
$\boxed{512.8245}$ & 471.8694 & 427.16517 & 382.6124 & 340.2777 & 304.8416 & 281.9340 \\
$\boxed{527.8588}$ & $\boxed{495.1438}$ & $\boxed{452.9768}$ & 406.1780 & 361.0795 & 323.0318 & 295.0521 \\
$\boxed{604.5873}$ & $\boxed{579.0504}$ & $\boxed{545.8858}$ & 489.0079 & $\boxed{437.0095}$ & 394.3032 & 363.2276 \\
$\boxed{647.6487}$ & $\boxed{601.5022}$ & $\boxed{545.9849}$ & $\boxed{507.2284}$ & 455.2797 & 409.6292 & 375.8176 \\
\hline
\end{tabular}%
}
\caption{Computed eigenvalues for  $\mathbf{K}^{-1}=10^{3}\mathbf{I}$, $k=1$, $N=8$, different values of $a$  and  $\texttt{a}_S=ak(k+1)$.}
\label{tb:sp_akk1}
\end{table}

\begin{table}[!htbp]
\centering
\small
\setlength{\tabcolsep}{3pt}
\renewcommand{\arraystretch}{1.1}
\resizebox{\linewidth}{!}{%
\begin{tabular}{|c|c|c|c|c|c|c|}
\hline
$2^{-3}/\kappa$ & $2^{-2}/\kappa$ & $2^{-1}/\kappa$ & $1/\kappa$ & $2/\kappa$ & $4/\kappa$ & $8/\kappa$ \\
\hline
152.9317 & 135.0139 & 117.4854 & 101.6098 & 88.6569 & 79.3912 & 73.5455 \\
352.9674 & 328.5764 & 296.4376 & 260.4908 & 231.3170 & 210.3561 & 197.0156 \\
375.0036 & 336.1620 & 301.1152 & 272.9582 & 247.7343 & 228.5324 & 215.9587 \\
395.6796 & 367.5911 & 332.7328 & 295.4713 & 262.6055 & 237.8065 & 221.6693 \\
404.6547 & 370.1419 & 336.7113 & 305.5210 & 277.5423 & 256.0172 & 241.7609 \\
465.4440 & 419.2579 & 374.9141 & 336.1838 & 303.5393 & 275.9383 & 258.0235 \\
471.8694 & 427.16517 & 382.6124 & 340.2777 & 304.8416 & 281.9340 & 267.0460 \\
$\boxed{495.1438}$ & $\boxed{452.9768}$ & 406.1780 & 361.0795 & 323.0318 & 295.0521 & 277.0053 \\
$\boxed{579.0504}$ & $\boxed{545.8858}$ & 489.0079 & $\boxed{437.0095}$ & 394.3032 & 363.2276 & 343.23766 \\
$\boxed{601.5022}$ & $\boxed{545.9849}$ & $\boxed{507.2284}$ & 455.2797 & 409.6292 & 375.8176 & 353.86167 \\
\hline
\end{tabular}%
}
\caption{Computed eigenvalues for  $\mathbf{K}^{-1}=10^{3}\mathbf{I}$, $k=1$, $N=8$, different values of $a$, and $\texttt{a}_S=a(k+1)^2$.}
\label{tb:sp_ak12}
\end{table}

From Tables \ref{tb:sp_ak2},\ref{tb:sp_akk1} and \ref{tb:sp_ak12} we observe that when $k=1$ there is presence of spurious eigenvalues for each definition of $\texttt{a}_S$. However, some differences are observed. First,  the amount of spurious eigenvalues when $\texttt{a}_S=ak^2$ is larger when we compare the results for $\texttt{a}_S\in\{ak(k+1),a(k+1)^2\}$. In fact, the spurious eigenvalues on Tables \ref{tb:sp_akk1} and \ref{tb:sp_ak12} are observable until $a=2/\kappa$ and $a=1/\kappa$, respectively, whereas if $\texttt{a}_S=ak^2$, the pollution remains until $a=4/\kappa$. This allows to conclude that the choice $\texttt{a}_S=a(k+1)^2$ is the one that best avoids spurious values. It is important to remark that there are some eigenvalues that may look as spurious but they are not. For example, in Table \ref{tb:sp_ak2}, when $a=2^{-3}/\kappa$, the first eigenvalue is $\lambda_{1,h}=183.0044$, whereas from \cite{MR5023159} we have that an extrapolated value is $\lambda_{1,h}=65.3658$. However, when we check the eigenfunctions, is observable that our first eigenvalue has physical  shape, identical to the obtained with the method of \cite{MR5023159} ,  but corresponds to a poor approximation of the expected one that is improving while the value of  $a$ increases. On the other tables, the phenomenon is similar.  Now we present the results when $k=2$.

\begin{table}[!htbp]
	\centering
	\begin{tabular}{|c|c|c|c|c|c|c|}
		\hline
		$2^{-3}/\kappa$ &$2^{-2}/\kappa$ & $2^{-1}/\kappa$ & $1/\kappa$ & $2/\kappa$ & $4/\kappa$ & $8/\kappa$\\
		\hline
		71.3356 & 69.9041 & 68.5097 & 67.3345 & 66.4736 & 65.9187 & 65.5954 \\
		
		195.8561 & 189.3534 & 183.1184 & 177.9085 & 174.0955 & 171.6300 & 170.1877 \\
		
		206.7639 & 201.5184 & 196.2738 & 191.7545 & 188.3936 & 186.2071 & 184.9262 \\
		
		219.4305 & 211.0060 & 202.9111 & 196.0763 & 191.0186 & 187.7242 & 185.7896 \\
		
		241.6589 & 233.9884 & 226.1142 & 219.1277 & 213.8084 & 210.2959 & 208.2211 \\
		
		249.6581 & 240.5175 & 231.7092 & 224.2521 & 218.7454 & 215.1747 & 213.0868 \\
		
		267.9827 & 258.2654 & 248.7457 & 240.5618 & 234.4415 & 230.4362 & 228.0801 \\
		
		273.7143 & 262.6006 & 251.7884 & 242.6016 & 235.7984 & 231.3745 & 228.7816 \\
		
		343.4864 & 327.8363 & 312.8303 & 300.1929 & 290.8733 & 284.8200 & 281.2717 \\
		
		347.2792 & 330.3779 & 314.4936 & 301.3281 & 291.7314 & 285.5458 & 281.9365 \\
		\hline
	\end{tabular}
	\caption{Computed eigenvalues for  $\mathbf{K}^{-1}=10^{3}\mathbf{I}$, $k=2$, $N=8$, different values of $a$, and $\texttt{a}_S=ak^2$.}
	\label{tb:sp2_ak2}
\end{table}

\begin{table}[!htbp]
	\centering
	\begin{tabular}{|c|c|c|c|c|c|c|}
		\hline
		$2^{-3}/\kappa$ &$2^{-2}/\kappa$ & $2^{-1}/\kappa$ & $1/\kappa$ & $2/\kappa$ & $4/\kappa$ & $8/\kappa$\\
		\hline
		70.5059 & 69.0707 & 67.7871 & 66.7912 & 66.1159 & 65.7072 & 65.4791 \\
		
		192.0728 & 185.6177 & 179.9126 & 175.5029 & 172.5071 & 170.6870 & 169.6680\\
		
		203.7399 & 198.4003 & 193.5051 & 189.6379 & 186.9857 & 185.3697 & 184.4645\\
		
		214.5253 & 206.1638 & 198.7153 & 192.8910 & 188.8982 & 186.4599 & 185.0913\\
		
		237.2580 & 229.3355 & 221.8569 & 215.7895 & 211.5510 & 208.9408 & 207.4705\\
		
		244.3358 & 235.2534 & 227.1314 & 220.7813 & 216.4452 & 213.8094 & 212.3351\\
		
		262.3443 & 252.5948 & 243.7356 & 236.7121 & 231.8647 & 228.8966 & 227.2292\\
		
		267.2633 & 256.1454 & 246.1518 & 238.3160 & 232.9499 & 229.6794 & 227.8469\\
		
		334.3688 & 318.8565 & 305.0674 & 294.3198 & 286.9755 & 282.5006 & 279.9921\\
		
		337.3877 & 320.8379 & 306.3847 & 295.2698 & 287.7443 & 283.1851 & 280.6377\\
		\hline
	\end{tabular}
	\caption{Computed eigenvalues for  $\mathbf{K}^{-1}=10^{3}\mathbf{I}$, $k=2$, $N=8$, different values of $a$, and $\texttt{a}_S=ak(k+1)$.}
	\label{tb:sp2_akk1}
\end{table}

\begin{table}[!htbp]
	\centering
	\begin{tabular}{|c|c|c|c|c|c|c|}
		\hline
		$2^{-3}/\kappa$ &$2^{-2}/\kappa$ & $2^{-1}/\kappa$ & $1/\kappa$ & $2/\kappa$ & $4/\kappa$ & $8/\kappa$\\
		\hline
		69.6587 & 68.2906 & 67.1651 & 66.3593 & 65.8499 & 65.5571 & 65.3990 \\
		
		188.2504 & 182.1451 & 177.1584 & 173.5883 & 171.3234 & 170.0168 & 169.3093\\
		
		200.6055 & 195.4378 & 191.0961 & 187.9443 & 185.9349 & 184.7744 & 184.1459\\
		
		209.5778 & 201.6404 & 195.0851 & 190.3423 & 187.3134 & 185.5601 & 184.6091\\
		
		232.6344 & 224.8364 & 218.0937 & 213.0898 & 209.8560 & 207.9745 & 206.9518\\
		
		238.9664 & 230.3232 & 223.1713 & 218.0110 & 214.7308 & 212.8396 & 211.8166\\
		
		256.6002 & 247.2334 & 239.3659 & 233.6201 & 229.9361 & 227.8004 & 226.6415\\
		
		260.7011 & 250.0828 & 241.2680 & 234.8895 & 230.8236 & 228.4743 & 227.2018\\
		
		325.1868 & 310.4778 & 298.3641 & 289.6295 & 284.0662 & 280.8510 & 279.1088\\
		
		327.5519 & 312.0288 & 299.4377 & 290.4574 & 284.7780 & 281.5093 & 279.7421\\
		\hline
	\end{tabular}
	\caption{Computed eigenvalues for  $\mathbf{K}^{-1}=10^{3}\mathbf{I}$, $k=2$, $N=8$, different values of $a$, and $\texttt{a}_S=a(k+1)^2$.}
	\label{tb:sp2_ak12}
\end{table}

Tables \ref{tb:sp2_ak2}, \ref{tb:sp2_akk1} and \ref{tb:sp2_ak12} reveal that when the polynomial degree increases, the computed spectrum does not present spurious eigenvalues, which is different compare with the lowest order case. In fact, for each of the choices of $\texttt{a}_S$ and each of the values of $a$, the method is robust for $k=2$. Moreover, each of the eigenvalues are well approximations, particularly the first one, contrary to what is observed from the tables for $k=1$. For completeness, we also present 
results for $k=3$ which are similar to the ones presented for $k=2$. Hence, and to avoid a large amount of tables, we only present in Table \ref{tb:sp3_ak2} results for $k=3$ and $\texttt{a}_S=a(k+1)^2$.

\begin{table}[!htbp]
	\centering
	\begin{tabular}{|c|c|c|c|c|c|c|}
		\hline
		$2^{-3}/\kappa$ &$2^{-2}/\kappa$ & $2^{-1}/\kappa$ & $1/\kappa$ & $2/\kappa$ & $4/\kappa$ & $8/\kappa$\\
		\hline
		65.5298 & 65.4680 & 65.4198 & 65.3875 & 65.3682 & 65.3576 & 65.3502 \\
		
		168.8659 & 168.4890 & 168.1990 & 168.0068 & 167.8930 & 167.8304 & 167.7975\\
		
		183.4184 & 183.1639 & 182.9671 & 182.8360 & 182.7580 & 182.7151 & 182.6925\\
		
		184.2141 & 183.6944 & 183.2919 & 183.0242 & 182.8654 & 182.7781 & 182.7321\\
		
		205.5887 & 205.1928 & 204.8855 & 204.6803 & 204.5583 & 204.4911 & 204.4557\\
		
		210.9531 & 210.4253 & 210.0133 & 209.7378 & 209.5739 & 209.4836 & 209.4361\\
		
		225.3380 & 224.7898 & 224.3645 & 224.0809 & 223.9123 & 223.8195 & 223.7706\\
		
		225.5351 & 224.9179 & 224.4435 & 224.1297 & 223.9443 & 223.8425 & 223.7890\\
		
		277.2956 & 276.0535 & 275.0971 & 274.4642 & 274.0904 & 273.8852 & 273.7775\\
		
		277.4253 & 276.2692 & 275.3811 & 274.7947 & 274.4488 & 274.2592 & 274.1597\\
		\hline
	\end{tabular}
	\caption{Computed eigenvalues for  $\mathbf{K}^{-1}=10^{3}\mathbf{I}$, $k=3$, $N=8$, different values of $a$, and $\texttt{a}_S=a(k+1)^2$.}
	\label{tb:sp3_ak2}
\end{table}

Let us remark that all the previous results were computed by considering $a=2^{i}/\kappa$ with $i\in\{-3,-2,-1,0,1,2,3\}$. This choice depends on the coefficient $\kappa$ associated to the permeability and for the experiments, this parameter was $\kappa=10^3$ which is considerable large. However, if we choose $a$ in such a way that does not depend on $\kappa$, we have observed that the method is more robust with respect to the spurious eigenvalues, where in fact, the computed spectrum is cleaner. To illustrate these phenomena, the following table presents the computed eigenvalues using $a=2^{i}$ and $\texttt{a}_S=a k^2$.

\begin{table}[!htbp]
	\centering
	\begin{tabular}{|c|c|c|c|c|c|c|}
		\hline
		$2^{-3}$ &$2^{-2}$ & $2^{-1}$ & $1$ & $2$ & $4$ & $8$\\
		\hline
		68.4306 & 67.4687 & 66.9776 & 66.7295 & 66.6047 & 66.5422 & 66.5109\\
		
		185.2321 & 183.0023 & 181.8621 & 181.2854 & 180.9954 & 180.8499 & 180.7771\\
		
		204.6880 & 202.5421 & 201.4434 & 200.8874 & 200.6077 & 200.4674 & 200.3972\\
		
		207.2811 & 204.5506 & 203.1538 & 202.4472 & 202.0919 & 201.9137 & 201.8245\\
		
		228.8493 & 226.3752 & 225.1064 & 224.4638 & 224.1404 & 223.9782 & 223.8970\\
		
		242.1016 & 239.0878 & 237.5473 & 236.7683 & 236.3767 & 236.1803 & 236.0820\\
		
		253.6358 & 251.0673 & 249.7498 & 249.0825 & 248.7466 & 248.5781 & 248.4937\\
		
		260.9349 & 257.8823 & 256.3199 & 255.5295 & 255.1319 & 254.9325 & 254.8326\\
		
		325.4272 & 322.0400 & 320.3057 & 319.4280 & 318.9865 & 318.7651 & 318.6542\\
		
		334.2269 & 330.4905 & 328.5776 & 327.6097 & 327.1228 & 326.8786 & 326.7563\\
		\hline
	\end{tabular}
	\caption{Computed eigenvalues for  $\mathbf{K}^{-1}=10^{3}\mathbf{I}$, $k=1$, $N=8$, different values of $a$, and  $a_S=ak^2$ with $a=2^{i}$ and $i\in\{-3,-2,-1,0, 1,2,3\}$. }
	\label{tb:sp1_ak2_a_sin_kappa}
\end{table}

When we compare the results of Table \ref{tb:sp_ak2} where the stabilization depends on $\kappa$ with Table \ref{tb:sp1_ak2_a_sin_kappa}, we observe that the choice of $\texttt{a}_S$ with no dependence on $\kappa$ leads to a method that does not introduce spurious eigenvalues for any value of $a$. This behavior has been also observed for other definitions of $\texttt{a}_S$ and polynomial degrees and  since  the results are similar, we do not incorporate on the manuscript.

Now, we analyze the stability in the computation of the spectrum for a large number of eigenvalues using the method. Here, the aim is to compare the choices of $a$ depending on whether or not this parameter relies on $\kappa$. To do this task, we compute 40 eigenvalues and consider $a_i = 8^{-1}i$  and $a_i = (8\kappa)^{-1}i$, with $i\in\{1,\dots,40\}$,  as values involved on the stabilization parameter $\texttt{a}_S$. The following plots report the results when $\mathbf{K}^{-1}=10^3\mathbf{I}$ and for $k=1,2,3$. For simplicity we report the results for the stabilization  $\texttt{a}_S=a_i k^2$.

%
%

\begin{figure}[!htbp]
	\centering
	\begin{subfigure}{0.4\textwidth}
		\centering
		\includegraphics[width=\linewidth]{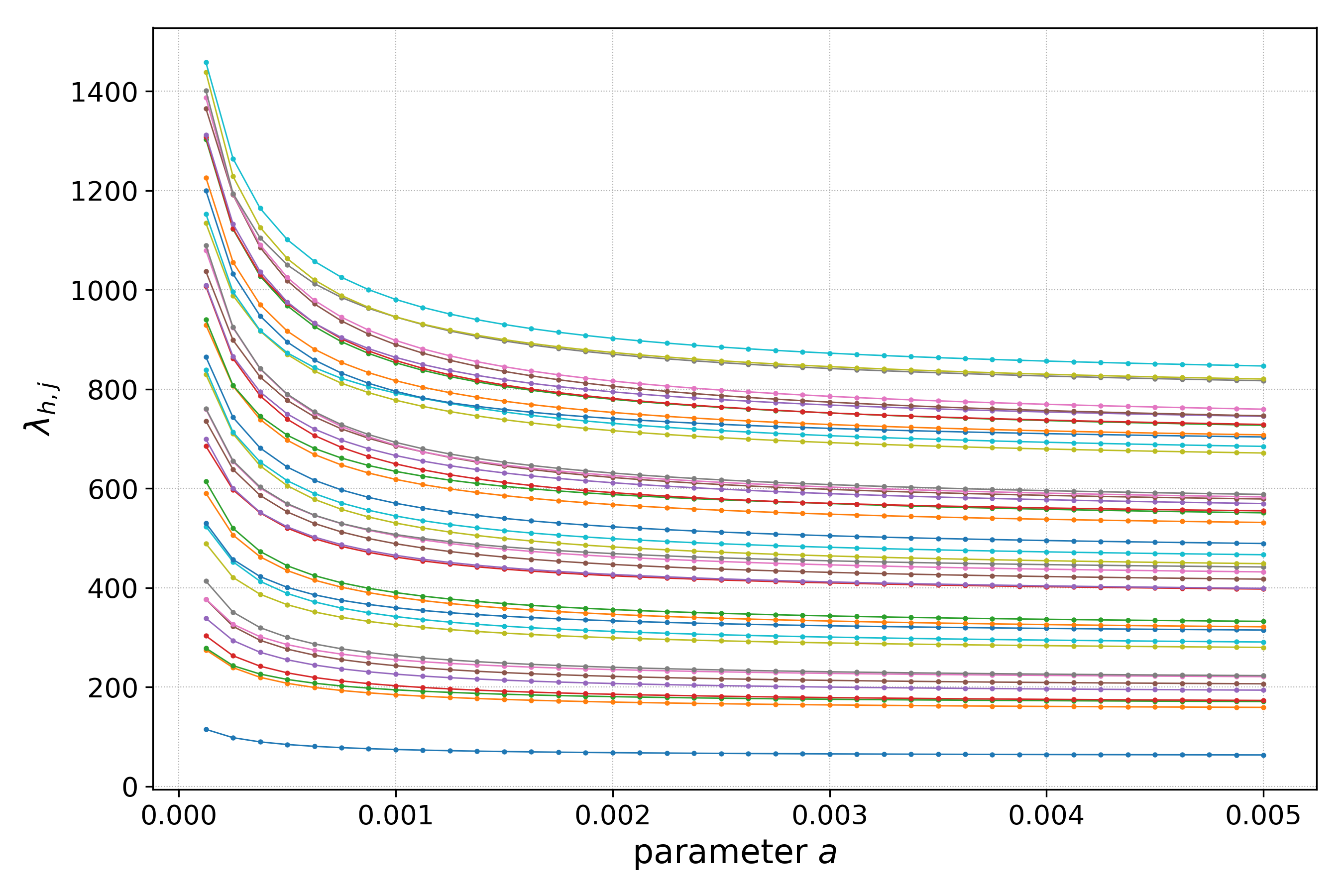}
		\caption{$a_i=(8\kappa)^{-1}i$}
		\label{fig:boffi_ak2_k1_a}
	\end{subfigure}
	\hfill
	\begin{subfigure}{0.4\textwidth}
		\centering
		\includegraphics[width=\linewidth]{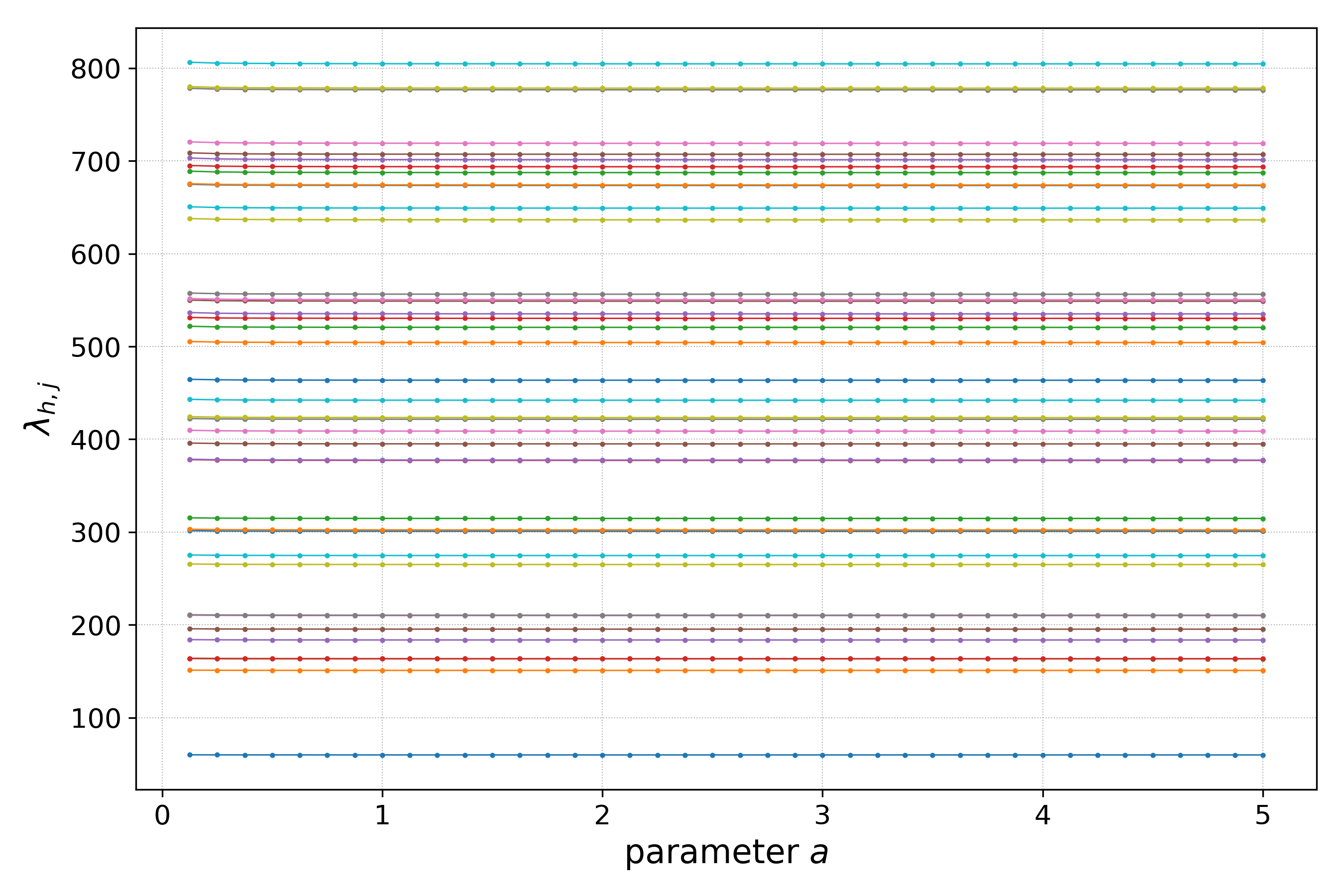}
		\caption{$a_i=8^{-1}i$}
		\label{fig:boffi_ak2_k1_b}
	\end{subfigure}
	\caption{Influence of  $a$ on the computation of the spectrum for ${\texttt{a}_S}=ak^2$ and $k=1$. Left:  40 eigenvalues computed when $a_i=(8\kappa)^{-1}i$. Right: 40 eigenvalues computed when $a_i=8^{-1}i$. }
	\label{fig:boffi_ak2_k1}
\end{figure}
\begin{figure}[!htbp]
	\centering
	\begin{subfigure}{0.4\textwidth}
		\centering
		\includegraphics[width=\linewidth]{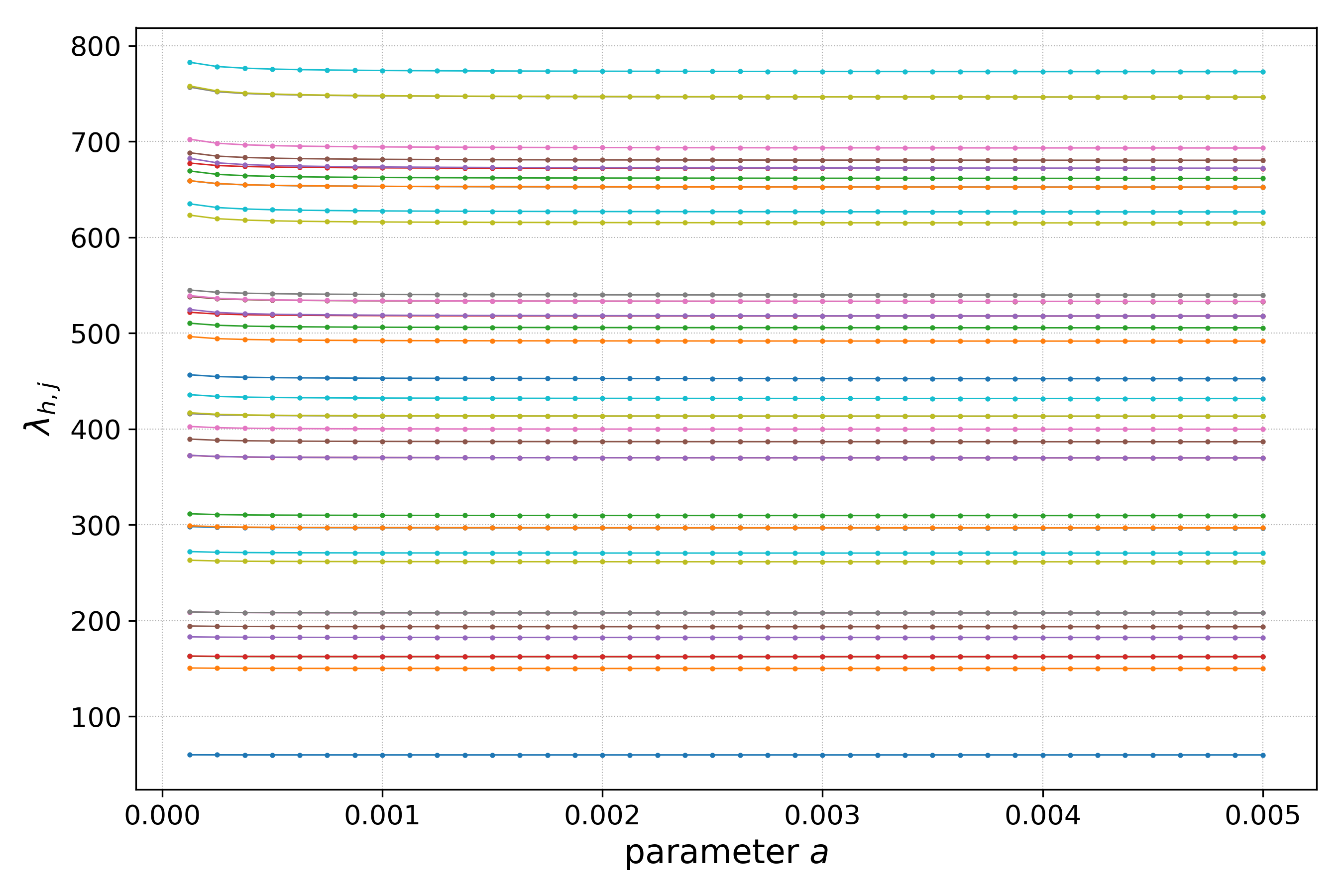}
		\caption{$a_i=(8\kappa)^{-1}i$}
		\label{fig:boffi_ak2_k2_a}
	\end{subfigure}
	\hfill
	\begin{subfigure}{0.4\textwidth}
		\centering
		\includegraphics[width=\linewidth]{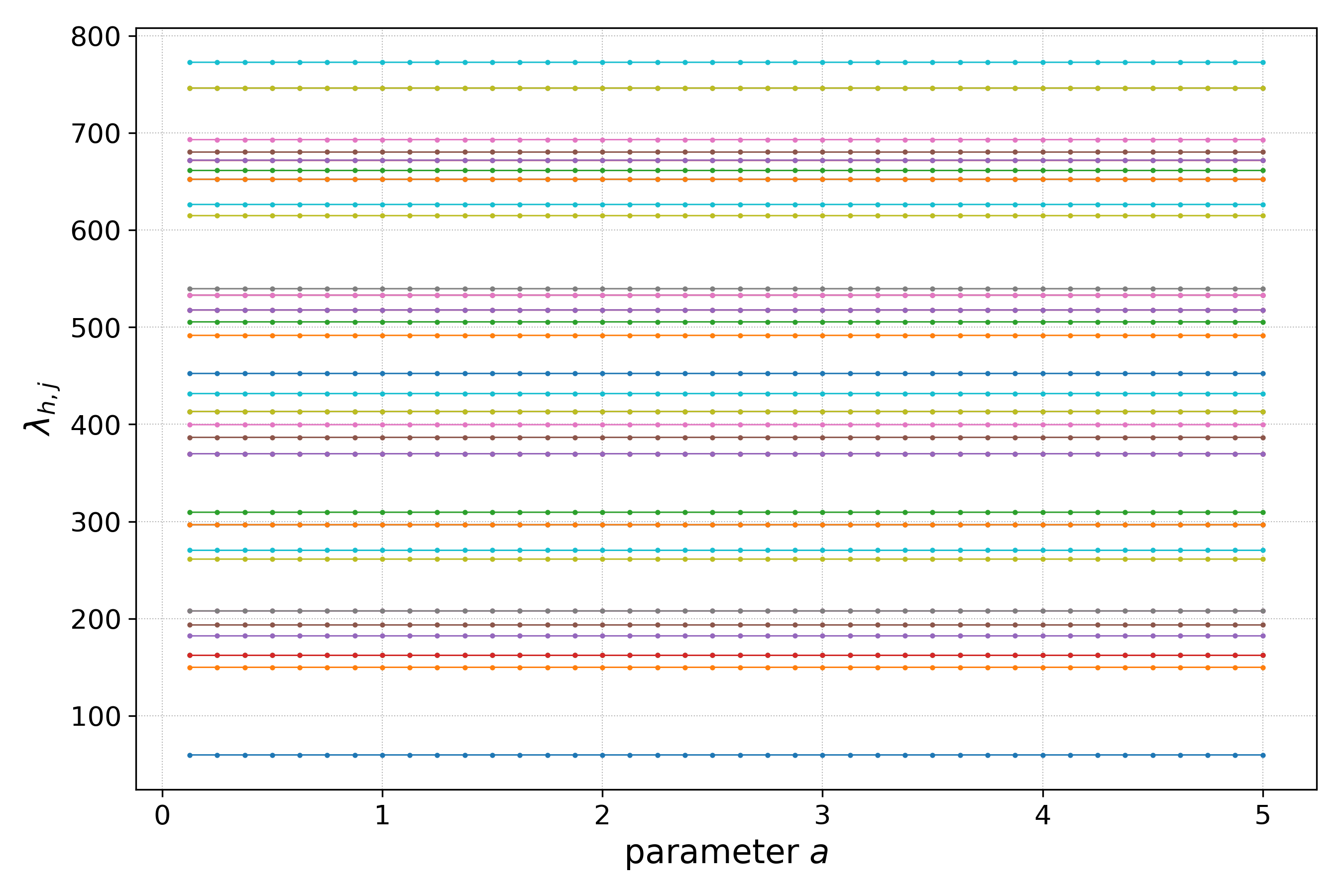}
		\caption{$a_i=8^{-1}i$}
		\label{fig:boffi_ak2_k2_b}
	\end{subfigure}
	\caption{Influence of  $a$ on the computation of the spectrum for ${\texttt{a}_S}=ak^2$ and $k=2$. Left:  40 eigenvalues computed when $a_i=(8\kappa)^{-1}i$. Right: 40 eigenvalues computed when $a_i=8^{-1}i$.}
	\label{fig:boffi_ak2_k2}
\end{figure}

\begin{figure}[!htbp]
	\centering
	\begin{subfigure}{0.4\textwidth}
		\centering
		\includegraphics[width=\linewidth]{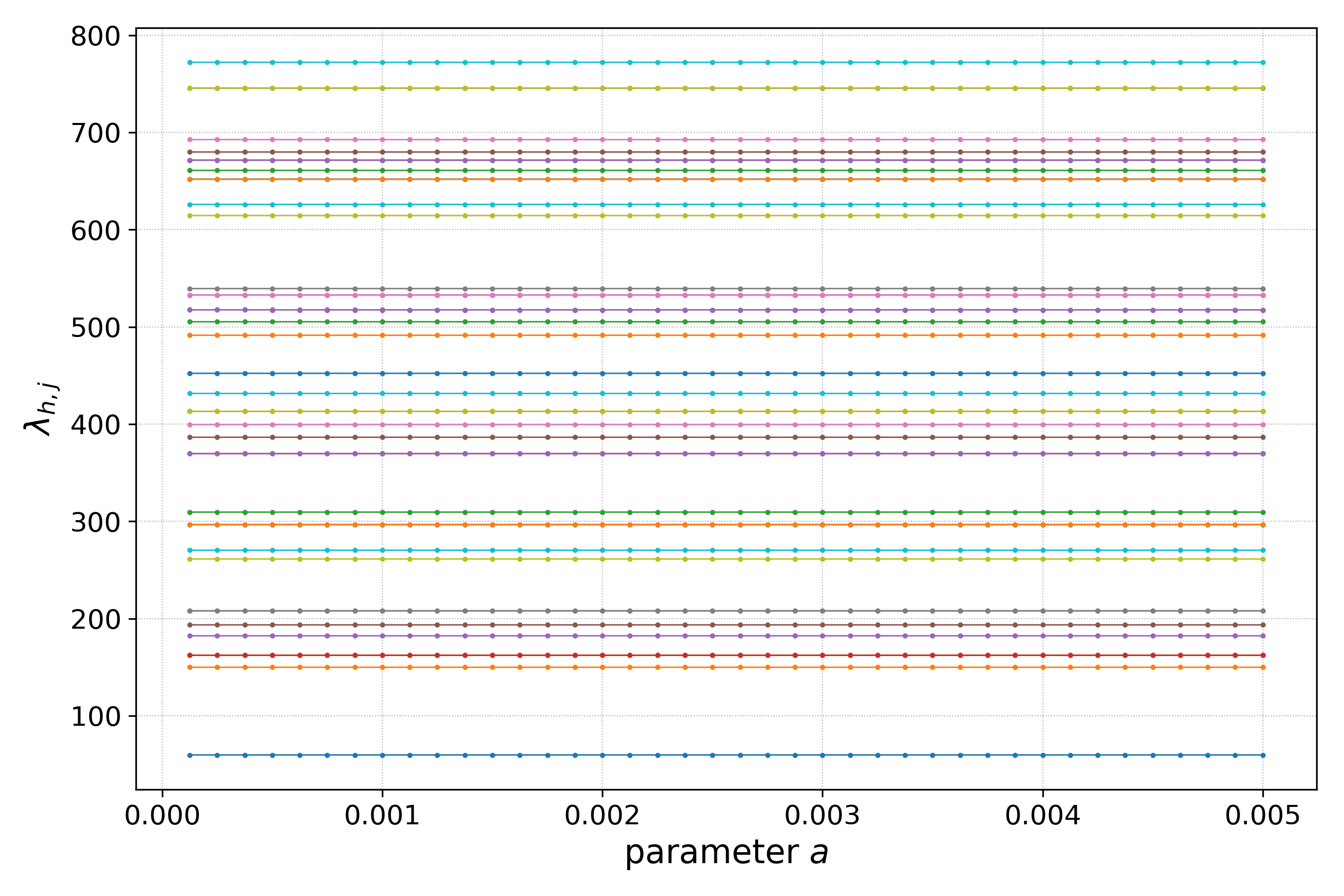}
		\caption{$a_i=(8\kappa)^{-1}i$}
		\label{fig:boffi_ak2_k3_a}
	\end{subfigure}
	\hfill
	\begin{subfigure}{0.4\textwidth}
		\centering
		\includegraphics[width=\linewidth]{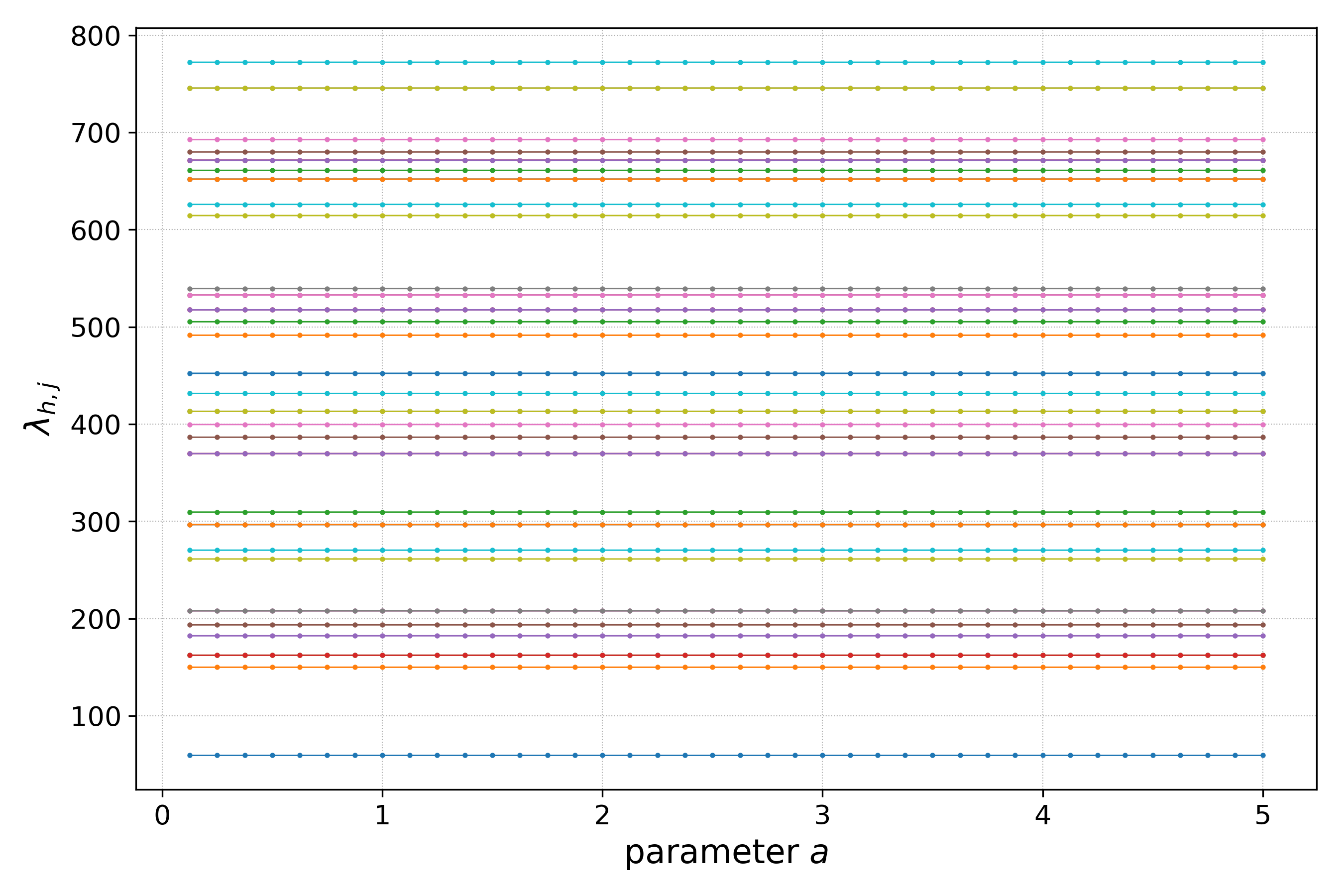}
		\caption{$a_i=8^{-1}i$}
		\label{fig:boffi_ak2_k3_b}
	\end{subfigure}
	\caption{Influence of  $a$ on the computation of the spectrum for ${\texttt{a}_S}=ak^2$ and $k=3$. Left:  40 eigenvalues computed when $a_i=(8\kappa)^{-1}i$. Right: 40 eigenvalues computed when $a_i=8^{-1}i$.}
	\label{fig:boffi_ak2_k3}
\end{figure}
\FloatBarrier

Figures \ref{fig:boffi_ak2_k1}, \ref{fig:boffi_ak2_k2}, and \ref{fig:boffi_ak2_k3} show the computation of 40 eigenvalues when $a_i=(8\kappa)^{-1}i$ and $a_i=8^{-1}i$. We observe from Figure \ref{fig:boffi_ak2_k1_a} that when we consider $a_i=(8\kappa)^{-1}i$, for $k=1$ the eigenvalues are well approximated when $a\approx 0.002$. On the other hand, when we compare the results for $a_i=8^{-1}i$ in Figure \ref{fig:boffi_ak2_k1_b}, this choice yields a more stable computation of the spectrum. However, when we increase the polynomial approximation, as we observe from Figures \ref{fig:boffi_ak2_k2_a} and \ref{fig:boffi_ak2_k3_a}, the method improves its capability to approximate the spectrum when $a_i=(8\kappa)^{-1}i$. Nevertheless,  $a_i=8^{-1}i$  still being  the best choice for the stabilization $\texttt{a}_S$, since Figures \ref{fig:boffi_ak2_k1_b}, \ref{fig:boffi_ak2_k2_b}, and \ref{fig:boffi_ak2_k3_b} reveal that the computed spectrum is stable for any choice of $a_i$. Hence, all our forthcoming experiments will be performed with stabilization values independent of $\kappa$.

\FloatBarrier
\subsection{Computation of rates of convergence}
The aim now is to compute experimental rates of convergence in order to confirm the theoretical result given in Theorem \ref{orddob}, To do this task, we consider the same configuration of the previous test, namely the unit square $\O:=(0,1)^2$, the inner obstacle $\O_D:=(3/8,5/8)^2$, null Dirichlet boundary conditions, and the same definition of $\mathbf{K}^{-1}$ given in \eqref{eq:permea_tests}. Since we have now a safe threshold of  $a$ for which the stabilizations $\texttt{a}_S$ does not introduce spurious eigenvalues, we will use the same as in the previous tests, specifically, values of $a$ that not depend on the parameter $\kappa$. 

To present the results in an organized manner, we first report the first four eigenvalues computed on a sufficiently refined mesh and compare these results with extrapolated values obtained via least-squares fitting. This fitting is such that, if $\lambda_h$ is a discrete eigenvalue, the order of convergence is calculated by fitting the following expression:
\begin{equation}
\label{eq:fitting}
\lambda_{h}\approx \lambda_{\text{extr}} + Ch^{\text{Order}},
\end{equation}
where $\lambda_{\text{extr}}$ is an extrapolated eigenvalue respect to the exact, $C>0$ is a constant associated to the error, and  "Order" represents the order of convergence of the method given by the fitting for the considered eigenvalue. This fitting leads to a system of equations where the unknowns are precisely $\lambda_{\text{extr}}$, constant $C$ and the order of convergence.   We remark that this is a classic non linear  fitting which is useful when the meshes are sufficiently refined. Let us remark that  all the reported  eigenvalues have been computed for different values of $\mathbf{K}^{-1}$.
\begin{table}[!htbp]
\centering
\scriptsize
\setlength{\tabcolsep}{4pt}
\renewcommand{\arraystretch}{1.10}
\caption{First four eigenvalues computed with a sufficiently refined mesh for different values of $\mathbf{K}^{-1}$. The results were computed for $k=1,2,3$ when $\texttt{a}_S=ak^2$ and $a=8$.}
\label{tab:eigen_all}
\begin{tabular}{@{}c|c|c|cccc@{}}
\toprule
$\mathbf K^{-1}$ & $k$ & Result
& $\lambda_1$ & $\lambda_2$ & $\lambda_3$ & $\lambda_4$ \\
\hline
\multirow{7}{*}{$10^{-8}\mathbf I$}
& \multirow{2}{*}{1}
& Finest mesh   & 52.4408 & 92.3778 & 92.4635 & 128.7809 \\
& & Extrap.& 52.3445 & 92.1215 & 92.1230 & 128.1973 \\
\cmidrule(lr){2-7}
& \multirow{2}{*}{2}
& Finest mesh   & 52.3448 & 92.1250 & 92.1253 & 128.2119 \\
& & Extrap. & 52.3446 & 92.1243 & 92.1243 & 128.2093 \\
\cmidrule(lr){2-7}
& \multirow{2}{*}{3}
& Finest mesh   & 52.3446 & 92.1243 & 92.1243 & 128.2095 \\
& & Extrap.& 52.3446 & 92.1243 & 92.1243 & 128.2095 \\
\cmidrule(lr){2-7}
&  & \cite{MR5023159}  & 52.3448 & 92.1243 & 92.1241 & 128.2085 \\
\hline
\multirow{6}{*}{$10^{-3}\mathbf I$}
& \multirow{2}{*}{1}
& Finest mesh   & 52.4948 & 92.5202 & 92.6542 & 129.1005 \\
& & Extrap. & 52.3395 & 92.1049 & 92.1037 & 128.1154 \\
\cmidrule(lr){2-7}
& \multirow{2}{*}{2}
& Finest mesh   & 52.3448 & 92.1250 & 92.1253 & 128.2119 \\
& & Extrap.& 52.3446 & 92.1243 & 92.1243 & 128.2093 \\
\cmidrule(lr){2-7}
& \multirow{2}{*}{3}
& Finest mesh   & 52.3447 & 92.1246 & 92.1246 & 128.2096 \\
& & Extrap.& 52.3447 & 92.1246 & 92.1246 & 128.2096 \\
\hline

\multirow{6}{*}{$10\mathbf I$}
& \multirow{2}{*}{1}
& Finest mesh   & 52.9621 & 95.0571 & 95.1903 & 129.8808 \\
& & Extrapolated & 52.8108 & 94.6316 & 94.6415 & 128.8998 \\
\cmidrule(lr){2-7}
& \multirow{2}{*}{2}
& Finest mesh   & 52.8153 & 94.6524 & 94.6528 & 128.9948 \\
& & Extrap.& 52.8152 & 94.6517 & 94.6517 & 128.9922 \\
\cmidrule(lr){2-7}
& \multirow{2}{*}{3}
& Finest mesh   & 52.8152 & 94.6517 & 94.6517 & 128.9925 \\
& & Extrap. & 52.8152 & 94.6518 & 94.6518 & 128.9924 \\
\midrule

\multirow{7}{*}{$10^{3}\mathbf I$}
& \multirow{2}{*}{1}
& Finest mesh   & 65.4125 & 168.4313 & 183.5071 & 183.5775 \\
& & Extrap. & 65.3662 & 167.6949 & 182.6188 & 182.6085 \\
\cmidrule(lr){2-7}
& \multirow{2}{*}{2}
& Finest mesh   & 65.3651 & 167.7495 & 182.6618 & 182.6622 \\
& & Extrap. & 65.3659 & 167.7479 & 182.6604 & 182.6604 \\
\cmidrule(lr){2-7}
& \multirow{2}{*}{3}
& Finest mesh   & 65.3649 & 167.7476 & 182.6598 & 182.6600 \\
& & Extrap. & 65.3642 & 167.7460 & 182.6583 & 182.6574 \\
\cmidrule(lr){2-7}
& & \cite{MR5023159} & 65.3633 & 167.7440 & 182.6572 & 182.6564 \\
\bottomrule
\end{tabular}
\end{table}
 The results reported in Table \ref{tab:eigen_all} show, in one hand, eigenvalues that are known from \cite{MR5023159}  when $\mathbf{K}^{-1}\in\{10^{-8}\mathbf{I}, 10^3\mathbf{I}\}$. Additionally, we present results when $\mathbf{K}^{-1}\in\{10^{-3}\mathbf{I},10\mathbf{I}\}$. We observe that our extrapolated values are closed when we compare them to the results on \cite{MR5023159}, whereas for the other cases, the eigenvalues when are computed for a sufficiently refined mesh are similar to the extrapolated ones. Let us remark that when $\texttt{a}_S\in\{ak(k+1),a(k+1)^2\}$ and $a=8$, the results are similar.

Next, we compute the experimental rates of convergence. These values are obtained using the fitting defined in \eqref{eq:fitting} for each of the three choices of $\texttt{a}_S$ with $a=8$ fixed.

\begin{table}[H]
\centering
\caption{Experimental convergence orders for
$k=1,2,3$, considering three choices of the
stabilization parameter $a_S$ with $a=8$.}
\label{tab:orders_all}


\begin{subtable}[t]{0.48\textwidth}
\centering
\caption{$k=1$}
\label{tab:orders_k1}
\resizebox{\linewidth}{!}{%
\begin{tabular}{@{}c c cccc@{}}
	\toprule
	$\mathbf K^{-1}$ & $a_S$
	& $\lambda_1$ & $\lambda_2$
	& $\lambda_3$ & $\lambda_4$\\
	\midrule
	& $8k^2$       & 1.99 & 1.98 & 1.99 & 1.96\\
	$10^{-8}\mathbf I$ & $8k(k+1)$ & 1.99 & 1.98 & 1.99 & 1.96\\
	& $8(k+1)^2$   & 1.99 & 1.98 & 1.99 & 1.96\\

\midrule
& $8k^2$       & 1.96 & 1.94 & 1.96 & 1.88\\
$10^{-3}\mathbf I$& $8k(k+1)$    & 1.96 & 1.94 & 1.96 & 1.88\\
& $8(k+1)^2$   & 1.96 & 1.94 & 1.96 & 1.88\\
\midrule
& $8k^2$       & 1.96 & 1.94 & 1.97 & 1.88\\
$10\mathbf I$& $8k(k+1)$    & 1.96 & 1.94 & 1.97 & 1.88\\
& $8(k+1)^2$   & 1.96 & 1.94 & 1.97 & 1.88\\

\midrule
& $8k^2$       & 2.00 & 1.86 & 1.91 & 1.90\\
$10^3\mathbf I$& $8k(k+1)$    & 2.00 & 1.86 & 1.91 & 1.90\\
& $8(k+1)^2$   & 2.00 & 1.86 & 1.91 & 1.90\\
\bottomrule
\end{tabular}%
}
\end{subtable}%
\hfill
\begin{subtable}[t]{0.48\textwidth}
\centering
\caption{$k=2$}
\label{tab:orders_k2}
\resizebox{\linewidth}{!}{%
\begin{tabular}{@{}c c cccc@{}}
\toprule
$\mathbf K^{-1}$ & $a_S$
& $\lambda_1$ & $\lambda_2$
& $\lambda_3$ & $\lambda_4$\\
\midrule

& $8k^2$       & 3.91 & 3.90 & 3.90 & 3.84\\
$10^{-8}\mathbf I$& $8k(k+1)$    & 3.90 & 3.89 & 3.90 & 3.84\\
& $8(k+1)^2$   & 3.57 & 3.99 & 3.95 & 3.20\\

\midrule
& $8k^2$       & 3.91 & 3.90 & 3.90 & 3.84\\
$10^{-3}\mathbf I$& $8k(k+1)$    & 3.91 & 3.90 & 3.90 & 3.84\\
& $8(k+1)^2$   & 3.90 & 3.90 & 3.90 & 3.84\\

\midrule
& $8k^2$       & 3.91 & 3.89 & 3.90 & 3.84\\
$10\mathbf I$& $8k(k+1)$    & 3.91 & 3.89 & 3.90 & 3.84\\
& $8(k+1)^2$   & 3.91 & 3.89 & 3.89 & 3.83\\

\midrule
& $8k^2$       & 3.46 & 3.93 & 3.99 & 3.96\\
$10^3\mathbf I$& $8k(k+1)$    & 3.46 & 3.93 & 3.99 & 3.95\\
& $8(k+1)^2$   & 3.45 & 3.94 & 4.00 & 3.96\\
\bottomrule
\end{tabular}%
}
\end{subtable}

\vspace{1.2em}


\begin{subtable}[t]{0.48\textwidth}
\centering
\caption{$k=3$}
\label{tab:orders_k3}
\resizebox{\linewidth}{!}{%
\begin{tabular}{@{}c c cccc@{}}
\toprule
$\mathbf K^{-1}$ & $a_S$
& $\lambda_1$ & $\lambda_2$
& $\lambda_3$ & $\lambda_4$\\
\midrule

& $8k^2$       & 5.97 & 5.99 & 5.97 & 5.88\\
$10^{-8}\mathbf I$& $8k(k+1)$    & 5.96 & 5.93 & 5.96 & 5.86\\
& $8(k+1)^2$   & 5.90 & 6.00 & 5.95 & 5.87\\

\midrule
& $8k^2$       & 5.97 & 5.99 & 5.97 & 5.88\\
$10^{-3}\mathbf I$& $8k(k+1)$    & 5.94 & 5.94 & 5.96 & 5.87\\
& $8(k+1)^2$   & 5.92 & 6.01 & 5.95 & 5.87\\

\midrule
& $8k^2$       & 6.05 & 6.05 & 6.03 & 5.90\\
$10\mathbf I$& $8k(k+1)$    & 6.00 & 5.98 & 6.02 & 5.88\\
& $8(k+1)^2$   & 5.99 & 6.06 & 6.00 & 5.89\\

\midrule
& $8k^2$       & 5.74 & 6.06 & 6.25 & 6.00\\
$10^3\mathbf I$& $8k(k+1)$    & 5.73 & 6.06 & 6.26 & 6.00\\
& $8(k+1)^2$   & 5.73 & 6.06 & 6.25 & 6.00\\
\bottomrule
\end{tabular}%
}
\end{subtable}

\end{table}

\begin{figure}[!htbp]
    \centering
    \begin{subfigure}{0.295\textwidth}
        \centering
        \includegraphics[width=\linewidth]{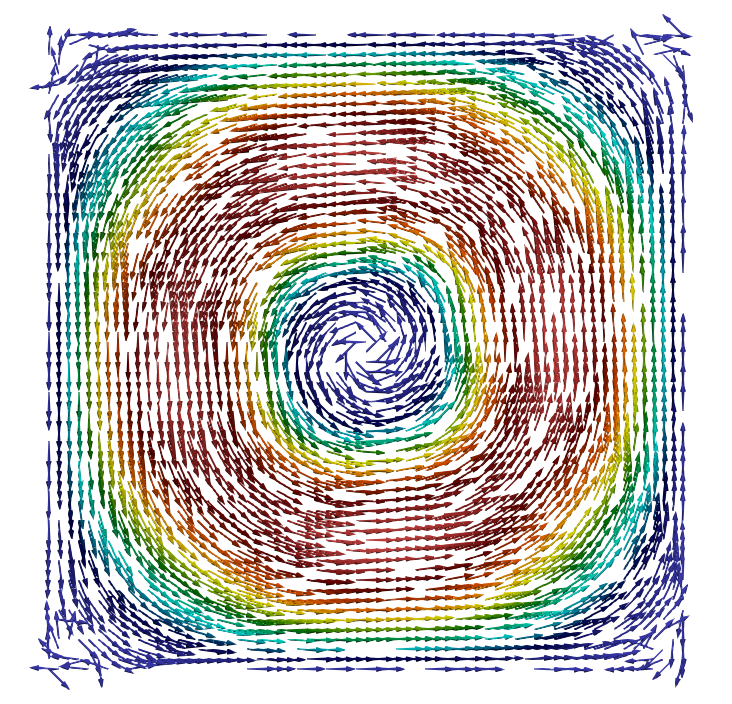}
    \end{subfigure}
    \hfill
    \begin{subfigure}{0.295\textwidth}
        \centering
        \includegraphics[width=\linewidth]{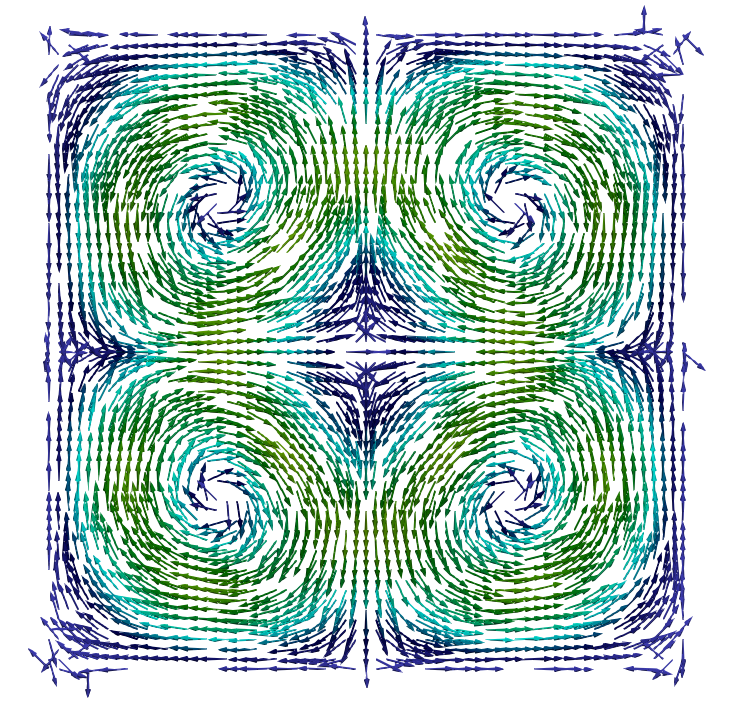}
    \end{subfigure}
    \hfill
    \begin{subfigure}{0.295\textwidth}
        \centering
        \includegraphics[width=\linewidth]{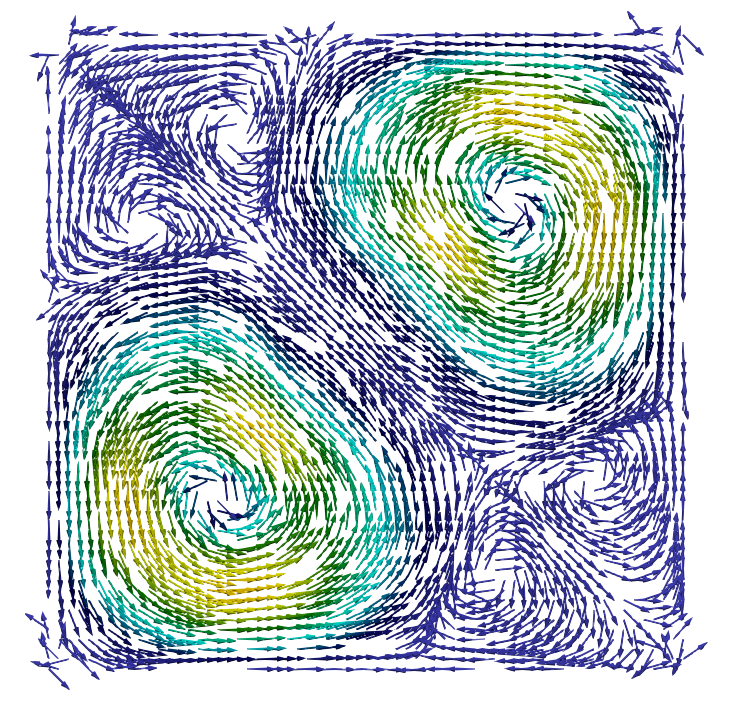}
    \end{subfigure}

    \par\medskip

    \begin{subfigure}{0.295\textwidth}
        \centering
        \includegraphics[width=\linewidth]{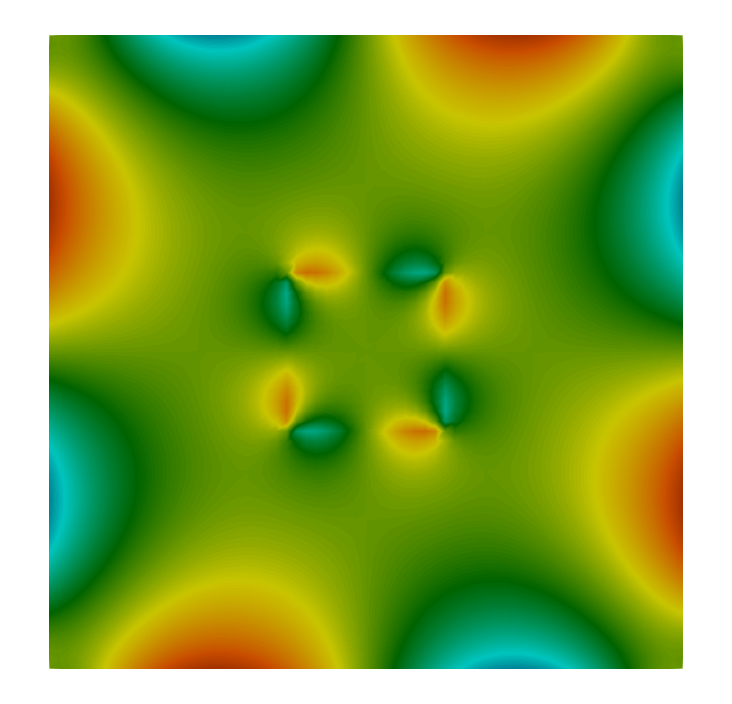}
    \end{subfigure}
    \hfill
    \begin{subfigure}{0.295\textwidth}
        \centering
        \includegraphics[width=\linewidth]{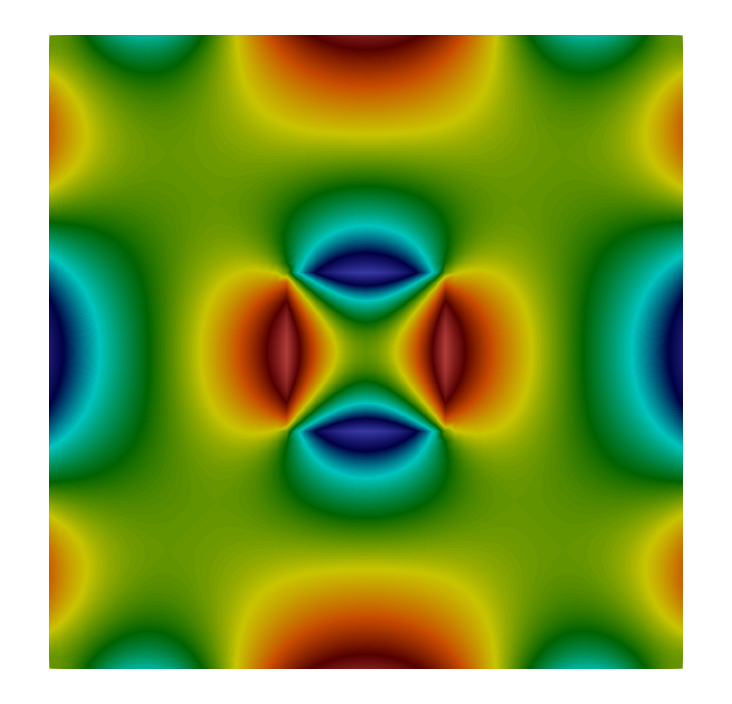}
    \end{subfigure}
    \hfill
    \begin{subfigure}{0.295\textwidth}
        \centering
        \includegraphics[width=\linewidth]{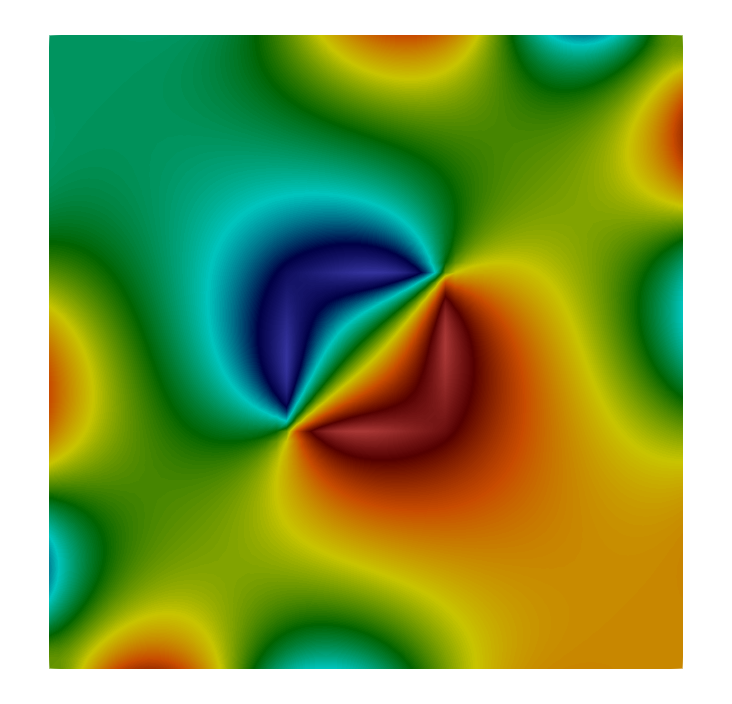}
    \end{subfigure}

    \caption{Plots of some computed fields. Top row: velocity fields associated to the first, second, and fourth eigenfunctions, respectively.
    Bottom row: pressure fluctuations associated to the first, second, and fourth eigenfunctions, respectively.
   These plots have been computed when 
    $\mathbf{K}^{-1}=10^{3}\mathbf{I}$,
    $k=1$ and  ${\texttt{a}_S}=ak^2$.}
  \label{fig:vel_pre_2}
\end{figure}

\begin{figure}[!htbp]
    \centering

    \begin{subfigure}{0.28\textwidth}
        \centering
        \includegraphics[width=\linewidth]{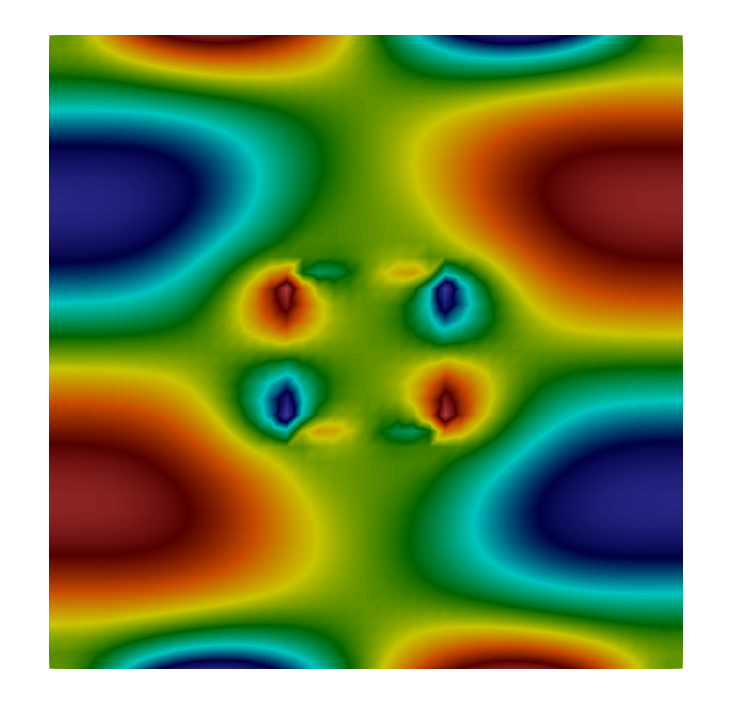}
        \caption{$\bsig_{h,(1,1)}$}
    \end{subfigure}
    \hspace{0.06\textwidth}
    \begin{subfigure}{0.28\textwidth}
        \centering
        \includegraphics[width=\linewidth]{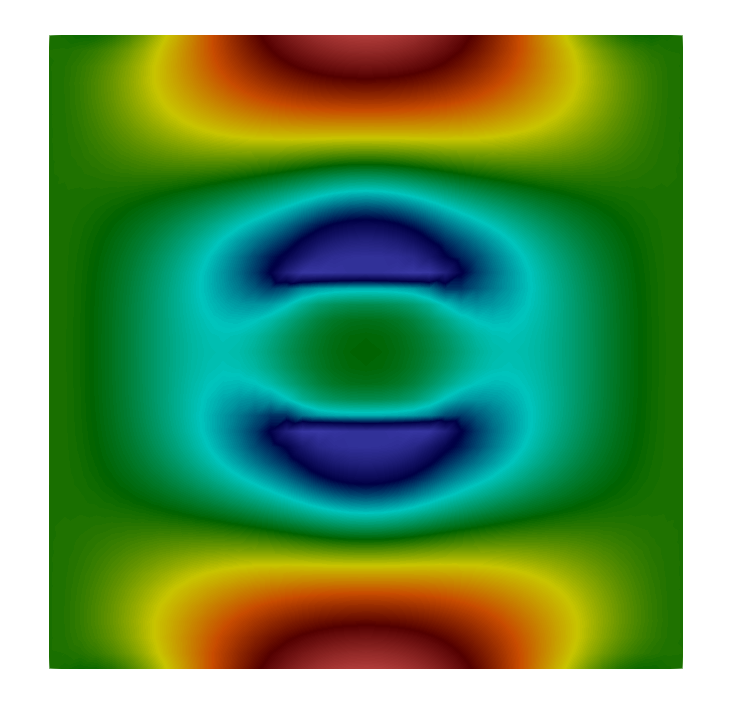}
        \caption{$\bsig_{h,(1,2)}$}
    \end{subfigure}

    \par\smallskip

    \begin{subfigure}{0.28\textwidth}
        \centering
        \includegraphics[width=\linewidth]{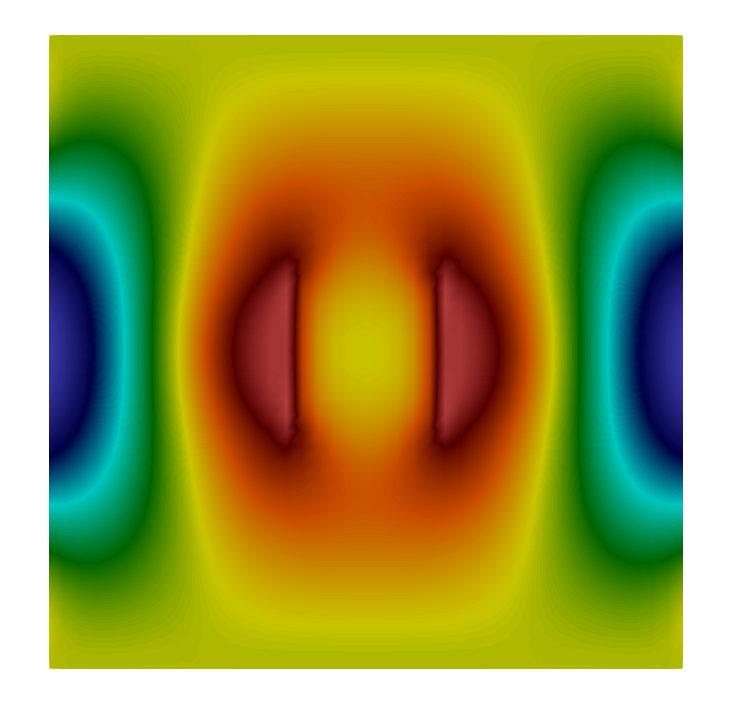}
        \caption{$\bsig_{h,(2,1)}$}
    \end{subfigure}
    \hspace{0.06\textwidth}
    \begin{subfigure}{0.28\textwidth}
        \centering
        \includegraphics[width=\linewidth]{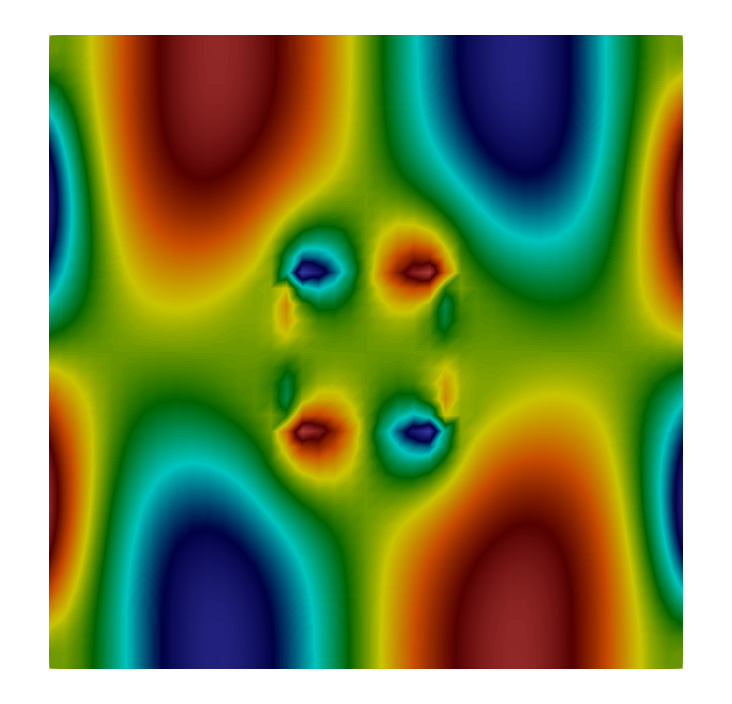}
        \caption{$\bsig_{h,(2,2)}$}
    \end{subfigure}

    \caption{Pseudostress components obtained for $\mathbf{K}^{-1}=10^{3}\mathbf{I}$,
    $k=1$ and ${\texttt{a}_S}=ak^2$.}

    \label{fig:comp_sig_103}
\end{figure}

Tables \ref{tab:orders_k1}, \ref{tab:orders_k2} and \ref{tab:orders_k3} show that the order of convergence for the eigenvalues is precisely $\mathcal{O}(h^{2\min\{r,k\}})$ as we expect when $k\in\{1,2,3\}$. Let us precise that when the value of $\mathbf{K}^{-1}$ increases, we observe that the order of convergence for the first eigenvalue shows some differences when $k>1$. Indeed, we expect this, since as the media is  more permeable, this involves the presence of an obstacle within the domain which affects the associated eigenfunctions and hence, this singularities that appear affects the regularity of the eigenfunctions which is reflected on the loss of the optimal order of convergence. This has been observed for the three stabilizations $\texttt{a}_S$ that we are considering. In Figure \ref{fig:vel_pre_2} we present plots of the computed  velocity and pressure fields whereas in Figure \ref{fig:comp_sig_103} the fourth components of the pseudostress are presented. All these plots have been obtained for $\mathbf{K}^{-1}=10^3\mathbf{I}$.

\FloatBarrier

\subsection{A three dimensional test}
As a final test, we compute experimental rates of convergence for a three dimensional domain. On this case, we consider the unit cube $\O:=(0,1)^3$ where within this domain, we define an additional cube defined by $\Omega_D=(3/8,5/8)^3$. For simplicity, we use $k=1$ as polynomial degree of approximation and $a=8$ for each of the selections of $\texttt{a}_S$ that we have used along our paper. In Table \ref{tb:cub_stabilization_comparison} we report the obtained results, first, for $\mathbf{K}^{-1}=10^{-8}\mathbf{I}$ which corresponds to a  Stokes regime. 
\begin{table}[!htbp]
\centering
\scriptsize
\setlength{\tabcolsep}{4pt}
\renewcommand{\arraystretch}{1.1}
\begin{tabular}{|c|c|c|c|c|c|c|}
\hline
$\texttt{a}_S$
& N
& $\lambda_{h,1}$
& $\lambda_{h,2}$
& $\lambda_{h,3}$
& $\lambda_{h,4}$
& $\lambda_{h,5}$ \\
\hline
\multirow{7}{*}{$ak^2$}
& 4 & 71.6619 & 73.0956 & 73.0956 & 109.4888 & 109.4888 \\
& 6 & 66.2997 & 66.9664 & 66.9664 & 100.1708 & 100.1708 \\
& 8 & 64.3842 & 64.8206 & 64.8206 & 96.7471 & 96.7471 \\
& 10 & 63.7444 & 64.0037 & 64.0037 & 95.0861 & 95.0861 \\
\cline{2-7}
& Order & 2.18 & 2.14 & 2.14 & 1.86 & 1.86 \\
\cline{2-7}
& Extrap. & 62.4407 & 62.4726 & 62.4726 & 91.8919 & 91.8919 \\
\cline{2-7}
&  \cite{MR5023159} & 62.1733 & 62.1733 & 62.1733 & 91.6298 & 91.6298 \\
\hline
\multirow{7}{*}{$ak(k+1)$}
& 4 & 71.5249 & 72.9424 & 72.9424 & 109.3409 & 109.3409 \\
& 6 & 66.2375 & 66.8937 & 66.8937 & 100.0893 & 100.0893 \\
& 8 & 64.3523 & 64.7819 & 64.7819 & 96.6989 & 96.6989 \\
& 10 & 63.7204 & 63.9757 & 63.9757 & 95.0513 & 95.0513 \\
\cline{2-7}
& Order & 2.18 & 2.14 & 2.14 & 1.86 & 1.86 \\
\cline{2-7}
& Extrap. & 62.4354 & 62.4649 & 62.4649 & 91.8803 & 91.8803 \\
\cline{2-7}
&  \cite{MR5023159} & 62.1733 & 62.1733 & 62.1733 & 91.6298 & 91.6298 \\
\hline
\multirow{7}{*}{$a(k+1)^2$}
& 4 & 71.4562 & 72.8655 & 72.8655 & 109.2666 & 109.2666 \\
& 6 & 66.2064 & 66.8574 & 66.8574 & 100.0486 & 100.0486 \\
& 8 & 64.3363 & 64.7626 & 64.7626 & 96.6747 & 96.6747 \\
& 10 & 63.7084 & 63.9616 & 63.9616 & 95.0339 & 95.0339 \\
\cline{2-7}
& Order & 2.18 & 2.15 & 2.15 & 1.86 & 1.86 \\
\cline{2-7}
& Extrap. & 62.4327 & 62.4810 & 62.4810 & 91.8745 & 91.8745 \\
\cline{2-7}
& \cite{MR5023159} & 62.1733 & 62.1733 & 62.1733 & 91.6298 & 91.6298 \\
\hline
\end{tabular}
\caption{First four eigenvalues for $k=1$,
$\mathbf{K}^{-1}=10^{-8}\mathbf{I}$ and $a=8$,
considering different choices of the stabilization
parameter $\texttt{a}_S$.}
\label{tb:cub_stabilization_comparison}
\end{table}

\begin{figure}[!htbp]
    \centering
    \begin{subfigure}{0.34\textwidth}
        \centering
        \includegraphics[width=\linewidth]{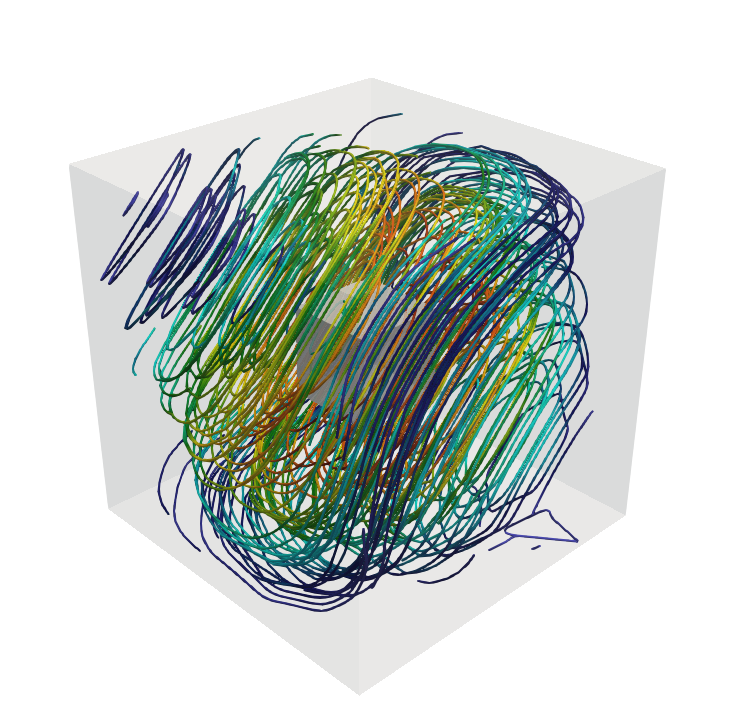}
    \end{subfigure}
    \hspace{0.06\textwidth}
  \begin{subfigure}{0.34\textwidth}
        \centering
        \includegraphics[width=\linewidth]{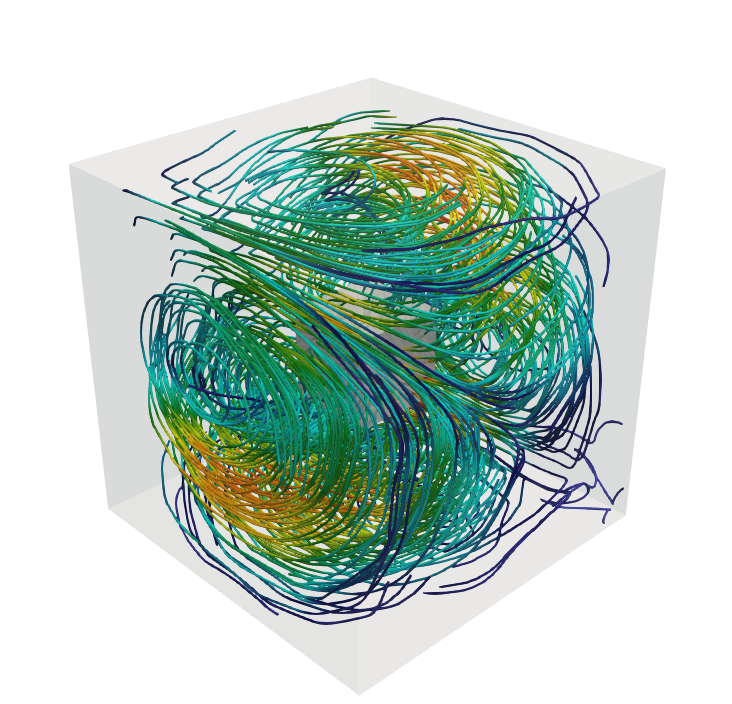}
    \end{subfigure}
    \caption{Representation of the velocity fields
    $\bu_{h,1}$ and $\bu_{h,4}$ for
    $\mathbf{K}^{-1}=10^{-8}\mathbf{I}$,
    $k=1$ and ${\texttt{a}_S}=ak^2$.
    For the rest of the choices of $\texttt{a}_S$,
    the figures are similar.}
    \label{fig:vel_cubo}
\end{figure}

Table \ref{tb:cub_stabilization_comparison} reveals that the method is capable to capture the expected spectrum reported in \cite{MR5023159}  for the chosen permeability. On this case, we also observe that the optimal order of convergence is attained for  each of the choices of $\texttt{a}_S$, where clearly the quadratic order of convergence is obtained when $a=8$, similarly  to the two dimensional case.  In Figure \ref{fig:vel_cubo} we present the first and fourth eigenvelocities for this configuration of the problem.

Now the aim is to consider a different permeability which escapes from a nearly Stokes regime. With this  aim in mind, we consider $\mathbf{K}^{-1}=10^3\mathbf{I}$ as permeability parameter and $a=8$ for the stabilization terms.
\begin{table}[!htbp]
\centering
\scriptsize
\setlength{\tabcolsep}{4pt}
\renewcommand{\arraystretch}{1.1}
\begin{tabular}{|c|c|c|c|c|c|c|}
\hline
$a_S$ & $N$ & $\lambda_{h,1}$ & $\lambda_{h,2}$
& $\lambda_{h,3}$ & $\lambda_{h,4}$ & $\lambda_{h,5}$\\
\hline

\multirow{7}{*}{$ak^2$}
& 9 & 71.4844 & 71.6208 & 71.6208 & 97.8144 & 97.8144\\
& 10 & 71.0276 & 71.1022 & 71.1022 & 97.0646 & 97.0646\\
& 11 & 70.5556 & 70.5782 & 70.5782 & 96.3024 & 96.3024\\
& 12 & 70.3297 & 70.3297 & 70.3371 & 95.9371 & 95.9371\\
\cline{2-7}
& Order & 1.82 & 1.97 & 2.06 & 1.90 & 1.90\\
\cline{2-7}
& Extrap. & 68.5869 & 68.5704 & 68.6841 & 93.2573 & 93.2573\\
\cline{2-7}
& \cite{MR5023159} 
& 69.1751 & 69.1685 & 69.1836 & 93.4348 & 93.4602\\
\hline

\multirow{7}{*}{$ak(k+1)$}
& 9 & 71.4844 & 71.6208 & 71.6208 & 97.8144 & 97.8144\\
& 10 & 71.0010 & 71.0713 & 71.0713 & 97.0274 & 97.0274\\
& 11 & 70.5556 & 70.5782 & 70.5782 & 96.3024 & 96.3024\\
& 12 & 70.3103 & 70.3103 & 70.3207 & 95.9118 & 95.9118\\
\cline{2-7}
& Order & 1.97 & 2.16 & 2.29 & 2.04 & 2.04\\
\cline{2-7}
& Extrap. & 68.7293 & 68.7476 & 68.8842 & 93.4609 & 93.4609\\
\cline{2-7}
& \cite{MR5023159} 
& 69.1751 & 69.1685 & 69.1836 & 93.4348 & 93.4602\\
\hline

\multirow{7}{*}{$a(k+1)^2$}
& 9 & 71.4346 & 71.5631 & 71.5631 & 97.7475 & 97.7475\\
& 10 & 70.9876 & 71.0559 & 71.0559 & 97.0089 & 97.0089\\
& 11 & 70.5268 & 70.5441 & 70.5441 & 96.2595 & 96.2595\\
& 12 & 70.3007 & 70.3007 & 70.3126 & 95.8992 & 95.8992\\
\cline{2-7}
& Order & 1.76 & 1.96 & 2.11 & 1.90 & 1.90\\
\cline{2-7}
& Extrap. & 68.5175 & 68.5693 & 68.7509 & 93.2615 & 93.2615\\
\cline{2-7}
&\cite{MR5023159} 
& 69.1751 & 69.1685 & 69.1836 & 93.4348 & 93.4602\\
\hline
\end{tabular}
\caption{First five eigenvalues for $k=1$,
$\mathbf{K}^{-1}=10^3\mathbf{I}$ and $a=8$,
considering three choices of the stabilization parameter
$\texttt{a}_S$.}
\label{tb:cub_k1_K3}
\end{table}

\begin{figure}[!htbp]
    \centering
    \begin{subfigure}{0.34\textwidth}
        \centering
        \includegraphics[width=\linewidth]{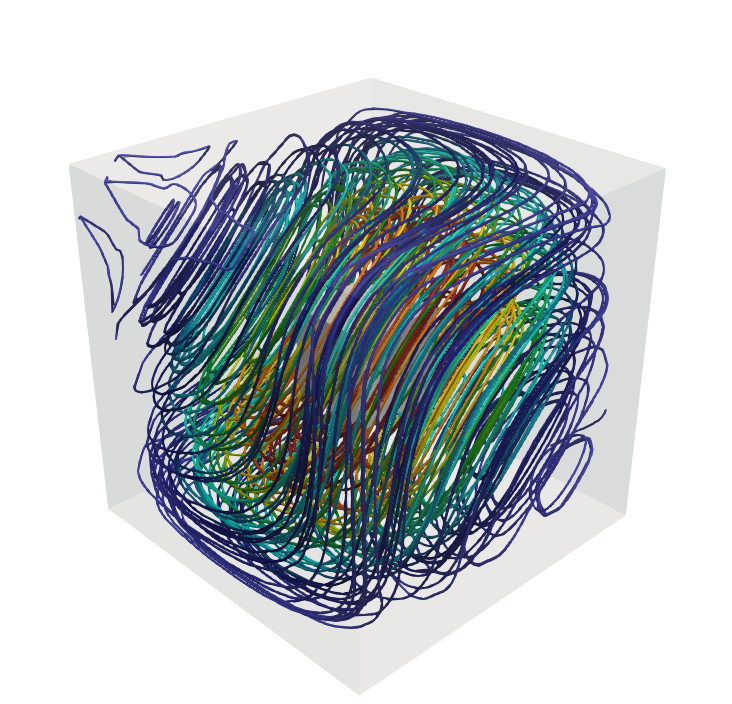}
    \end{subfigure}
    \hspace{0.06\textwidth}
    \begin{subfigure}{0.34\textwidth}
        \centering
        \includegraphics[width=\linewidth]{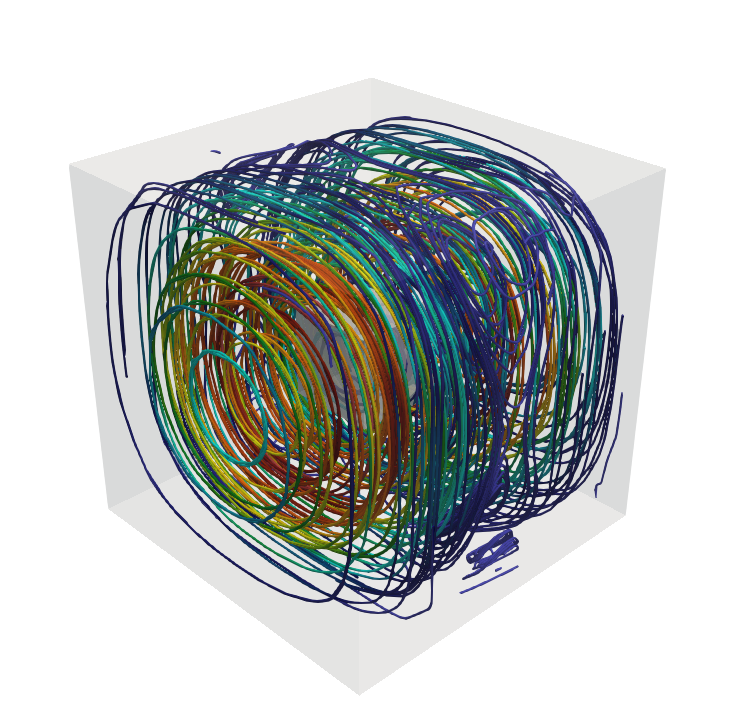}
    \end{subfigure}

    \caption{Representation of the velocity fields
    $\bu_{h,1}$ and $\bu_{h,4}$ for
    $\mathbf{K}^{-1}=10^{3}\mathbf{I}$,
    $k=1$ and ${\texttt{a}_S}=a(k+1)^2$.
    For the rest of the choices of $\texttt{a}_S$,
    the figures are similar.}

    \label{fig:vel_cubo_brinkman}
\end{figure}

The results reported in  Table \ref{tb:cub_k1_K3} show that for the lowest order of approximation $k=1$ the method is capable to obtain the desire order fo convergence when the domain presents a permeability far from the Stokes regime. We observe that for the first eigenvalue when we choose $\texttt{a}_S=a(k+1)^2$ the order of convergence is a bit poor compared with the other stabilizations but this can be improved with refined meshes. Also, the extrapolated values  of our method can be improved with more refined meshes in order to get more close to the reference values of \cite{MR5023159}. We end this test with plots of the velocity fields in Figure \ref{fig:vel_cubo_brinkman}.

\FloatBarrier

\textbf{Declaration of competing interest.}
The authors declare that they have no known competing financial interests or personal relationships that could have appeared
to influence the work reported in this paper.

\textbf{Acknowledgments.}
This paper is part of Isabella Estrada's Master thesis in the Master of  Applied Mathematics program  at Universidad del Bío-Bío, Chile.

\textbf{Data availability.}
Data will be made available on request.


\bibliographystyle{siamplain}

\end{document}